\documentclass[3p,preprint,12pt]{elsarticle}

\makeatletter
\def\ps@pprintTitle{%
  \let\@oddhead\@empty
  \let\@evenhead\@empty
  \let\@oddfoot\@empty
  \let\@evenfoot\@oddfoot
}
\makeatother
\usepackage{amsmath}
\usepackage{amssymb}
\usepackage{algorithmicx}
\usepackage{algpseudocode}
\usepackage{algorithm}

\usepackage{pgfplots}
\usepackage[toc,page]{appendix}
\usepackage{float}
\usepackage{booktabs}
\usepackage{bm}
\usepackage[parfill]{parskip}
\usepgfplotslibrary{groupplots}
\usepackage{booktabs}
\usepackage{makecell}

\usepackage{enumitem}
\usepackage{mathtools}
\usepackage{url}
\usepackage{amsthm}
\theoremstyle{definition}

\newcommand{\RR}[0]{\mathbb{R}}

\usepackage[utf8]{inputenc} 
\usepackage[T1]{fontenc}    
\usepackage{hyperref}       
\usepackage{cleveref}

\usepackage{booktabs}       
\usepackage{amsfonts}       
\usepackage{nicefrac}       
\usepackage{microtype}      
\usepackage{graphicx}
\usepackage{caption}
\usepackage{multirow}
\usepackage{wrapfig}
\usepackage{tabulary,xcolor}
\usepackage[normalem]{ulem}
\usepackage{pifont}
\usepackage[nointegrals]{wasysym}

\usetikzlibrary{patterns}
\definecolor{BurntOrange}{rgb}{0.8, 0.33, 0.0}
\usepackage{subcaption}

\usepackage{tikz}

\usepackage[utf8]{inputenc}
\usepgfplotslibrary{groupplots,dateplot}
\usetikzlibrary{patterns,shapes.arrows}
\pgfplotsset{compat=newest}
\usepackage{pgfplotstable}

\pgfplotstableset{format=inline}
            
\algdef{SE}[DOWHILE]{Do}{doWhile}{\algorithmicdo}[1]{\algorithmicwhile\ #1}%
\usepackage[utf8]{inputenc}

\begin{document}
\begin{frontmatter}

\title{Trust Region Method for Coupled Systems of PDE Solvers and Deep Neural Networks}

\author[label1]{Kailai Xu}
\ead{kailaix@stanford.edu}
\address[label1]{Institute for Computational and Mathematical Engineering, Stanford University, Stanford, CA, 94305}
 
\author[label1,label2]{Eric Darve}
\ead{darve@stanford.edu}
\address[label2]{Mechanical Engineering, Stanford University, Stanford, CA, 94305}

\begin{abstract}
Physics-informed machine learning and inverse modeling require the solution of ill-conditioned non-convex optimization problems. First-order methods, such as SGD and ADAM, and quasi-Newton methods, such as BFGS and L-BFGS, have been applied with some success to optimization problems involving deep neural networks in computational engineering inverse problems. However, empirical evidence shows that convergence and accuracy for these methods remain a challenge. Our study unveiled at least two intrinsic defects of these methods when applied to coupled systems of partial differential equations (PDEs) and deep neural networks (DNNs):  
(1) convergence is often slow with long plateaus that make it difficult to determine whether the method has converged or not; 
(2) quasi-Newton methods do not provide a sufficiently accurate approximation of the Hessian matrix; this typically leads to early termination (one of the stopping criteria of the optimizer is satisfied although the achieved error is far from minimal).
Based on these observations, we propose to use trust region methods for optimizing coupled systems of PDEs and DNNs. Specifically, we developed an algorithm for second-order physics constrained learning, an efficient technique to calculate Hessian matrices based on computational graphs. We show that trust region methods overcome many of the defects and exhibit remarkable fast convergence and superior accuracy compared to ADAM, BFGS, and L-BFGS. 
\end{abstract}

\begin{keyword}

\end{keyword}

\end{frontmatter}

\section{Introduction}

In machine learning applications, first order optimization methods, such as stochastic gradient descent \cite{ruder2016overview} and ADAM \cite{kingma2014adam}, have proven very successful. However, these methods usually suffer from slow convergence and weak performance (in terms of accuracy) when applied to physics-informed learning applications due to ill-conditioned optimization problems \cite{xu2020physics}. There is another wide class of optimization techniques that use curvature information, i.e., Hessians of loss functions. These optimization methods are  called \textit{second order optimization methods} \cite{agarwal2017second,luenberger1984linear,nocedal2006numerical}. Second order optimization methods, due to more exploitation of local information, lead to fast convergence for many problems. For example, Newton-Raphson methods \cite{ypma1995historical,ben1966newton} are an example of second order methods. If the loss function is convex in the neighborhood of the local minimum, and the Hessian matrix does not vanish, the Newton-Raphson method enjoys a convergence rate of order two, i.e., each iteration the scheme converges approximately to two significant digits. Gradient descent methods usually converge much slower and suffer for landscapes that contain many saddle points and flat areas. Thus, it is very promising and interesting to investigate these methods for solving inverse problems.

Particularly, we propose trust region methods \cite{conn2000trust,chen2018stochastic,kouri2018inexact,bui2014pde} for solving the optimization problem involving PDE solvers and DNNs \cite{huang2019predictive,xu2020learning,fan2020solving,xu2019neural}. Trust region methods possess some desirable features that are suitable for our problems: firstly, as we show in \Cref{sect:second-order-numerical-benchmarks}, Hessian matrices in our problems are semidefinite positive or indefinite. This fact rules out many second order optimization choices, such as Newton's method \cite{galantai2000theory}, which requires positive definiteness of Hessian matrices. Methods, such as BFGS \cite{dai2013perfect} and LBFGS \cite{liu1989limited,zhu1997algorithm,skajaa2010limited}, which maintain the positive definiteness of approximated Hessians,  are ultimately inconsistent with the characteristics of our problems. Secondly, line search-based optimizers \cite{more1994line} may converge fast but get stuck at a bad local minimum. These methods are too aggressive in searching for optimal steps in one direction per iteration. One consequence is that the activation values get saturated quickly and loss functions are no longer sensitive to changes of many weights and biases. The saturation reduces ``effective degrees of freedom'', defined in \Cref{sect:second-order-numerical-benchmarks}, and degrades the approximation capability of DNNs. 

\begin{figure}[htbp]
  \centering
  \includegraphics[width=0.6\textwidth]{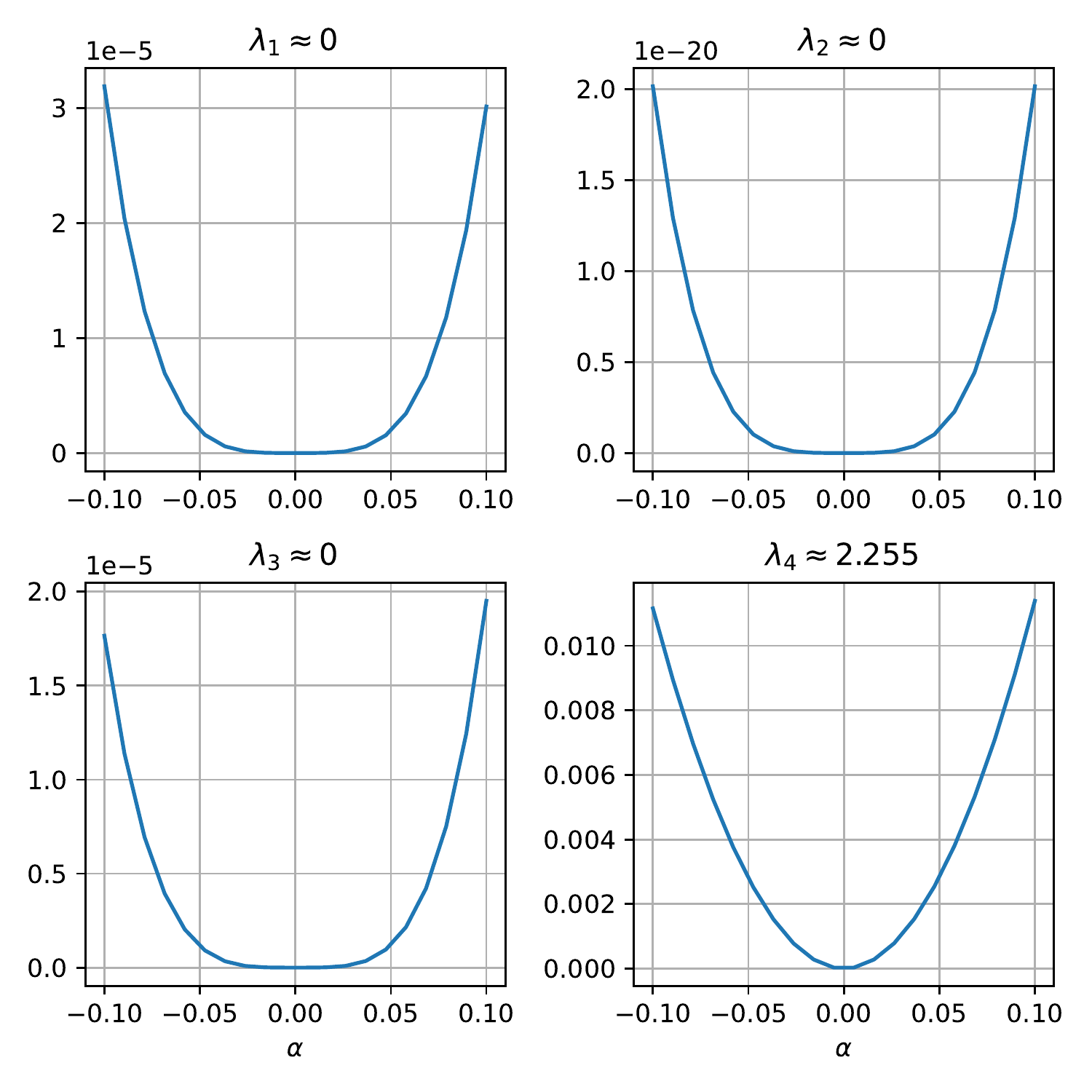}
  \caption{Perturbed loss function $\mathsf{L}(\alpha) = L(w^* + \alpha v_i)$ in the direction of eigenvectors $v_i$. $\lambda_i$ is the associated eigenvalue. Three eigenvalues are close to zero (with a scale of $10^{-13}$.}
  \label{fig:second-order-one-layer}
\end{figure}

Through the lens of Hessian eigenvalues, we reveal that over-parametrization of DNNs positions minimizers on a relatively higher dimensional manifold of the parameter space, and thus makes optimization easier. This can be illustrated through a simple example: consider a simple one layer neural network $y = w_2\tanh(w_1 x+b_1) + b_2$, and we fit the DNN with a pair $(x_0, y_0) = (0.5, \sin(0.5\pi))$. At convergence, we compute the Hessian of $L(w_1, b_1, w_2, b_2) = (w_2\tanh(w_1 x_0+b_1) + b_2 - y_0)^2$ and its eigenvalues/eigenvectors $(\lambda_i, v_i)$. We perturb the minimizer $\theta^*=(w_1^*, b_1^*,w_2^*, b_2^*)$ in the direction $v_i$ and compute the loss function $\mathsf{L}(\alpha) = L(\theta^* + \alpha v_i)$. From \Cref{fig:second-order-one-layer}, we can see that three eigenvalues are zero, and the loss function is flat in the corresponding eigenvector direction. That means that the minimizer almost lies on a 3-dimensional manifold in the 4-dimensional parameter space. We show in \Cref{sect:second-order-numerical-benchmarks} that this effects are prevailing in different problems, even when we couple a PDE solver with DNNs.

However, calculating Hessians, especially in a coupled system of PDEs and DNNs, is a very challenging task. At first glance, Hessians are just second-order partial derivatives of loss functions, and since automatic differentiation (AD) \cite{phipps2009sacado,revels2016forward,hogan2017adept,gower2010hessian,baydin2018automatic} enjoys great success in calculating gradients, we can apply AD twice to get Hessians. Indeed,
this approach has been implemented in many AD software. For example, in TensorFlow \cite{abadi2016tensorflow}, \texttt{tf.hessians} calculates Hessians by back-propagating the gradients from each component of the gradient. However, one fatal disadvantage of this approach is that it requires one gradient back-propagation per component of the gradient, which can be very expensive. In fact, for a function $f:\RR^d \rightarrow \RR$, when we calculate its gradients, TensorFlow constructs a computational graph for $g := \frac{\partial f}{\partial \theta}\in \RR^d$; here $x\in \RR$ is the input to $f$. due to the nature of reverse-mode automatic differentiation, it is in general impossible to reuse calculations for calculating the gradients of any two components in the Hessian
$$\frac{\partial g}{\partial \theta_i}, \quad \frac{\partial g}{\partial \theta_j}, \qquad 1\leq i < j \leq d$$
Thus, we need to construct at least $1 + 2 + \ldots + d = \frac{d(d+1)}{2}$ computational graphs for obtaining $\frac{\partial ^2 f}{\partial \theta_i \partial \theta_j}$, $\forall i, j$.
The second disadvantage of this approach is that we need to implement the adjoint update rule for each adjoint variable. For simple operators, such as addition, subtraction, $\sin$, $\cos$, etc., this task is not difficult. However, because in our problems there are many sophisticated operators related to numerical PDE solvers, adjoint update rules for adjoints further add unnecessary complications to implementations. 

Note that if what we need is Hessian vector dot product, we can do the calculation quite efficient using reverse-mode automatic differentiation. Let's consider a loss function $f(\theta)$, where $f: \RR^d \rightarrow \RR$, and we want to calculate $\nabla^2 f(\theta) p$, where $p\in \RR^d$, we can take gradient of $p^T\frac{\partial f}{\partial \theta}$ with respect to $\theta$; using TensorFlow, we have
\begin{verbatim}
u = tf.reduce_sum(tf.gradients(f, theta)[0] * p)
hessian_vector_dot_product = tf.gradients(u, theta)[0]
\end{verbatim}
However, the ability to calculate Hessians paves the way to more sophisticated and efficient optimization techniques, instead of restricting us to matrix-free approaches.

Forward-mode AD has also been used for calculating Hessians. For example, hyper-dual numbers \cite{fike2016derivative,fike2011development} are used for calculating second (or higher) order derivatives. In this approach, a new number system is introduced, with three units $1, \epsilon_1, \epsilon_2$, and a hyper-dual number $a$ have one real part and three non-real parts:
\begin{align*}
    a &= a_0 + a_1 \epsilon_1 + a_2 \epsilon_2 + a_3 \epsilon_1 \epsilon_2\\ 
    \epsilon_1^2 &=\epsilon_2^2 = 0 \\ 
    \epsilon_1 &\neq \epsilon_2 \neq 0\\ 
    \epsilon_1\epsilon_2 &= \epsilon_2 \epsilon_1 \neq 0
\end{align*}
Taylor series truncates exactly at second order derivative term 
\begin{equation}\label{equ:second-order-dual}
    f(\theta + h_1 \epsilon_1 + h_2\epsilon_2 + 0 \epsilon_1\epsilon_2) = f(\theta) + h_1 f'(\theta) \epsilon_1 + h_2 f'(\theta) \epsilon_2 + h_1 h_2 f''(x)\epsilon_1 \epsilon_2
\end{equation}
In implementations, each intermediate value is represented by a hyper dual number, and the ``forward-propagation'' rule is given by \Cref{equ:second-order-dual}. This formulation leads to no truncation error and no subtractive-cancellation error. However, the disadvantage is that we need to refurbish all existing operators and lift them to hyper dual number functions. 

In this chapter, we consider a more systematic approach that can achieve the same purpose but allows for simpler implementation and harnessing sparse structures of Hessian matrices. We derive an update rule for Hessian matrices  and leverage automatic differentiation for calculations. Our approach can be seen as an extension of the \texttt{edge\_pusing} algorithm proposed in \cite{gower2010hessian,gower2014computing,wang2016edge} to coupled systems of DNNs and PDEs in an AD computing environment. Note that in this paper, we focus on coupled system of PDE solvers and DNNs. There are many other frameworks for blending physical knowledge and deep neural networks, such as physics informed neural network \cite{raissi2019physics,meng2020ppinn,jagtap2020conservative,pang2020npinns} and deep learning approach to PDEs \cite{fan2019solving,fan2019solving}. These frameworks can also benefit from our approach and this line of research is left to further exploration.

\section{Trust Region Methods for Inverse Problems}

First order methods (e.g., ADAM) and quasi-Newton methods (e.g., BFGS and LBFGS) have been used for optimization in inverse problems. They both calculate a search direction and apply a step size adaptation or line search methods to generate a step size. For example, the ADAM optimizer calculates the search direction using past gradients and past square gradients and then applies adaptive learning rates for each parameter. BFGS uses secant equations to generate search directions and performs a line search to identify the optimal next step. However, these methods fail to exploit or only partially exploit curvature information of optimization problems. As we have pointed out, second order methods use local Hessian information to find a new minimizer and thus improve convergence and accuracy. 

Trust region methods are a class of second order optimization methods, which define a region around the current iterate within which they focus on finding a suitable next minimizer. At each iterate $f(x_k)$, the objective function $f(x_k+p)$ is approximated by a model quadratic function $m_k(p)$
$$m_k(p) = f_k + g_k^T p + \frac{1}{2}p^T B_kp$$
where $f_k = f(x_k)$, $g_k = \nabla f(x_k)$, and $B_k = \nabla^2 f(x_k+tp)$ is the Hessian matrix at $x_k + tp$. In some applications, $\nabla^2 f(x_k+tp)$ is approximated by a symmetric matrix $B_k$; for example, $B_k$ can be calculated using secant equations as in BFGS. However, despite the challenges of calculating Hessians, we do not use approximations in this thesis for two reasons: firstly, when we use the exact Hessian, the approximation error in the model $m_k$ is $\mathcal{O}(\|p\|^3)$ while with an approximation $B_k$, the error is $\mathcal{O}(\|p\|^2)$; secondly, the exactly Hessian matrices for our problems are usually indefinite or semi-positive definite, while some commonly used approximations, such as in BFGS, the resultant approximations are positive definite. Both lead to inaccurate local approximations, in which case trust region methods are ultimately self-defeating. Additionally, the indefiniteness/semi-positive definiteness also rules out some other second order methods, such as Newton's methods, which require that Hessian matrices are positive definite.

\begin{figure}[htbp]
  \centering
  \includegraphics[width=0.6\textwidth]{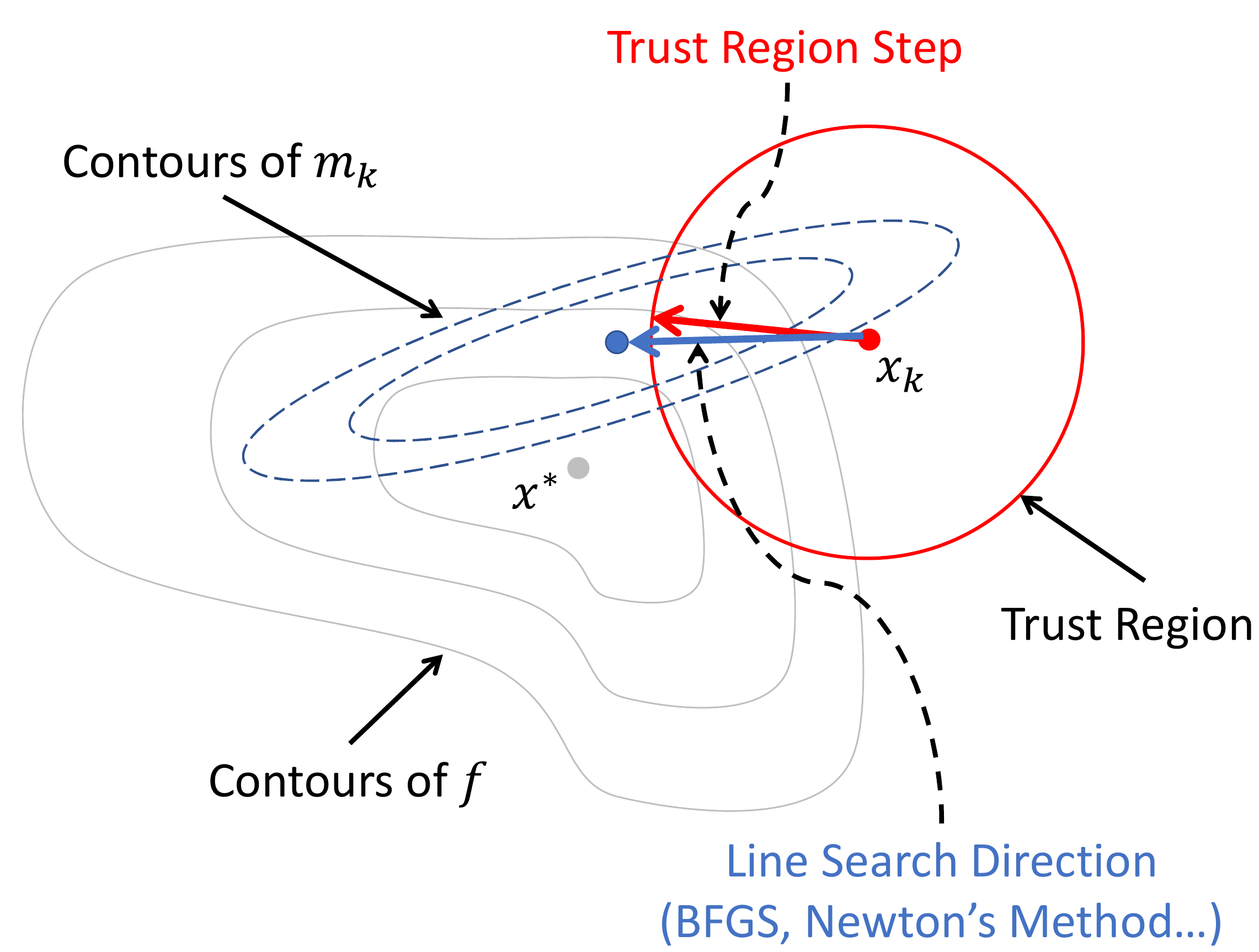}
  \caption{An illustration of trust region methods. $f$ is approximated by $m_k$ at $x_k$. Trust region methods find a new minimizer within the trust region (the region within the red circle).}
  \label{fig:second-order-trust-region}
\end{figure}

Upon local approximations, trust region methods seek a solution of the subproblem (see \Cref{fig:second-order-trust-region})
\begin{equation}\label{equ:second-order-local-problem}
    \min_{p} \ m_k(p) = f_k + g_k^T p + \frac{1}{2}p^T B_k p\quad \text{s.t.}\ \|p\| \leq \Delta_k
\end{equation}
where $\Delta_k$ is the trust region radius, which is adaptive as optimization proceeds. There exist many strategies for solving \Cref{equ:second-order-local-problem}, such as the dogleg method and two-dimensional subspace minimization method. However, these methods typically require positive definite Hessian matrices, which are not the case of our problems. Here, we use the nearly exact trust region method proposed in \cite{conn2000trust} (Chapter 7). This approach is implemented in the \texttt{scipy} library \cite{virtanen2020scipy}, which we use to solve \Cref{equ:second-order-local-problem}. 

\section{Hessian Calculation}

Consider a twice differentiable function, whose calculation can be expressed in the form (the notation is a slightly modified version from \cite{grievank2000principles})
\begin{equation}\label{equ:second-order-vm}
    v_m = P(\Phi_m(\Phi_{m-1}(\cdots(\Phi_1(v)))))
\end{equation}
Here $v_m$ is a scalar variable, $v$ is a vector, and $\Phi_i$ are intermediate functions, mapping vectors to vectors. $P$ is a linear operator, which returns the last component of the output of $\Phi_m$. In the following, we describe an approach \cite{gower2014computing,wang2016edge} to calculate the Hessian while reusing calculations from reverse-mode automatic differentiation. 

Let us consider an intermediate step $\Phi_i$, we define
$$f:= P(\Phi_m (\Phi_{m-1} (\cdots (\Phi_{i+1}(v))))),\qquad x:=\Phi_{i-1}(\cdots(\Phi_1(v)))$$
For simplicity, we also drop the subscript in $\Phi_i$ and write $G:=\Phi_i$ \Cref{equ:second-order-vm} can be rewritten as 
\begin{equation}\label{equ:algo-second-order-fG}
    v_m = f(G(x))
\end{equation}
We define 
\begin{align*}
    f_{,k}(y) &= \frac{\partial f(y)}{\partial y_k}, \quad f_{,kl}(y) = \frac{\partial^2 f(y)}{\partial y_k \partial y_l} \\ 
    G_{k,l}(x) &=\frac{\partial G_k(x)}{\partial x_l},\quad G_{k, lr}(x) = \frac{\partial^2 G_k(x)}{\partial x_l\partial x_r}
\end{align*}
Here $y_k$ (or $x_k$) is the $k$-th component of $y$ (or $x$) and $G_k$ is the $k$-th component of $G$. 

\begin{figure}[htbp]
  \centering
  \includegraphics[width=0.6\textwidth]{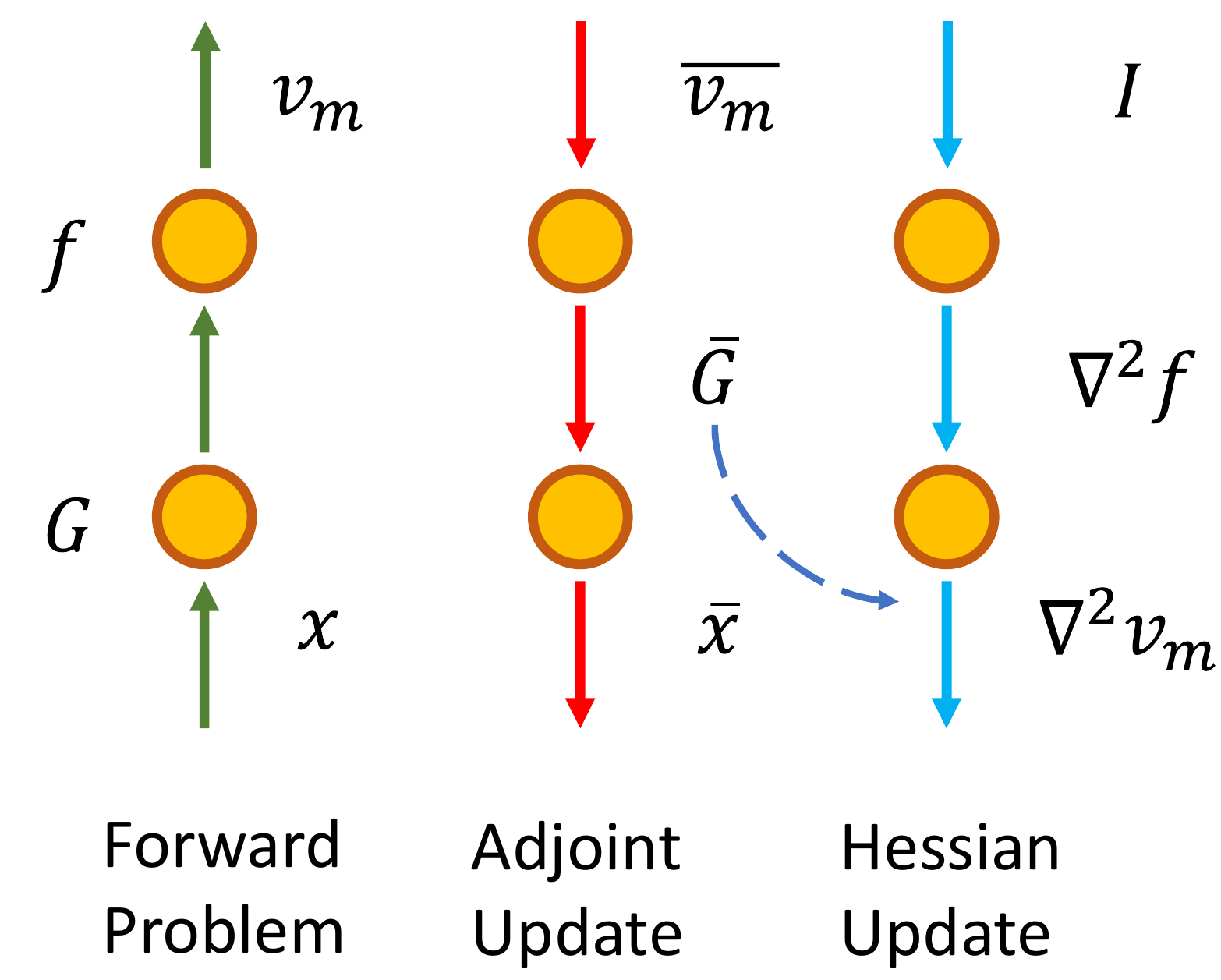}
  \caption{Second order physics constrained learning for updating Hessian matrices \Cref{equ:second-order-update}. Adjoint updates and Hessian updates are both in the reverse order of forward computation. When calculating $\nabla^2 v_m$, the adjoint variable $\bar G$ is used.}
  \label{fig:second-order-hessian}
\end{figure}

Then we have
$$\frac{\partial v_m}{\partial x_i} = f_{,k}G_{k,i}$$
Here we used the Einstein notation: the same index indicates a tensor contraction on the index. We take the derivative with respect to $x_j$ on both sides and get 
$$\frac{\partial^2 v_m}{\partial x_i \partial x_j} = f_{,kr} G_{k,i} G_{r,j} + f_{,k} G_{k,ij}$$
This equation can be rewritten in the vector/matrix form (\Cref{fig:second-order-hessian})
\begin{equation}\label{equ:second-order-update}
    \boxed{\nabla^2 v_m = (\nabla G)^T \nabla^2 f (\nabla G) + \nabla^2 (\bar G^T G)}
\end{equation}
Here $\bar G:=\nabla f$, which is the adjoint variable of $G(x)$ and in calculating $\nabla^2 (\bar G^T G)$, $\bar G$ is assumed to be independent of $x$. The gradient (Jacobian) $\nabla$ and Hessian $\nabla^2$ should be considered taken with respect to $x$. \Cref{equ:second-order-update} is called the \textit{Hessian update rule}.

Several remarks are in order: firstly, the adjoint variable $\bar G$ is already available if we have performed gradient back-propagation. Thus, we can reuse this value instead of recomputing it. Secondly, $G$ only depends on $v_k, k\prec i-1$, and $F$ only depends on $v_k = G(x)_k, k\prec i$. Therefore, many entries in $\nabla^2 F$ and $\nabla G$ are zero. We can ``condense'' \Cref{equ:second-order-update} with 
\begin{equation}\label{equ:second-order-update2}
\nabla^2 v_m = (\nabla G_{\mathcal{I}})^T \nabla^2 F_{\mathcal{J}} (\nabla G_{\mathcal{I}}) + \nabla^2 (\overline{ G_{\mathcal{I}}}^T G_{\mathcal{I}})\qquad \mathcal{I} = \{k: k \prec i-1\},  \mathcal{J} = \{k: k \prec i\}
\end{equation}
Here the subscript indicates only keeping the components of this set of $v_k$'s. The gradient (Jacobian) $\nabla$ and Hessian $\nabla^2$ should be considered taken with respect to $x_{\mathcal{I}}:=\{x_k: k\in \mathcal{I}\}$. Using \Cref{equ:second-order-update2} saves much space and computation than \Cref{equ:second-order-update}. Lastly, in general $\nabla^2 v_m$ is a dense matrix due to the existence of DNNs (as well as PDEs in many cases). But we can still exploit the sparsity when calculating $\nabla G$ or $\nabla G_{\mathcal{I}}$. For example, if $G(x) = 2x$, then $\nabla G(x) = 2I$, where $I$ is the identity matrix, which is sparse. 

We summarize the Hessian update algorithm in \Cref{algo:second-order-algo}. This algorithm only requires one backward pass and constructs the Hessian iteratively. Additionally, we can leverage the symmetry of the Hessian when we do the calculations in \Cref{algo:second-order-algo-update-line}. This algorithm also doesn't require looping over each component of the gradient. 

However, the challenge here is that we need to calculate $\nabla G$ and $\nabla^2(\bar G^T G)$. Developing complete support of such calculations for all operators can be a laborious but rewarding task. 

\begin{algorithm}
\caption{Second order physics constrained learning based on edge pushing method}
\label{algo:second-order-algo}
\begin{algorithmic}[1]
  \State Initialize $H \gets 0$
  \For{$k=m-1, m-2, \ldots, 1$}
  \State Define $f:=  P(\Phi_m (\Phi_{m-1} (\cdots (\Phi_{k+1}(\cdot)))))$, $G:= \Phi_k$
  \State Calculate the gradient (Jacobian) $J \gets \nabla G$
  \State Extract $\bar G$ from the saved gradient back-propagation data. 
  \State Calculate $Z = \nabla^2(\bar G^T G)$ \label{algo:second-order-algo-update-line}
  \State Update $H \gets J^THJ + Z$
  \EndFor
\end{algorithmic}
\end{algorithm}

\section{Developing Second Order PCL for a Sparse Linear Solver}\label{sect:sparse_linear_solver}

Here we consider an application of second order PCL for a sparse solver. We focus on the operator that takes the sparse entries of a matrix $A\in\mathbb{R}^{n\times n}$ as input and outputs $u$
$$Au = f$$
Here $f$ is a given vector. Denote $A = [a_{ij}]$, where some of $a_{ij}$ are zero. According to second order PCL, we need to calculate $\frac{\partial u_k}{\partial a_{ij}}$ and $\frac{\partial^2 (y^T u)}{\partial a_{ij} \partial a_{rs}}$ for an adjoint variable $y\in \RR^n$. This problem is related to \Cref{equ:algo-second-order-fG} if we let
$$x = \{a_{ij}\}, \ G(x) = A^{-1}f$$

We consider two multi-indices $l$ and $r$, which are tuples of two integers (e.g., $l = (1,2)$). We take the gradient with respect to $a_l$ on both sides of
$$a_{i1}u_1 + a_{i2}u_2 + \ldots + a_{in}u_n = f_i$$
which leads to
\begin{equation}\label{equ:second-order-aij}
    a_{i1,l}u_1 + a_{i2,l}u_2 + \ldots + a_{in,l}u_n + a_{i1}u_{1,l} + a_{i2}u_{2,l} + \ldots + a_{in}u_{n,l} = 0
\end{equation}
Here the second index in the subscript indicates the derivative. \Cref{equ:second-order-aij} leads to
\begin{equation}\label{equ:second-order-ul}
u_{,l} = -A^{-1}A_{,l} u    
\end{equation}
Note at most one entry in $A_{,l}$ is nonzero, and therefore at most one entry in $A_{,l} u$ is nonzero. Thus to calculate \Cref{equ:second-order-ul}, we can calculate the inverse $A^{-1}$ first, and then $u_{,l}$ can be obtained cheaply by taking a column from $A^{-1}$. The complexity will be $\mathcal{O}(n^3)$---the cost of inverting $A$.

Now take the derivative with respect to $a_r$ on both sides of \Cref{equ:second-order-aij}, we have
$$\begin{aligned} a_{l,i1}u_{1,r} + a_{l,i2}u_{2,r} + \ldots + &\\ a_{l,in}u_{n,r} + \ a_{r,i1}u_{1,l} + a_{r,i2}u_{2,l} + \ldots + &\\ a_{r,in}u_{n,l} +\ a_{i1}u_{1,rl} + a_{i2}u_{2,rl} + \ldots + a_{in}u_{n,rl} = 0 &\end{aligned}$$
which leads to
$$Au_{,rl} = -A_{,l} u_{,r} - A_{,r} u_{,l}$$
Therefore,
$$(y^Tu)_{,rl} = - y^TA^{-1}(A_{,l} u_{,r} + A_{,r} u_{,l})$$
We can calculate $z^T = y^TA^{-1}$ first with a cost $\mathcal{O}(n^2)$. Because $u_{,r}$, $u_{,l}$ has already been calculated and $A_{,l}$, $A_{,r}$ has at most one nonzero entry, $A_{,l} u_{,r} + A_{,r} u_{,l}$ has at most two nonzero entries. The calculation $z^T(A_{,l} u_{,r} + A_{,r} u_{,l})$ can be done in $\mathcal{O}(1)$ and therefore the total cost is $\mathcal{O}(d^2)$, where $d$ is the number of nonzero entries.

Upon obtaining $\frac{\partial u_k}{\partial a_{ij}}$ and $\frac{\partial^2 (y^T u)}{\partial a_{ij} \partial a_{rs}}$, which corresponds to $\nabla G$ and $\nabla^2(\bar G^T G)$ in \Cref{equ:second-order-update}, respectively, we can apply the recursive the formula to ``back-propagate'' the Hessian matrix.

\section{Numerical Benchmarks}\label{sect:second-order-numerical-benchmarks}

The numerical benchmarks are performed in the dimensional space. We consider Poisson's equations and heat equations. Particularly, we consider two types of Poisson's equation: the diffusivity coefficient is spatially-varying or dependent on the state variable. The computational domain is $[0,1]\times [0,1]$. We consider both finite difference and finite element methods. The finite element case is more technically challenging than the finite-difference case because it involves matrix/vector assembly, for which we need to apply the Hessian update rule \Cref{equ:second-order-update2}. In what follows, the deep neural network is a fully-connected deep neural network, with 3 hidden layers, each layer has 20 neurons, and the activation functions are $\tanh$. The numerical PDE solver and optimizers are implemented using the ADCME and AdFem libraries \cite{xu2020adcme}, which are available from 
\begin{center}
    \url{https://github.com/kailaix/ADCME.jl}\\ \url{https://github.com/kailaix/AdFem.jl}
\end{center}

\subsection{Static Poisson's Equation: Residual Minimization}\label{sect:second-order-static}

In this example, we consider the Poisson's equation
\begin{equation}\label{equ:second-order-ex1}
\begin{aligned}
     \nabla \cdot (\kappa(u) \nabla u)) &= f(x) & x\in \Omega\\ 
     u &= 0 & x\in  \partial\Omega
\end{aligned}
\end{equation}
Here the diffusivity coefficient $\kappa$ is given by 
$$\kappa(u) = 2.0 - (1.4-3u)\sin(18u)$$
and the exact solution $u$ is given by
$$u(x,y) = x (1-x) (1-y)^2\sin(y)$$

We discretize \Cref{equ:second-order-ex1} on a uniform grid with a grid size $h = 0.1$. Assume we can observe the values of $u$ and $f$ on the grid points, denoted by $u_{i,j}$ and $f_{i,j}$. We want to estimate $\kappa$ from these observations. To make the problem more challenging, we add $10\%$ uniform random noise to the observations, i.e., 
\begin{equation}
    \begin{aligned}
    \hat u_{i,j} &= u_{i,j} (1+0.1z_{i,j})\\ 
    \hat f_{i,j} &= f_{i,j} (1+0.1w_{i,j})\\ 
    \end{aligned}
\end{equation}
where $z_{i,j}\sim_{\text{i.i.d}}\mathcal{U}(0,1)$, $w_{i,j}\sim_{\text{i.i.d}}\mathcal{U}(0,1)$. 

We use a deep neural network $\kappa_\theta(u)$, where $\theta$ is the weights and biases of the DNN. 
We formulate the inverse problem as an optimization problem, which minimizes the total residuals
\begin{equation}\label{equ:second-order-ex1-optimization}
    \min_\theta\; \sum_{i,j} (F_{i,j}(\hat u; \theta) - \hat f_{i,j})^2 
\end{equation}
where ($h$ is the step size)
\begin{equation*}
\begin{aligned}
    F_{i,j}(\hat u; \theta) &= \frac{1}{h^2}\Bigg( \frac{\kappa_\theta(\hat u_{i+1,j}) + \kappa_\theta(\hat u_{i,j})}{2} (u_{i+1,j} - u_{i,j}) \\
    & - \frac{\kappa_\theta(u_{i,j}) + \kappa_\theta(u_{i-1,j})}{2} (\hat u_{i,j} - \hat u_{i-1,j}) \\ 
    &+  \frac{\kappa_\theta(\hat u_{i,j+1}) + \kappa_\theta(\hat u_{i,j})}{2} (\hat u_{i,j+1} - \hat u_{i,j})\\
    &-  \frac{\kappa_\theta(\hat u_{i,j}) + \kappa_\theta(\hat u_{i,j-1})}{2} (\hat u_{i,j} - \hat u_{i,j-1})\Bigg) 
\end{aligned}
\end{equation*}

\begin{figure}[htbp]
  \centering
  \includegraphics[width=0.45\textwidth]{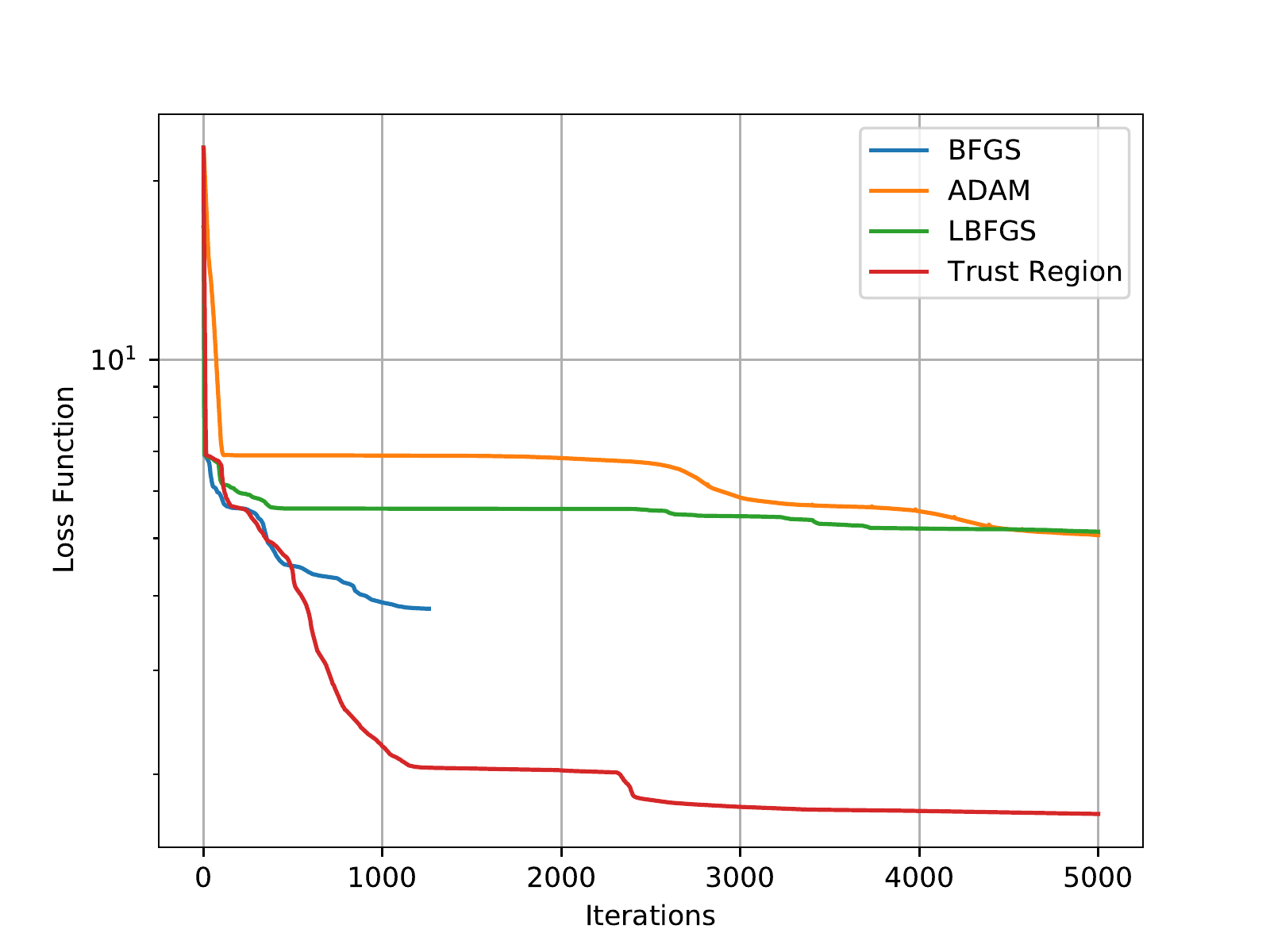}~
  \includegraphics[width=0.45\textwidth]{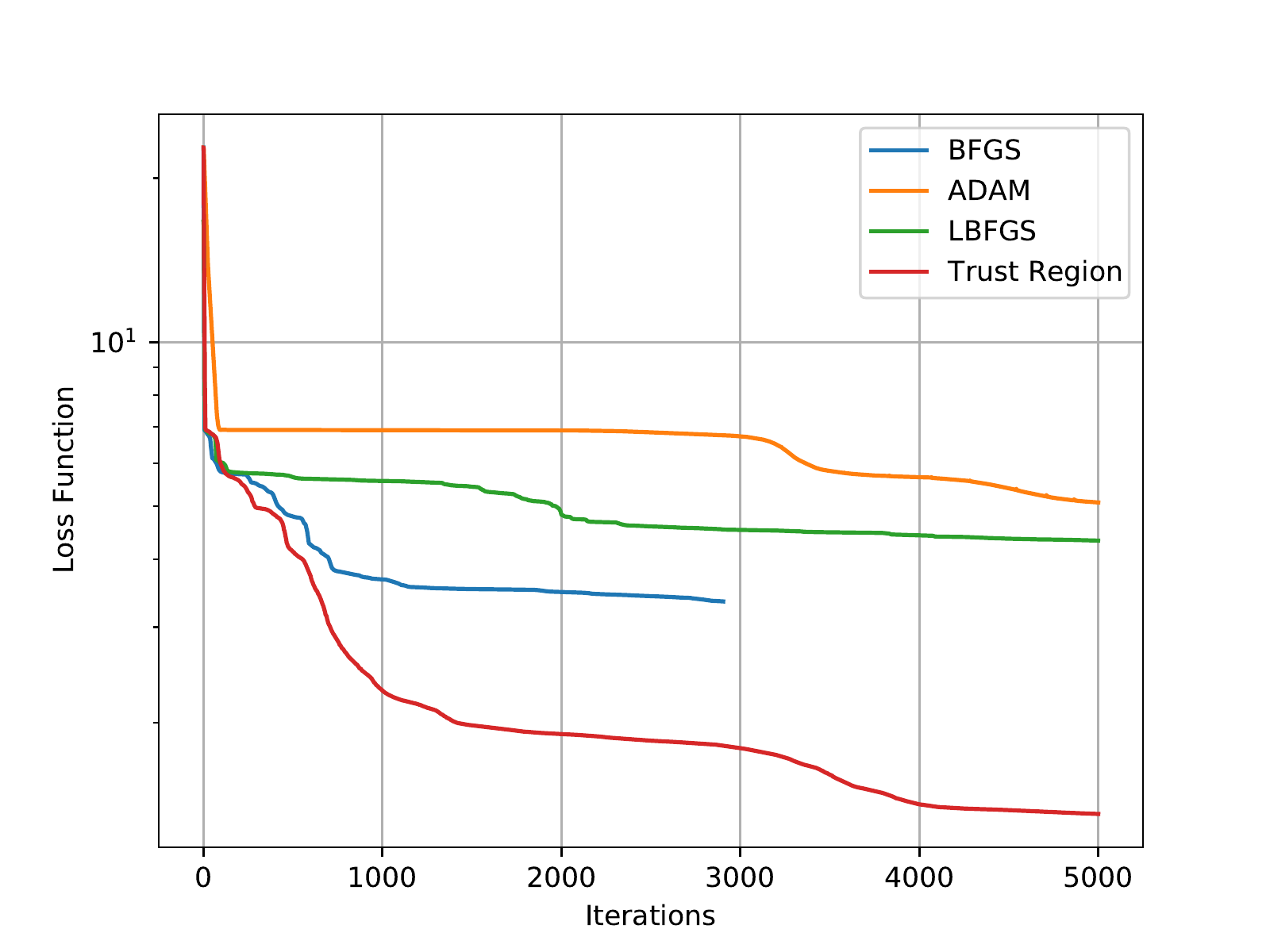}
   \includegraphics[width=0.45\textwidth]{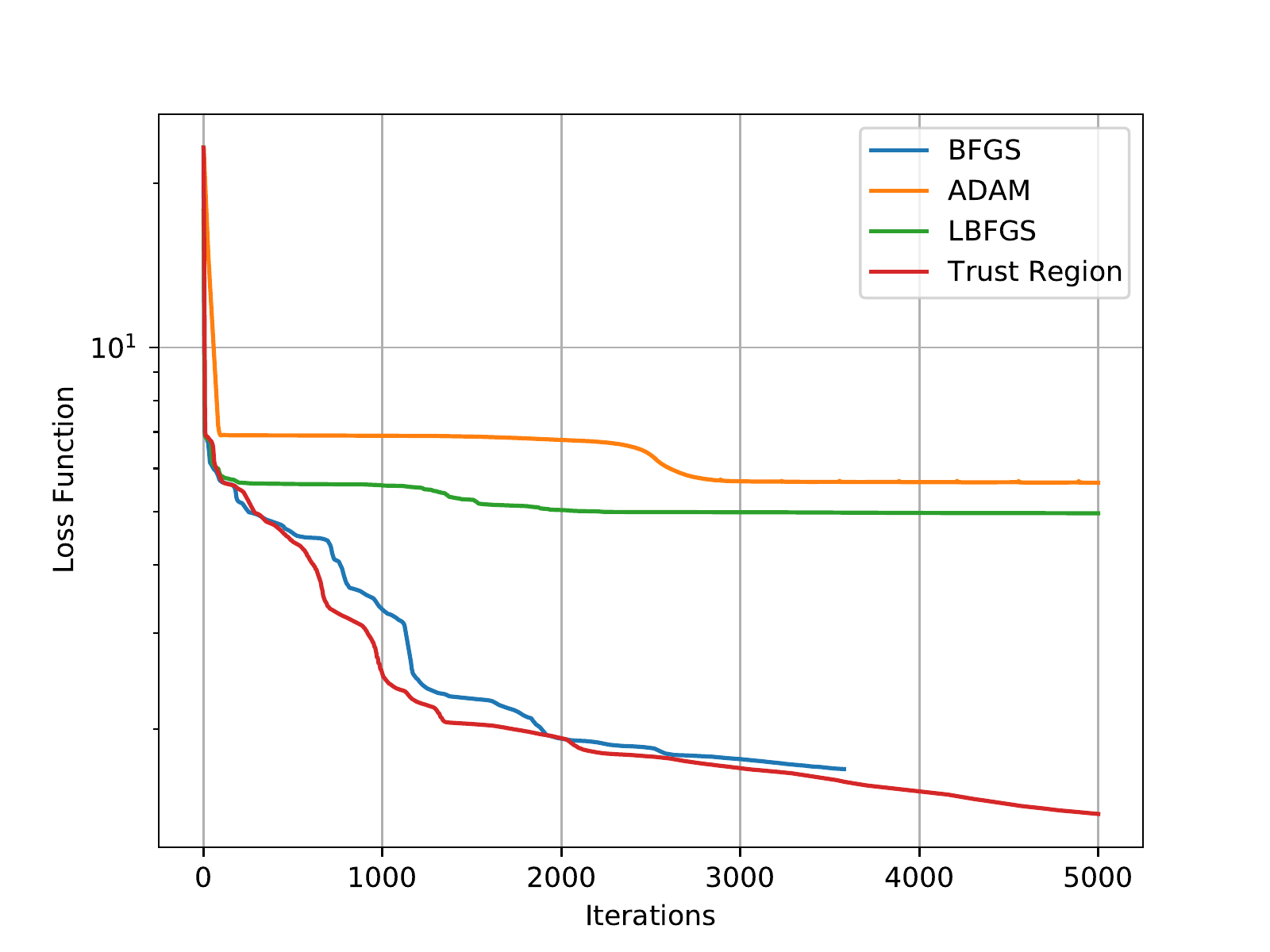}~
  \includegraphics[width=0.45\textwidth]{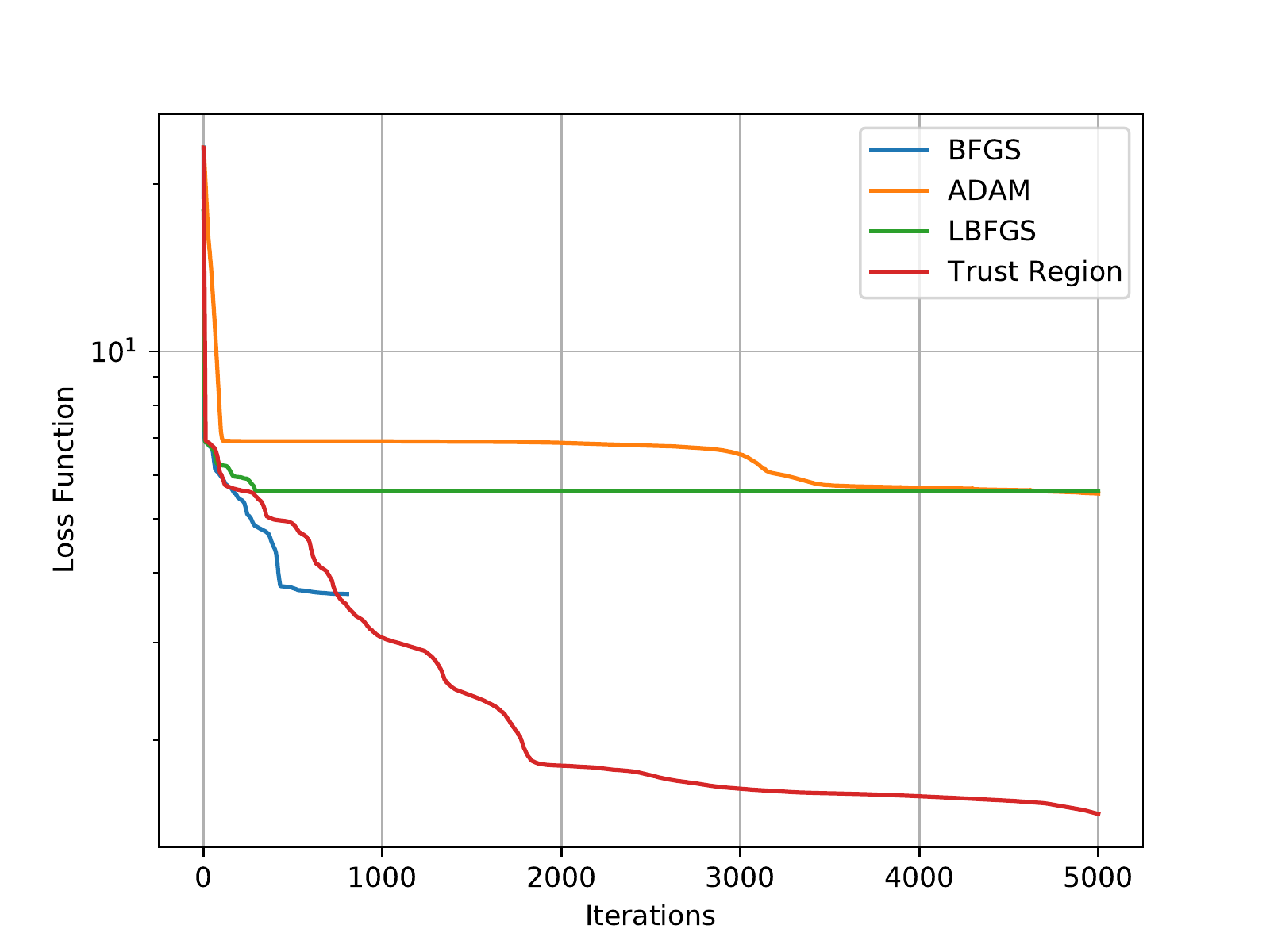}
  \caption{Loss functions for solving \Cref{equ:second-order-ex1-optimization} using different optimizers. Trust region methods perform consistently better than all other three methods.}
  \label{fig:second-order-ex1-loss}
\end{figure}

We apply four optimizers to solve \Cref{equ:second-order-ex1-optimization}. We use the default hyperparameters for optimizers in the ADCME library. The maximum number of iterations is capped at 5000.  Because the optimization results depend on the initialization of the deep neural network, we use four different initial guesses for DNNs. \Cref{fig:second-order-ex1-loss} shows that trust region methods perform consistently better than all other three methods. Despite its fast speed per iteration, ADAM converges much slower and may get stuck at a bad local minimum and never make any progress further. For this reason, we only discuss BFGS, LBFGS, and trust region methods in the following. 

\begin{figure}[htbp]
  \centering
  \includegraphics[width=0.45\textwidth]{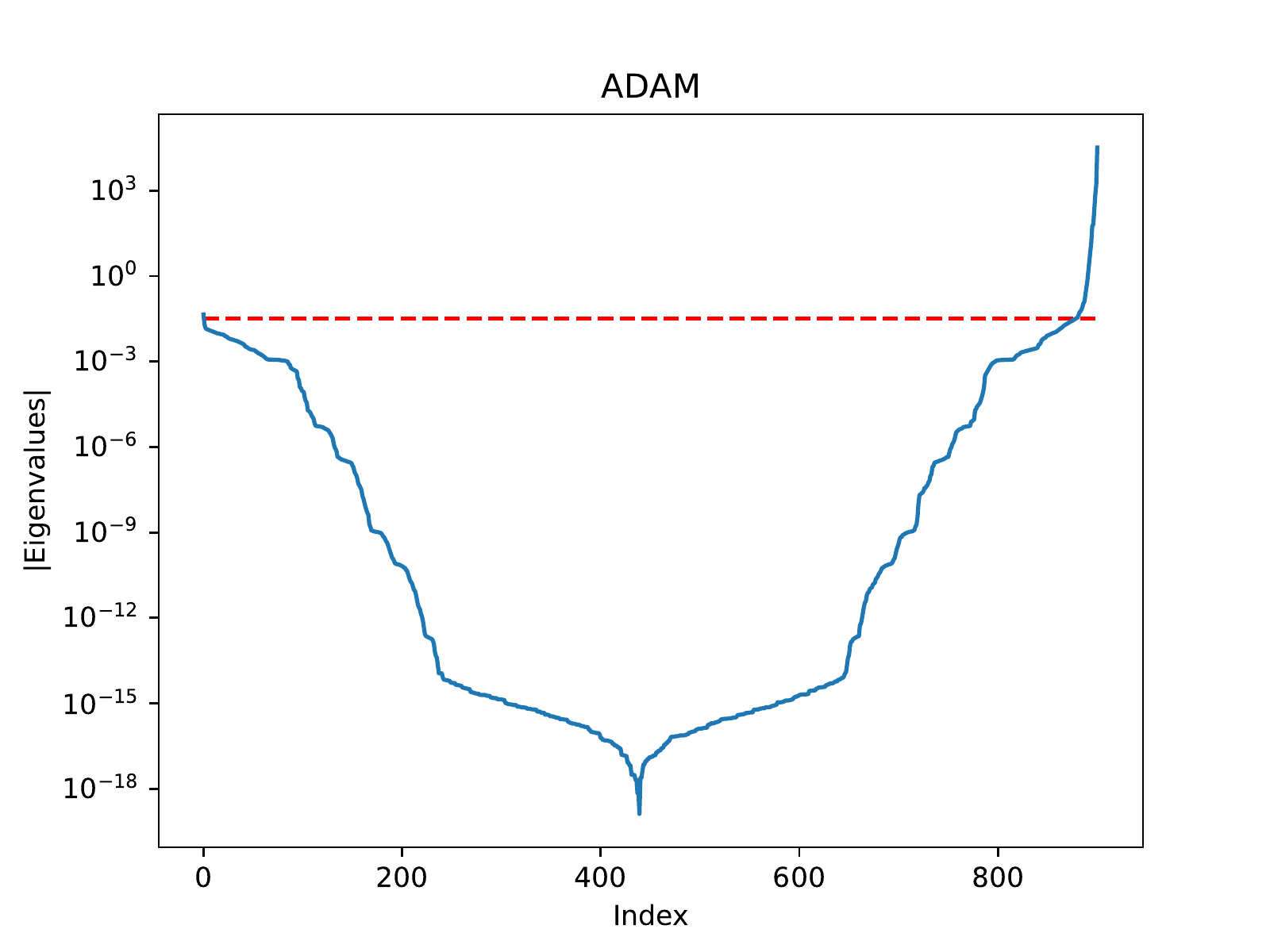}~
  \includegraphics[width=0.45\textwidth]{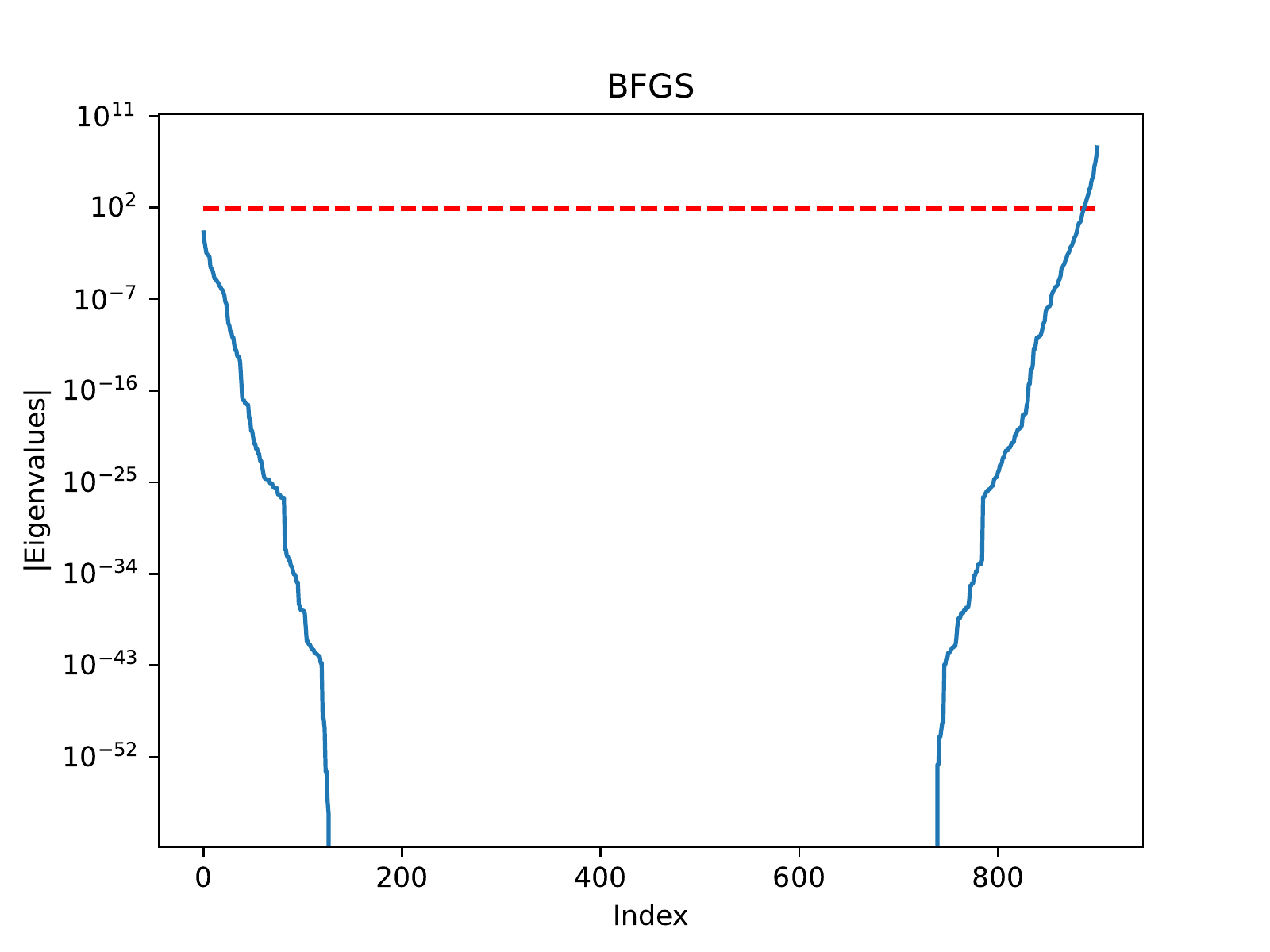}
   \includegraphics[width=0.45\textwidth]{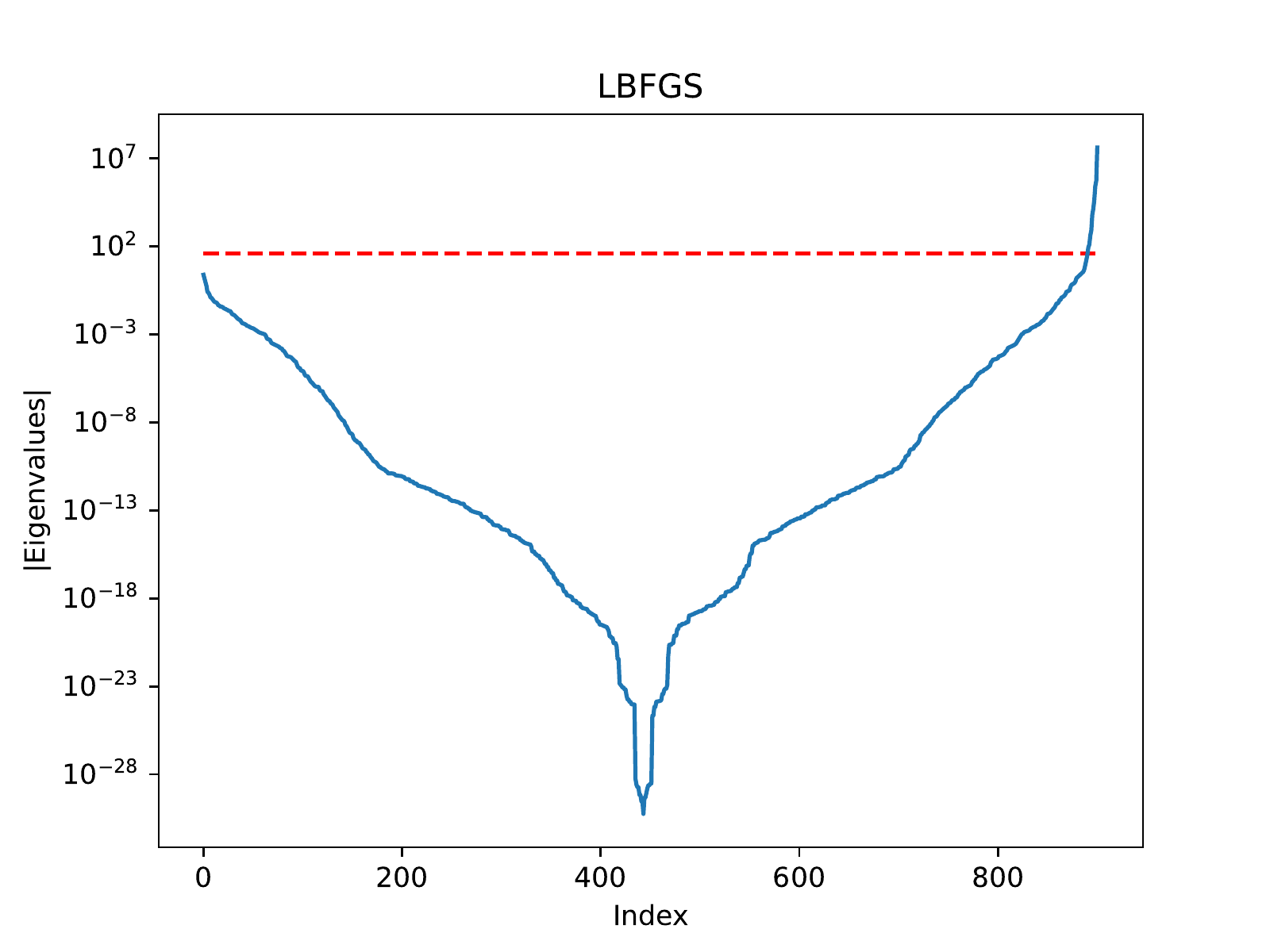}~
  \includegraphics[width=0.45\textwidth]{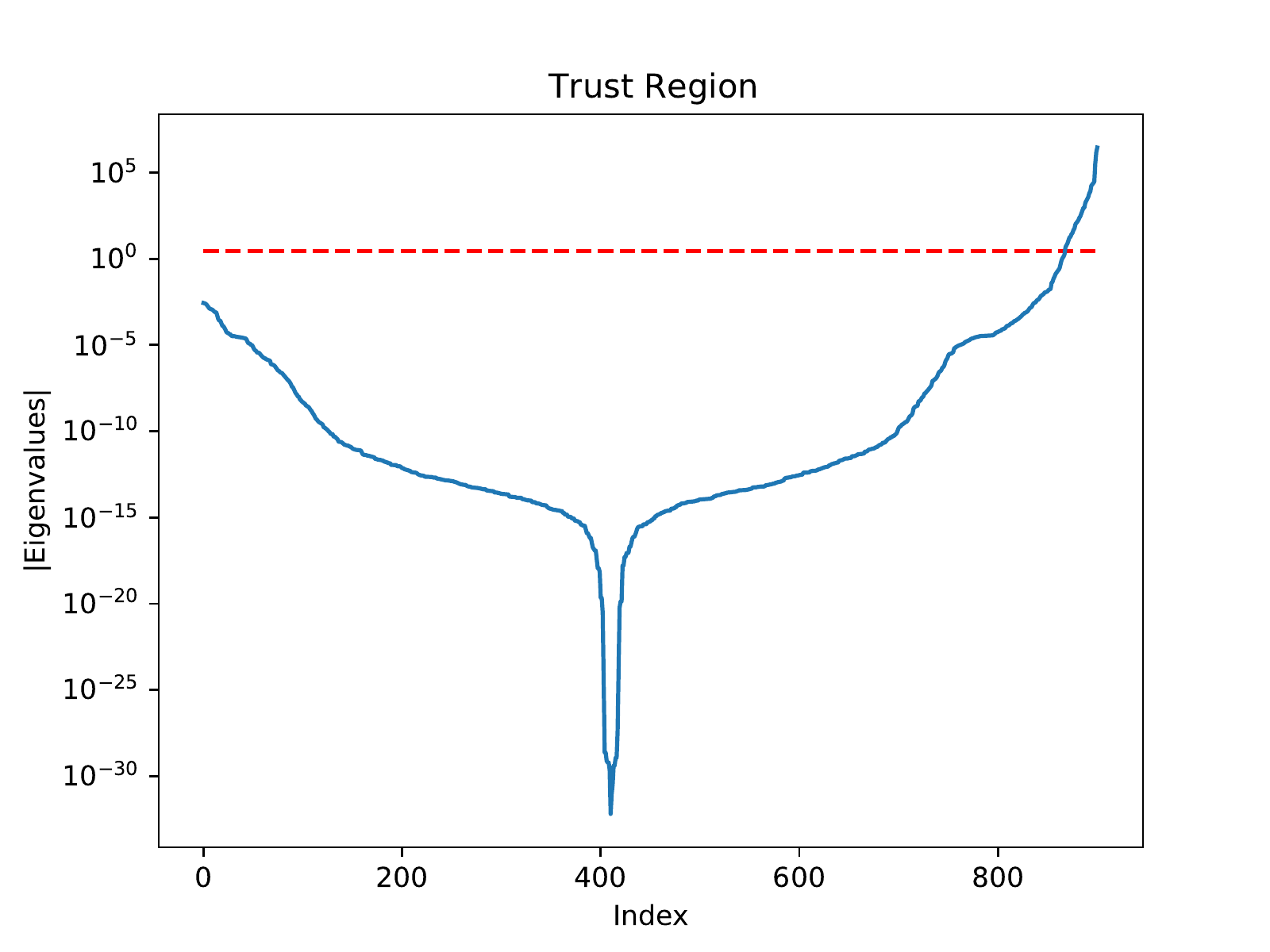}
  \caption{Magnitudes of eigenvalues at the last iteration for \Cref{equ:second-order-ex1-optimization}. The left decreasing values correspond to negative eigenvalues, and the right increasing ones correspond to positive eigenvalues. The red dashed line correspond to $10^{-6} \lambda_{\max}$, where $\lambda_{\max}$ is the maximum eigenvalue.}
  \label{fig:second-order-ex1-eigenvalues}
\end{figure}

\Cref{fig:second-order-ex1-eigenvalues} shows the magnitude of eigenvalues at the last iteration for the second case in \Cref{fig:second-order-ex1-loss}. Here we use a threshold $\epsilon=10^{-6}$ for classifying the sign of eigenvalues: for a given eigenvalue $\lambda$, it is treated as "positive" if $\lambda>\epsilon \lambda_{\max}$, and ``negative'' if $\lambda < - \epsilon \lambda_{\max}$, otherwise zero. Here $\lambda_{\max}$ is the maximum eigenvalue. We see that the Hessian matrices are all semi-positive definite, indicating that the stationary points are all local minimum. Also, we can see only a small fraction of eigenvalues are positive. We interpret the associated eigenvectors as ``effective degrees of freedoms (DOFs)'' because only perturbation in those directions changes the loss function values. We will revisit effective DOFs in \Cref{sect:second-order-dynamic}. 

\begin{figure}[htbp]
  \centering
  \includegraphics[width=0.45\textwidth]{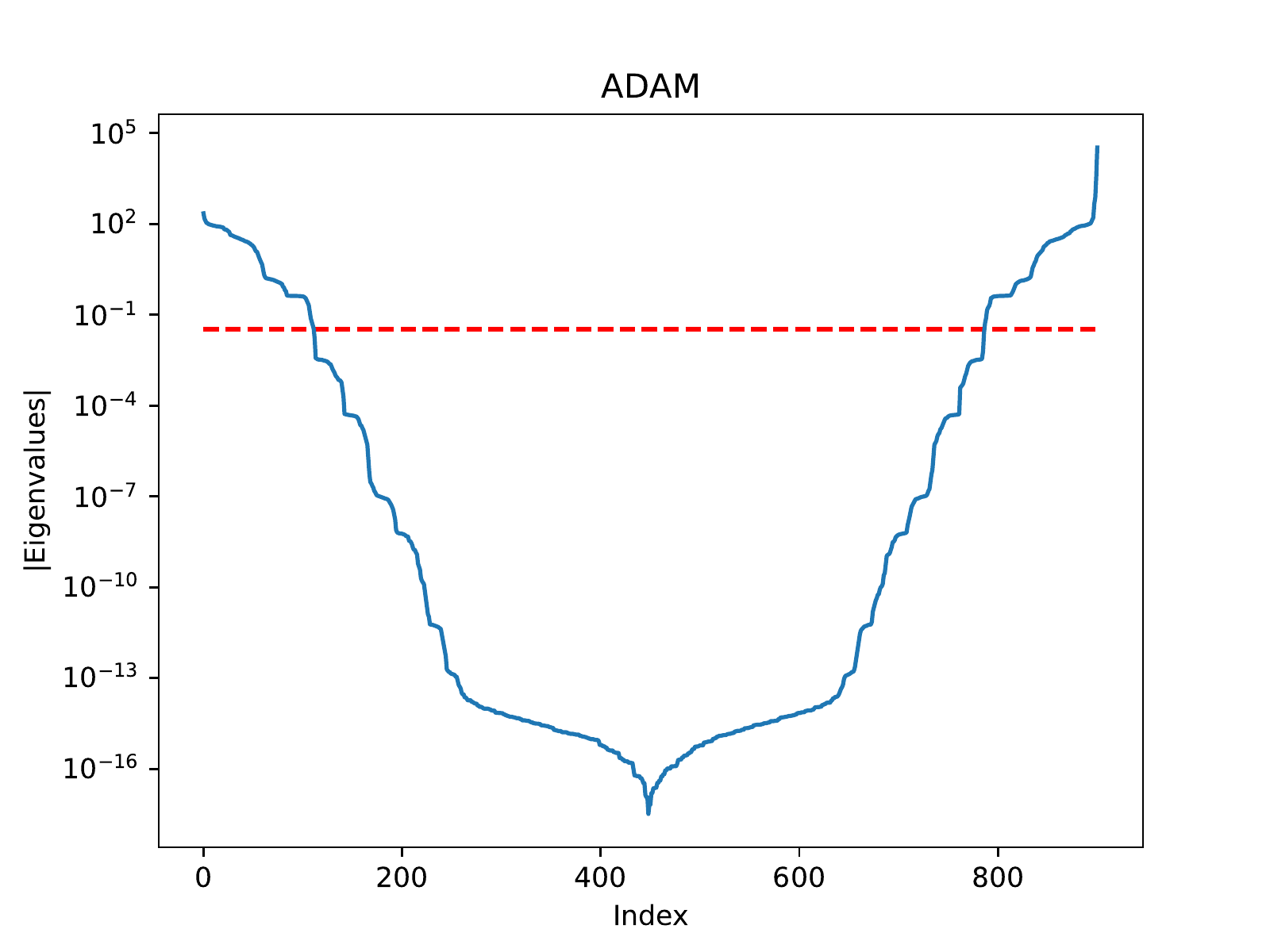}~
  \includegraphics[width=0.45\textwidth]{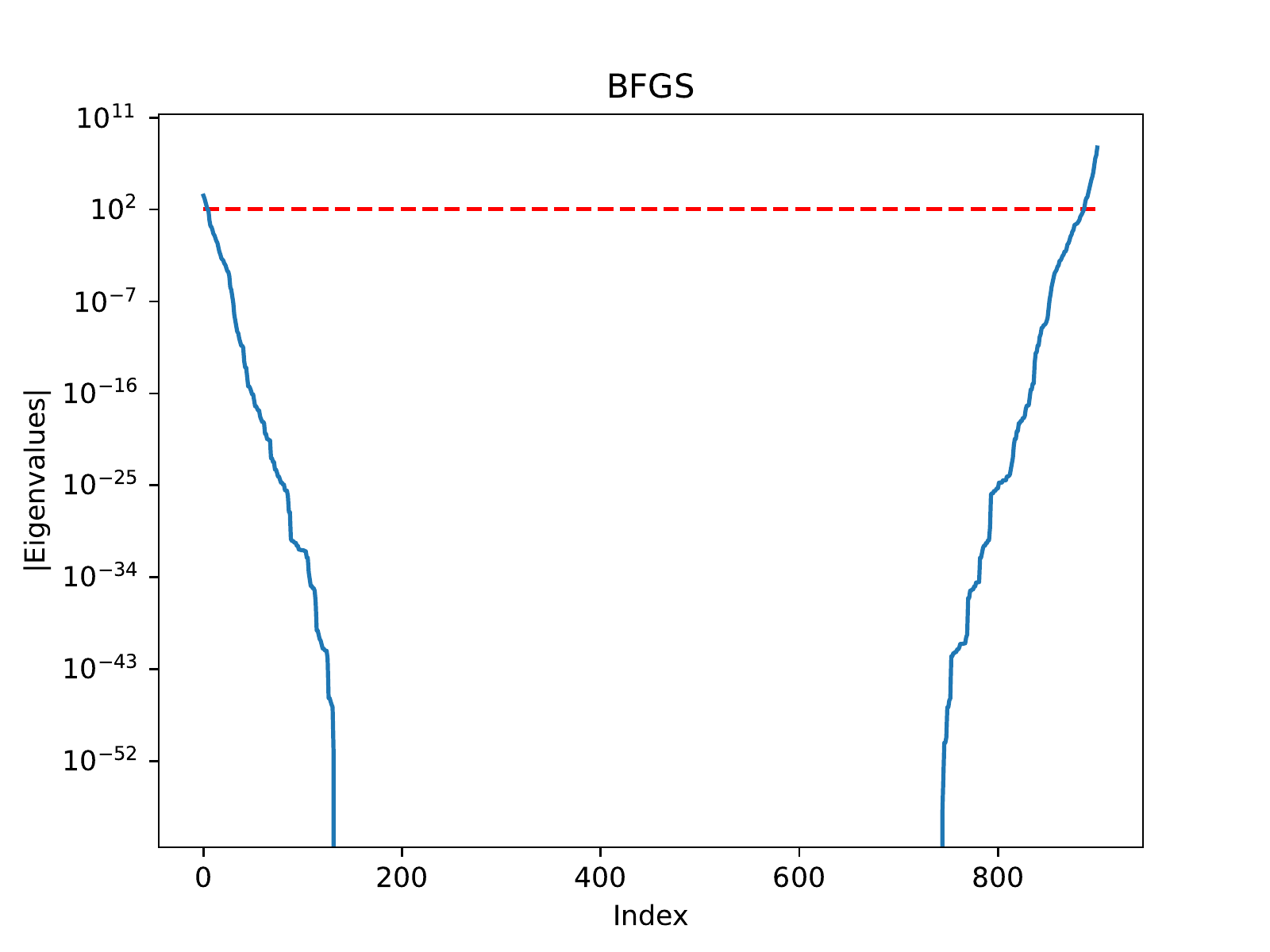}
   \includegraphics[width=0.45\textwidth]{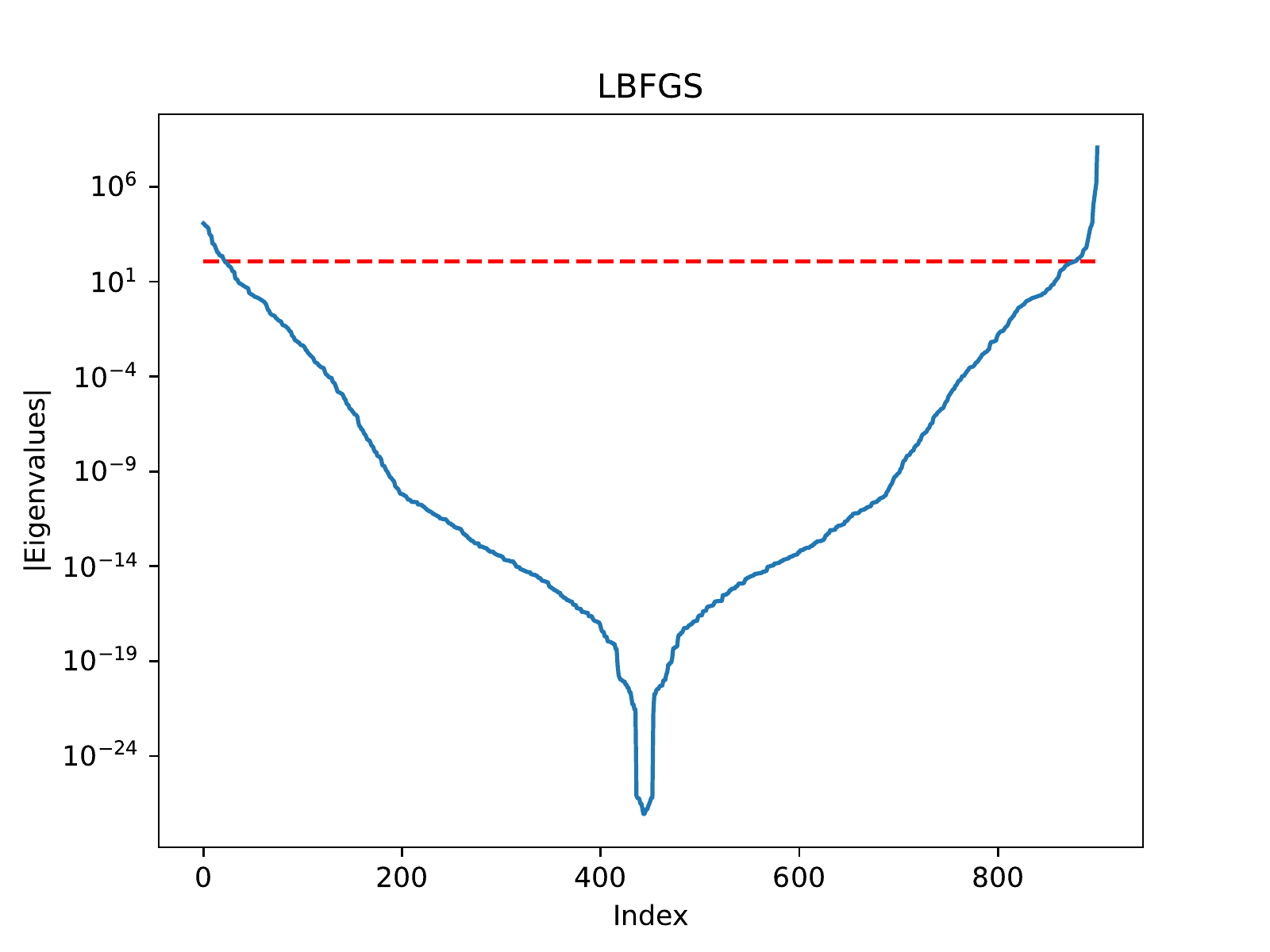}~
  \includegraphics[width=0.45\textwidth]{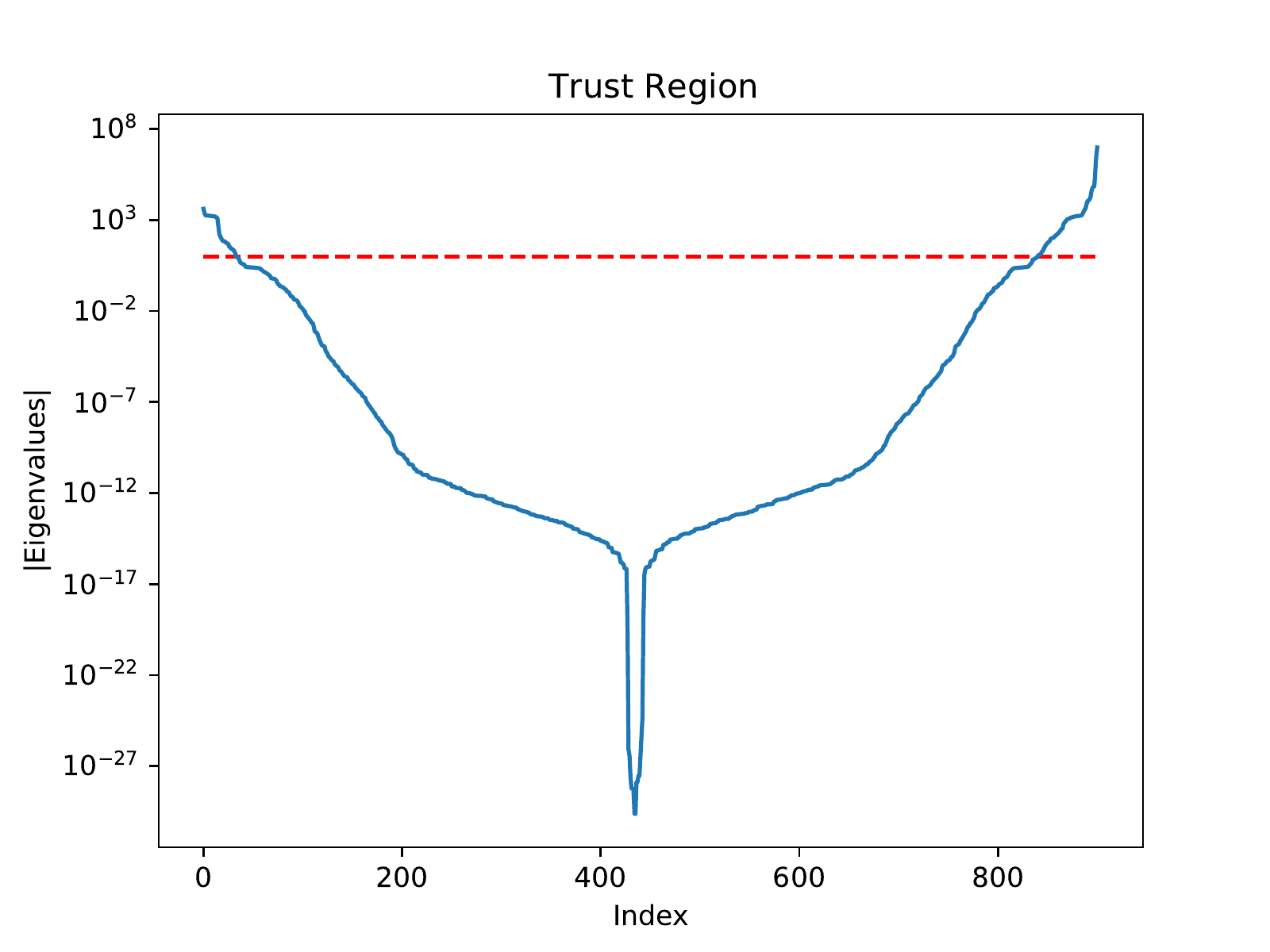}
  \caption{Magnitudes of eigenvalues at the last iteration for \Cref{equ:second-order-ex1-dnn-loss}. This plot should be compared to \Cref{fig:second-order-ex1-eigenvalues}, which combines the effect of DNNs and PDE solvers. The difference is the result of PDE constraints.}
  \label{fig:second-order-ex1-eigenvalue-dnn}
\end{figure}

We also investigate the effect of PDEs on effective DOFs. To this end, we calculate the Hessians of the ``DNN loss function''
\begin{equation}\label{equ:second-order-ex1-dnn-loss}
    l(\theta) = \sum_{i,j} (\kappa_\theta(\hat u_{i,j}) - \kappa(\hat u_{i,j}))^2
\end{equation}
The weights and biases are the converged values for \Cref{equ:second-order-ex1-optimization}, not a local minimizer of $\min l(\theta)$, so that we can inspect the effect of PDE constraints. 
The difference between \Cref{equ:second-order-ex1-optimization} and \Cref{equ:second-order-ex1-dnn-loss} is whether the PDE solver is considered in the loss function. The result is shown in \Cref{fig:second-order-ex1-eigenvalue-dnn}. We can see that the Hessian possesses some negative eigenvalues. This implies that the DNN and DNN-PDE loss functions indeed have different curvature structures at the local minimum. The structure is altered by the PDE constraint. Additionally, we can see that there are more effective DOFs for \Cref{equ:second-order-ex1-dnn-loss}, which indicates the PDE solver reduces effective DOFs and constrains the loss function on a smaller latent space. 

\begin{figure}[htbp]
  \centering
  \includegraphics[width=0.45\textwidth]{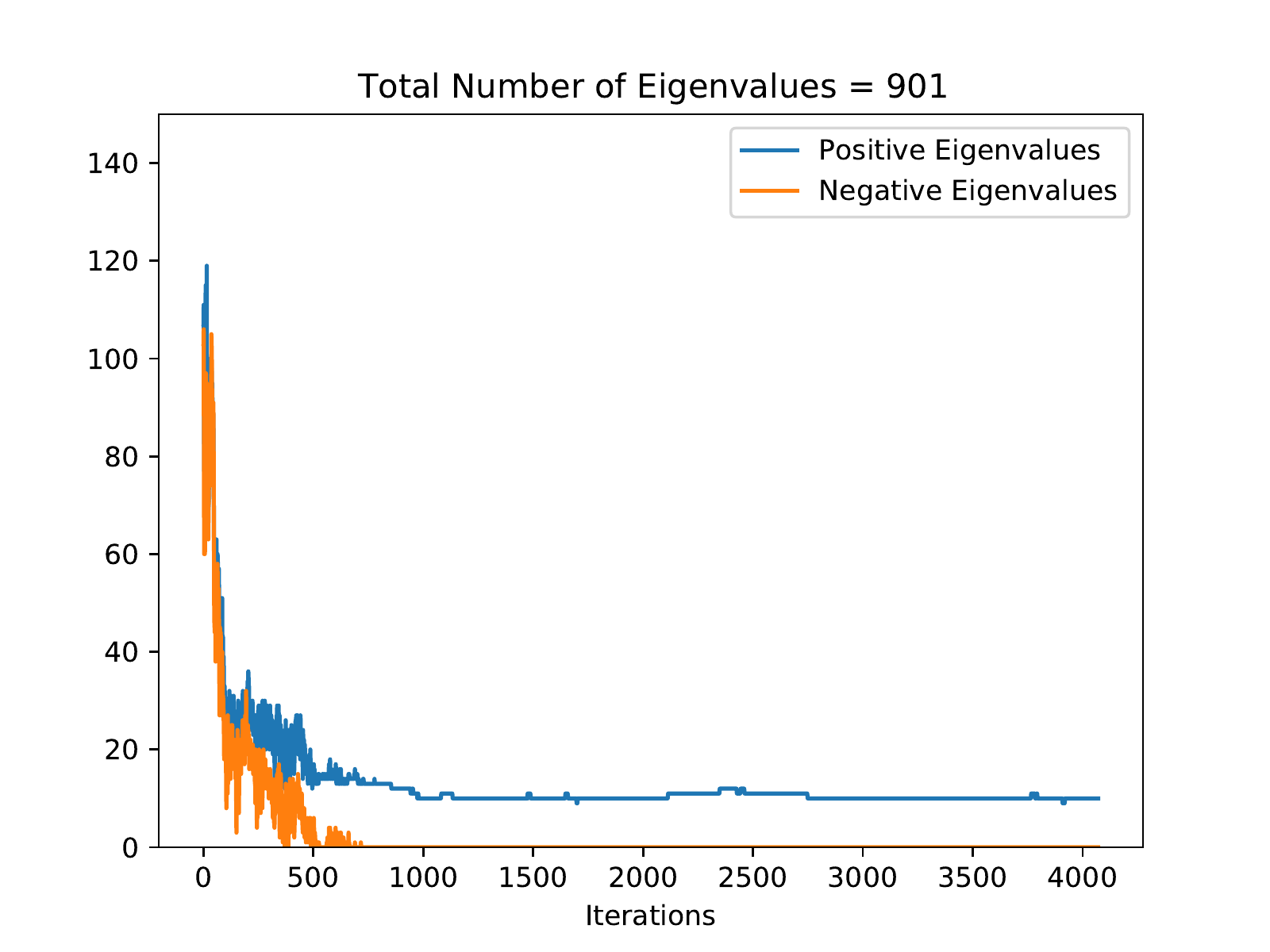}~
  \includegraphics[width=0.45\textwidth]{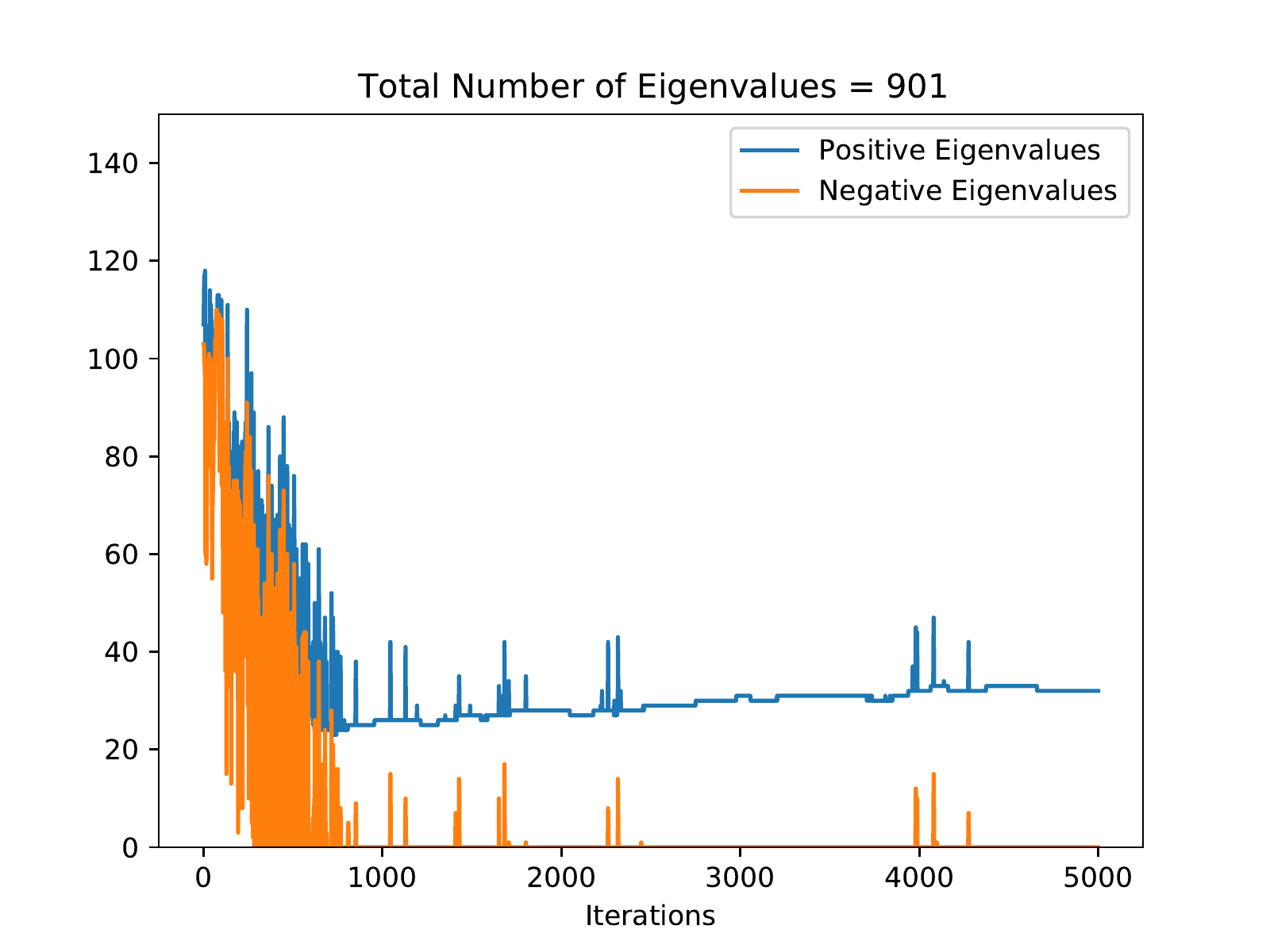}
  \caption{The number of positive and negative eigenvalues for the BFGS (left) and trust region (right) optimizer. We can see that the number of positive eigenvalues stays at around 18 and 30 for BFGS and trust region methods after a sufficient number of iterations. The spikes in the right panel are because, in trust region methods, the optimal solution to the local subproblem does not necessarily correspond to a better minimizer for the objective function. For example, the loss function values might increase.}
  \label{fig:second-order-ex1-eigenvalue-pn}
\end{figure}

\Cref{fig:second-order-ex1-eigenvalue-pn} shows the number of positive and negative eigenvalues for the BFGS and trust region optimizer. We can see that the number of positive eigenvalues stays at around 18 and 30 for BFGS and trust region methods after a sufficient number of iterations. The number of negative eigenvalues is nearly zero. This means that both optimizers converge to points with positive semi-definite Hessian matrices. Stationary points are true local minima, instead of saddle points.

\begin{figure}[htbp]
  \centering
  \includegraphics[width=0.45\textwidth]{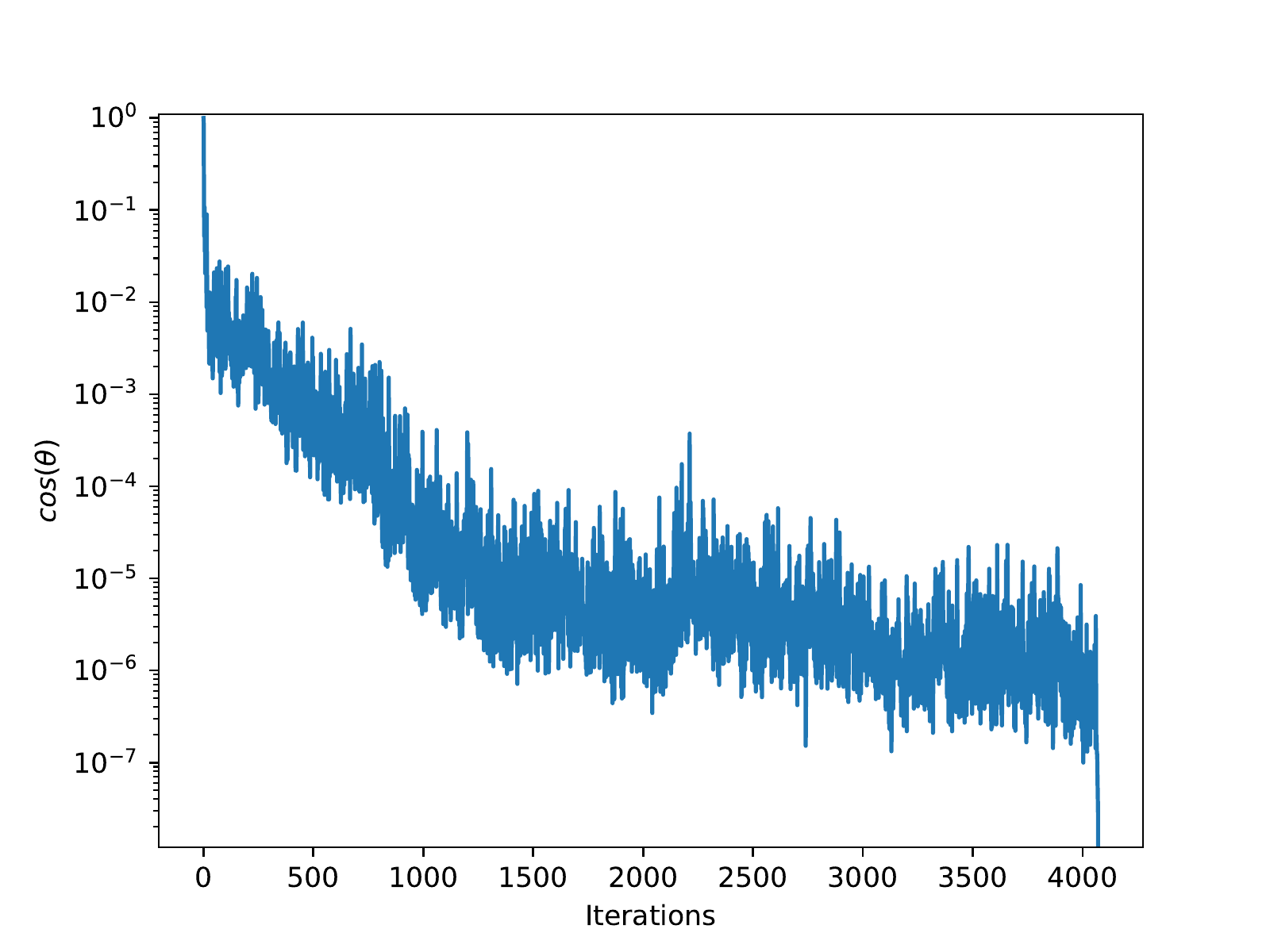}~
  \includegraphics[width=0.45\textwidth]{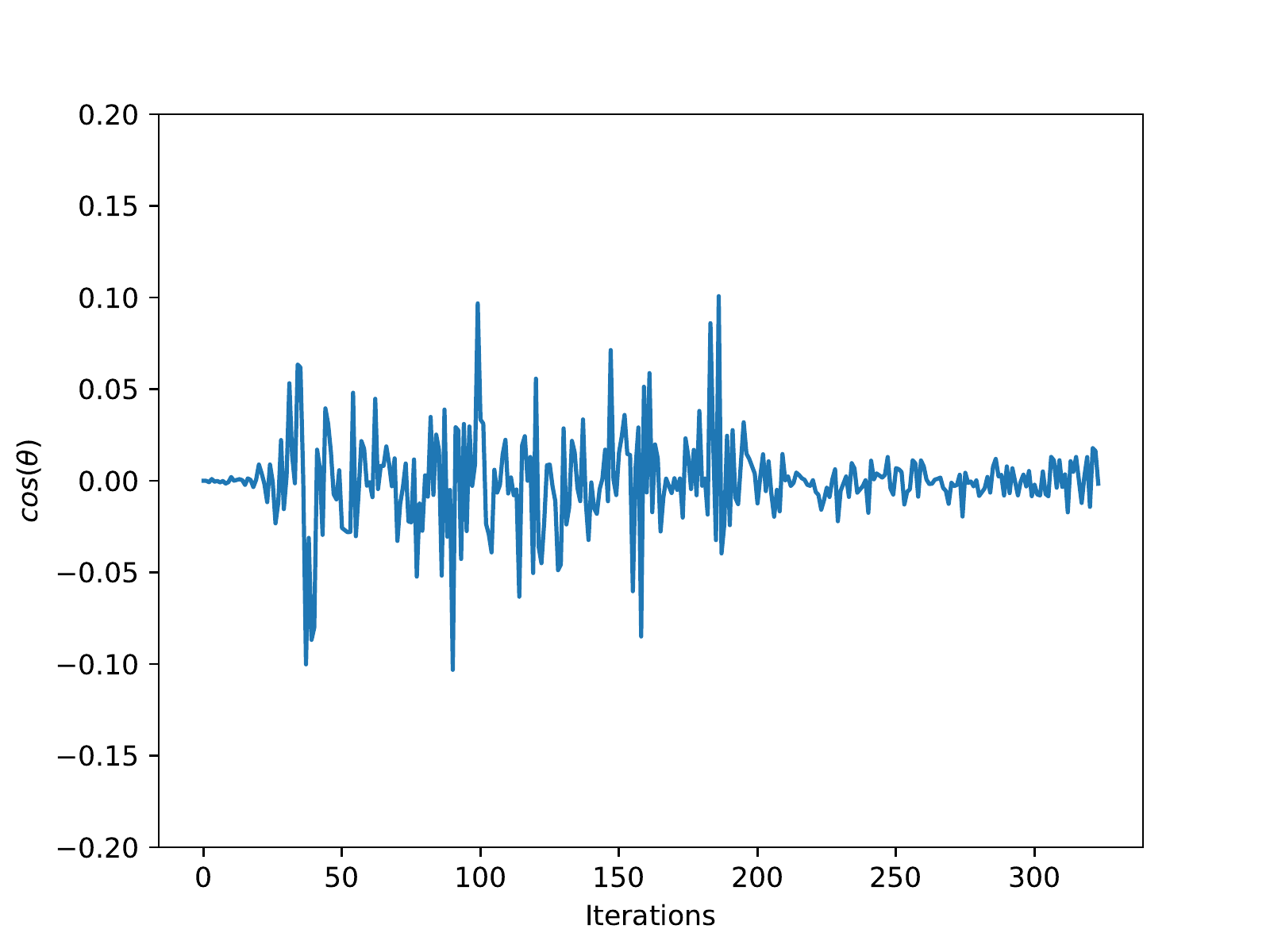}
  \caption{Left: the cosine values of angles between negative gradient directions and search directions in BFGS ($\cos\theta_1$). Right: the cosine values of angles between the Newton step and search directions in BFGS ($\cos\theta_2$). The right plot terminates at around 300 iterations because the Hessian matrix becomes ill-conditioned after a few iterations and the Newton step is not properly defined.}
  \label{fig:second-order-ex1-angles}
\end{figure}

We also analyze the direction of the search direction $p_k$ in the BFGS optimizer. We consider two values
$$\begin{aligned}\cos(\theta_1) &= \frac{-p_k^T g_k}{|p_k||g_k|} \qquad  \cos(\theta_2) &= \frac{p_k^T q_k}{|p_k||q_k|}\end{aligned}$$
Here $q_k$ is the direction for the Newton step
$$q_k = -H_k^{-1}g_k$$
The two quantities are shown in \Cref{fig:second-order-ex1-angles} (note after some number of iterations, the Hessian $H_k$ becomes semi-positive definite, making it impossible to calculate $q_k$, and thus $\cos(\theta_2)$ only has limited data points). There are two conclusions to draw from the plots
\begin{enumerate}
    \item The search direction of the BFGS optimizer deviates from the gradient descent method.
    \item The search direction of the BFGS optimizer is not very correlated with the Newton step; this indicates the search direction poorly recognizes the negative curvature directions.
\end{enumerate}

In sum, we conclude that trust region methods are very valuable for our inverse problems, which possess special characteristics, such as the indefiniteness of Hessian matrices. In the following examples, we will reveal more about the structure of our problem and how trust region methods are suitable for exploiting these structures.

\subsection{Heat Equation}\label{sect:second-order-dynamic}

Next, we consider a time-dependent PDE problem, the heat equation, whose governing equation is given by 
\begin{equation}\label{equ:second-order-ex2}
\begin{aligned}
     \frac{\partial u}{\partial t} &= \nabla \cdot (\kappa(x, y) \nabla u))+  f(x, y) & x\in \Omega\\ 
     u(x,y,0) &= x(1-x)y^2(1-y)^2 & (x,y)\in \Omega\\
     u(x,y,t) &= 0 & (x,y)\in  \partial \Omega
\end{aligned}
\end{equation}
Here the spatially-varying diffusivity coefficient $\kappa(x,y)$ is given by 
$$\kappa(x,y) = 2x^2 - 1.05x^4 + x^6 +xy+y^2$$
The exact solution is given by 
$$u(x,y,t) = x(1-x)y^2(1-y)^2 e^{-t}$$

We assume that we can observe the full field data of $u$ with no noise as snapshots. We want to use these observations to estimate $\kappa(x,y)$. We use a DNN $\kappa_\theta(x,y)$, where $\theta$ is the weights and biases, to approximate $\kappa$. 

We again apply the residual minimization method to train the deep neural network. The optimization problem can be written as 
\begin{equation}\label{equ:second-order-ex2-optimization}
    \min_\theta\; \sum_n\sum_{i,j} \left(\frac{u_{i,j}^{n+1} - u_{i,j}^n}{\Delta t} - F_{i,j}( u^{n+1}; \theta) -  f^{n+1}_{i,j}\right)^2 
\end{equation}
Here the superscript $n$ denotes time step, and $\Delta t > 0$ is the time step size.

\begin{figure}[htbp]
  \centering
  \includegraphics[width=0.45\textwidth]{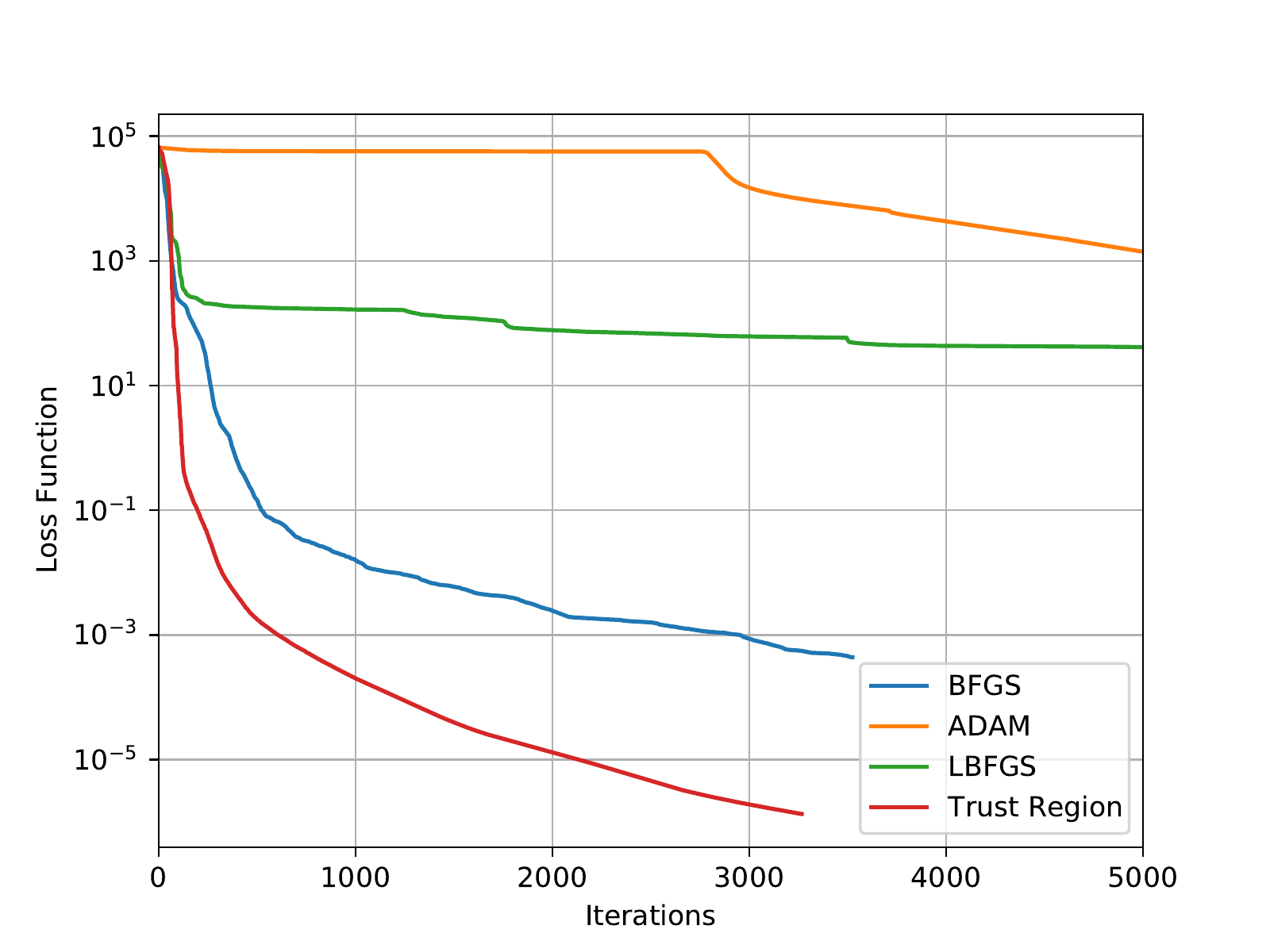}~
  \includegraphics[width=0.45\textwidth]{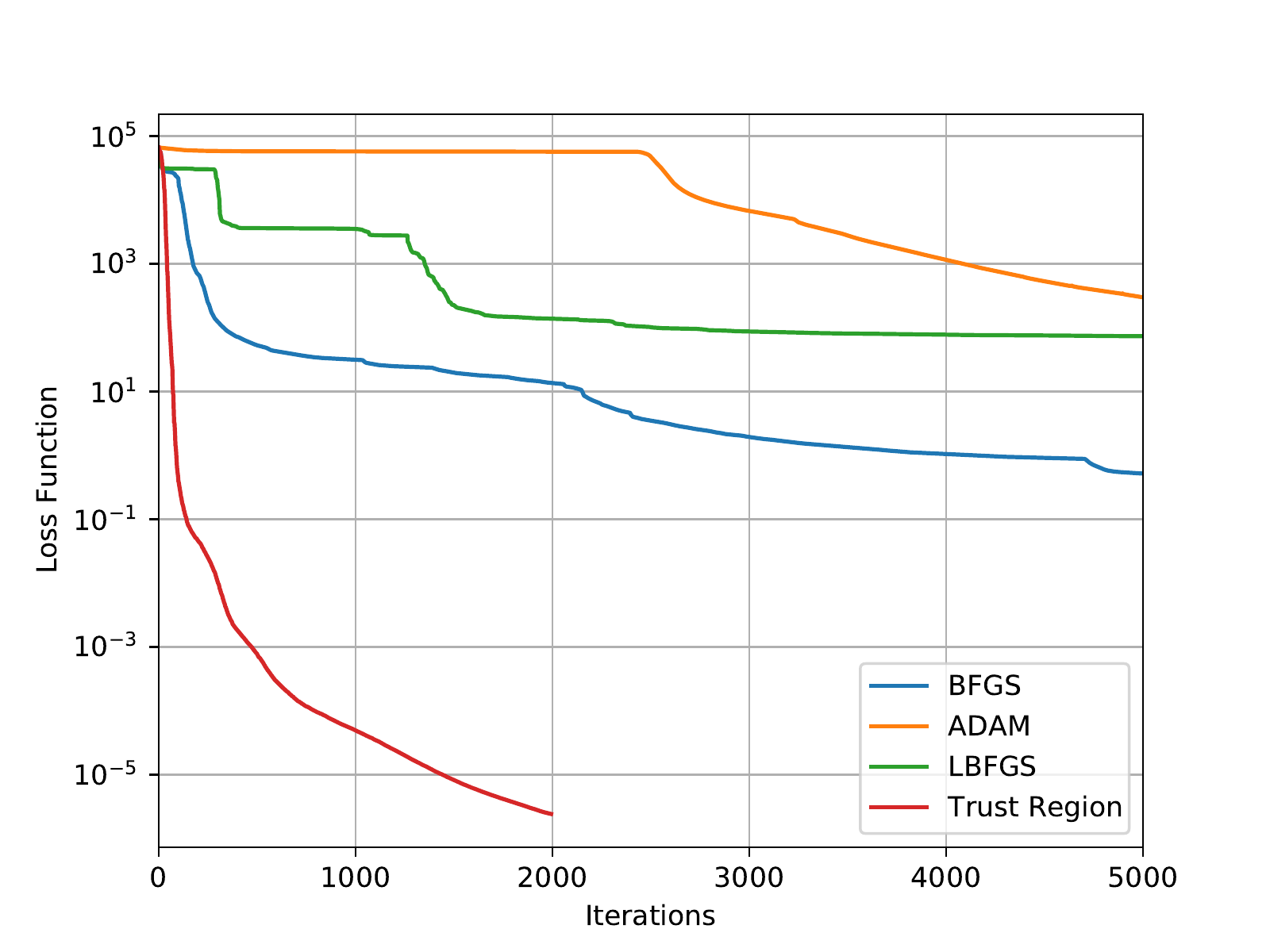}
   \includegraphics[width=0.45\textwidth]{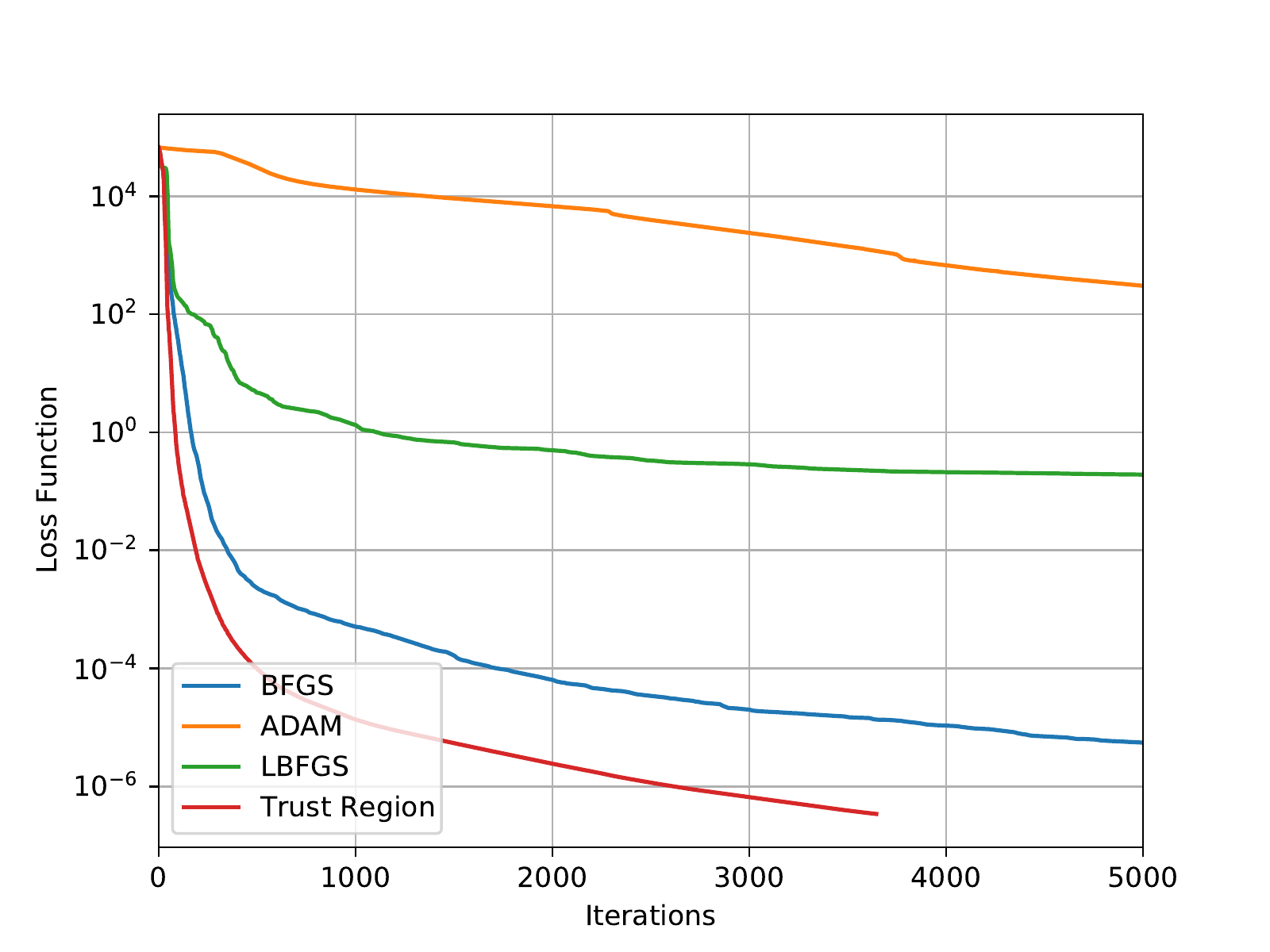}~
  \includegraphics[width=0.45\textwidth]{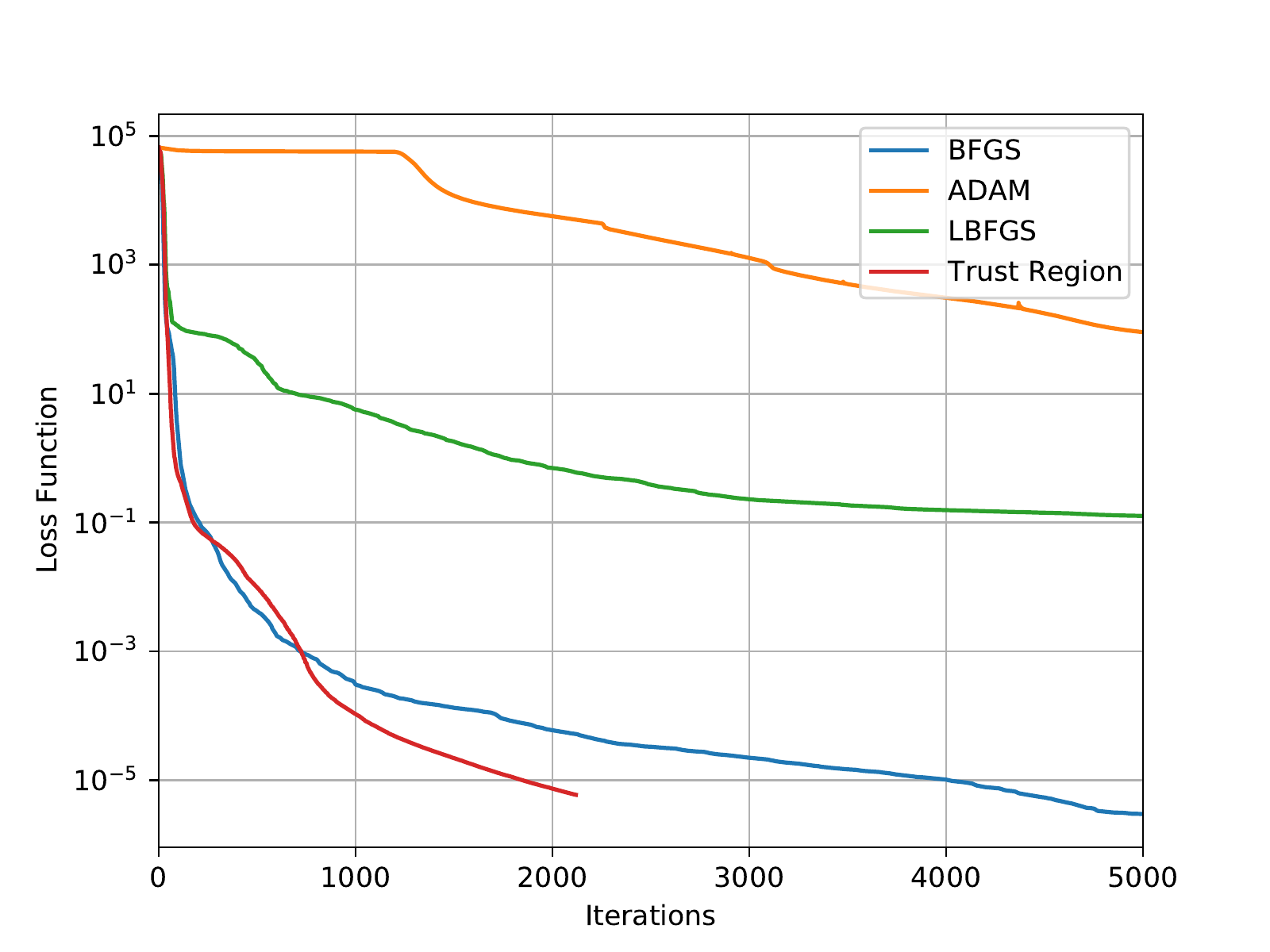}
  \caption{Loss functions for solving \Cref{equ:second-order-ex2-optimization} using different optimizers. Trust region methods perform consistently better than all other three methods.}
  \label{fig:second-order-ex2-loss}
\end{figure}

\paragraph{Convergence} \Cref{fig:second-order-ex2-loss} shows the convergence plots for different initial guesses of the DNNs.
We see that the trust region methods are more competitive than the other methods, which is consistent with our finding in \Cref{sect:second-order-static}. 

\begin{figure}[htbp]
  \centering
  \includegraphics[width=0.33\textwidth]{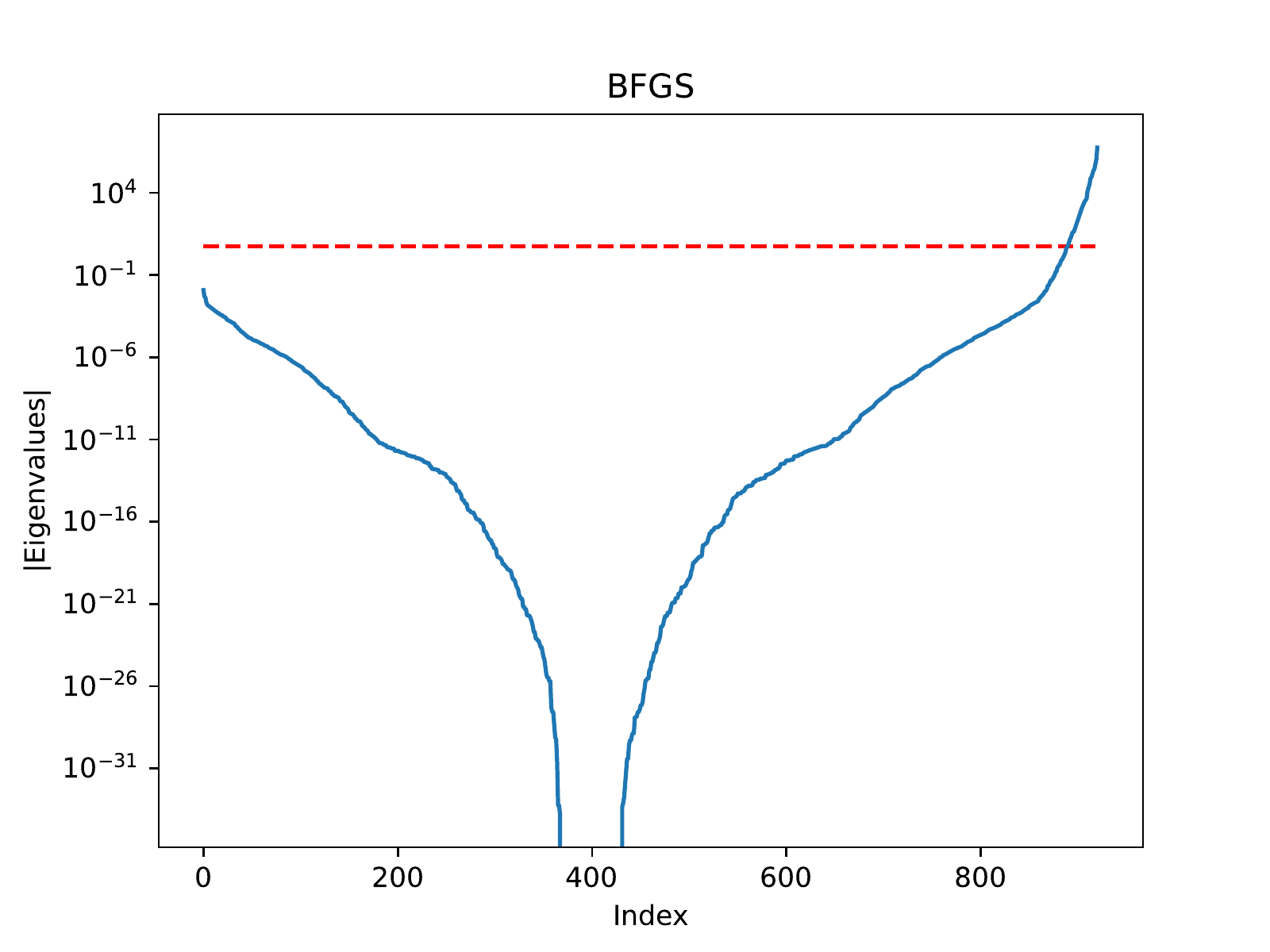}~
  \includegraphics[width=0.33\textwidth]{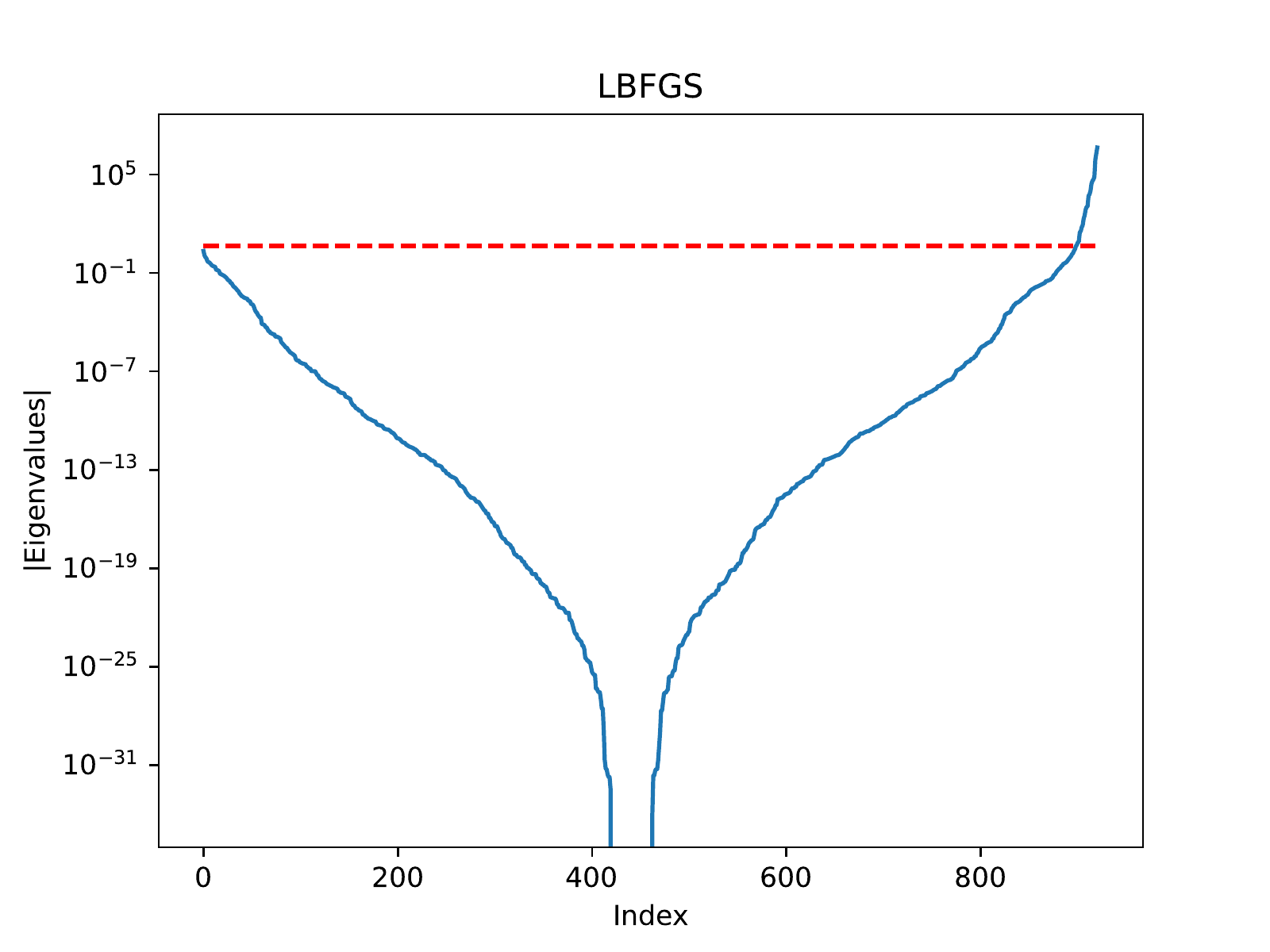}~
   \includegraphics[width=0.33\textwidth]{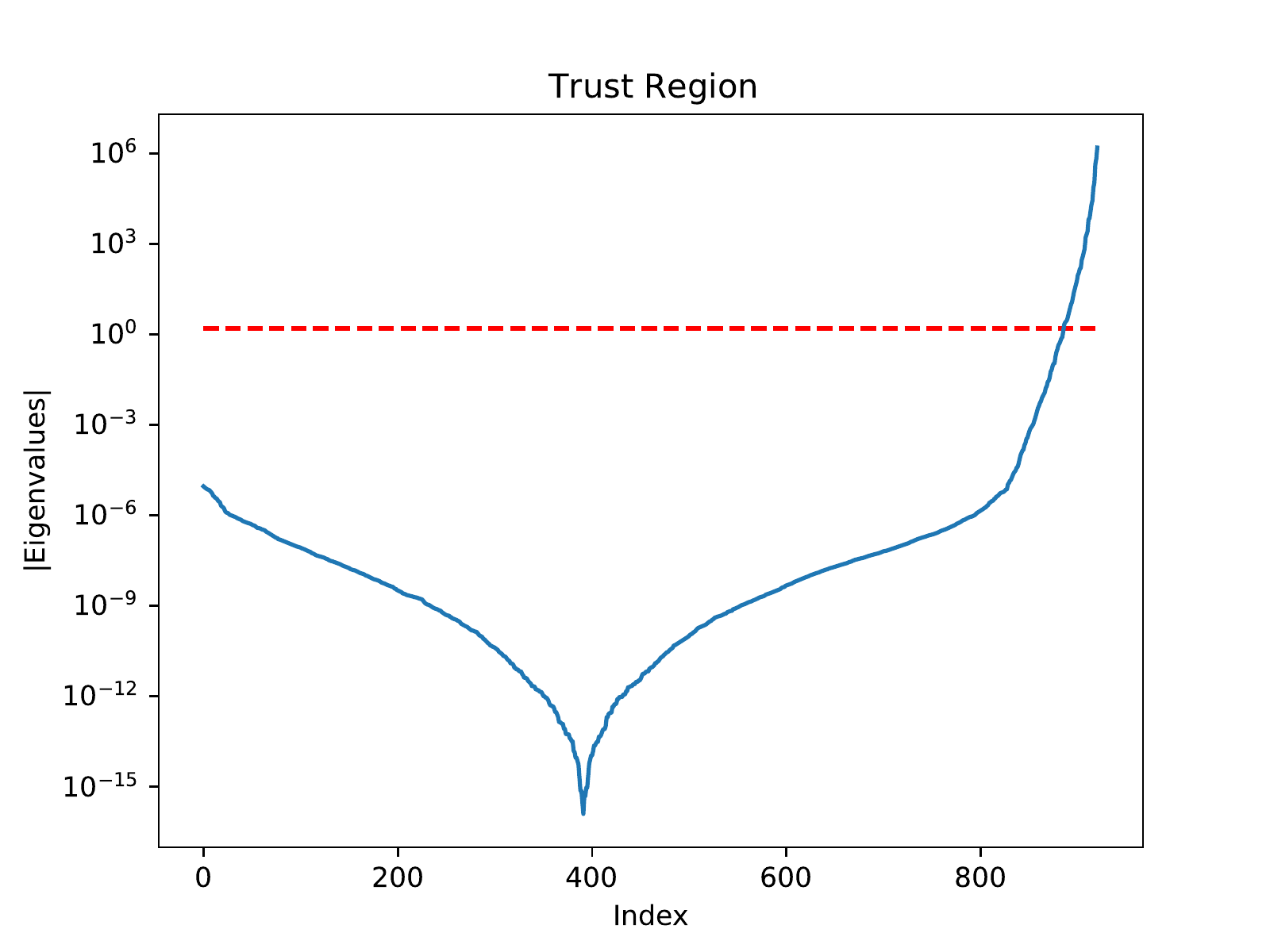}
   \includegraphics[width=0.33\textwidth]{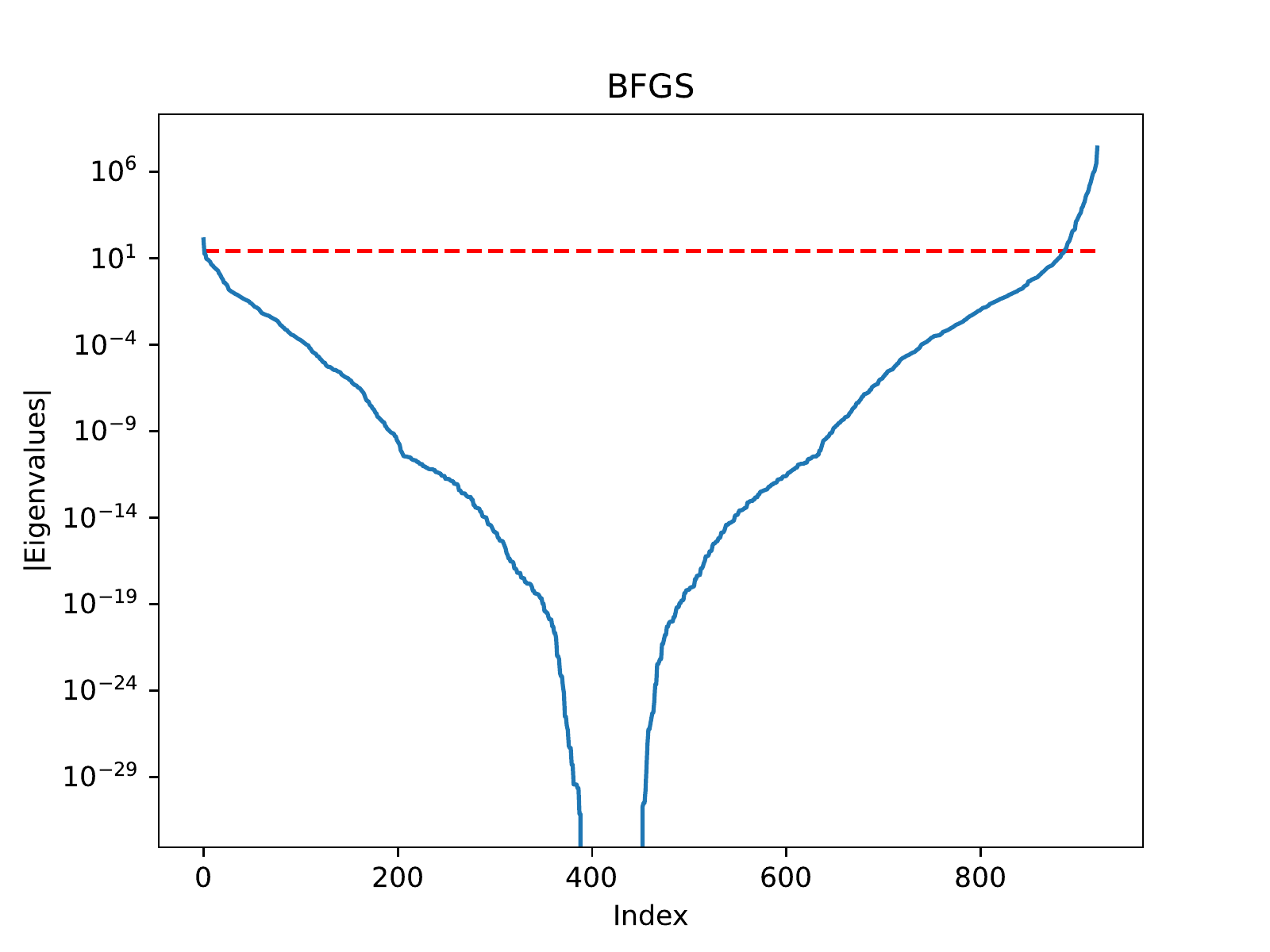}~
  \includegraphics[width=0.33\textwidth]{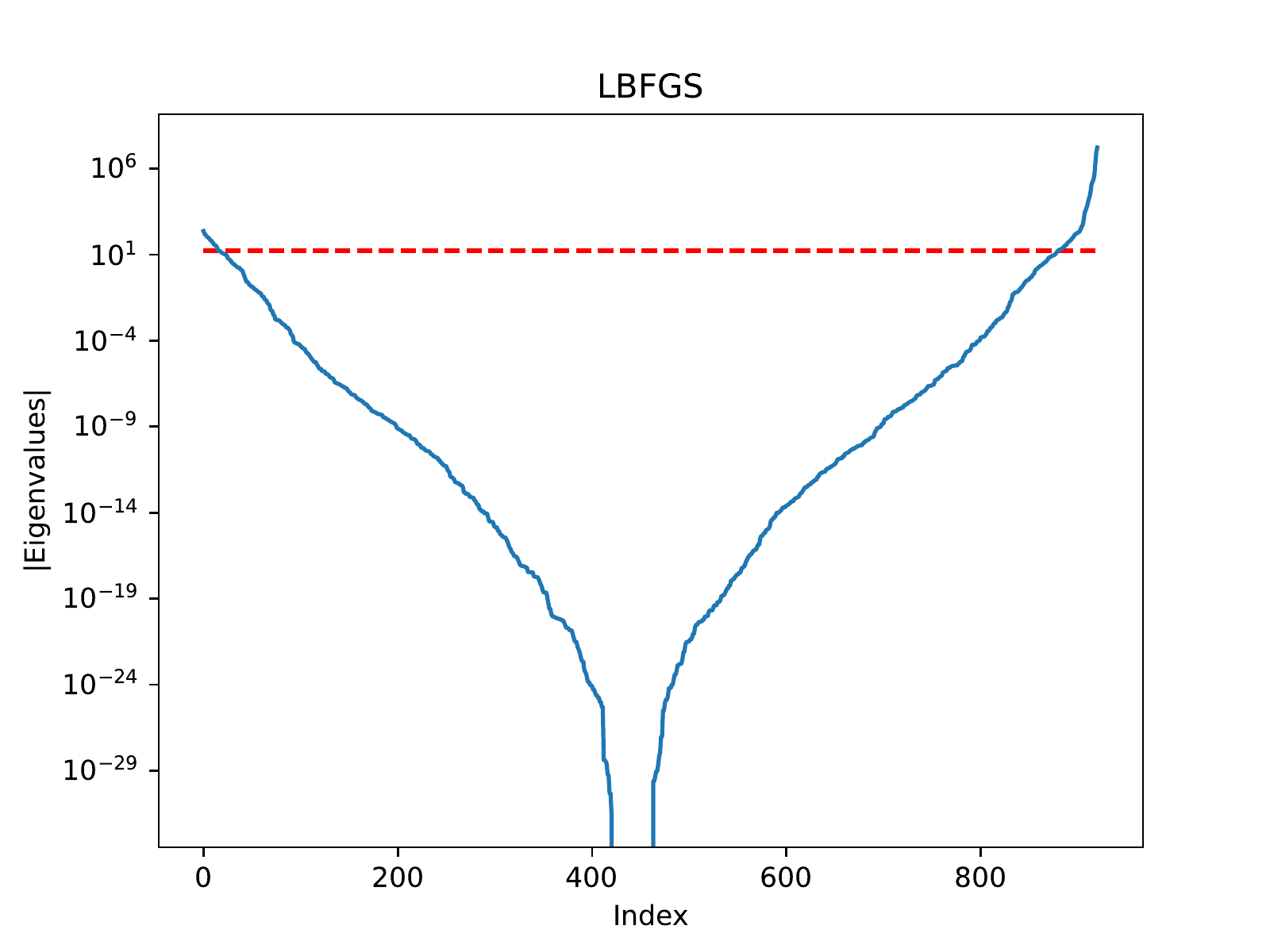}~
   \includegraphics[width=0.33\textwidth]{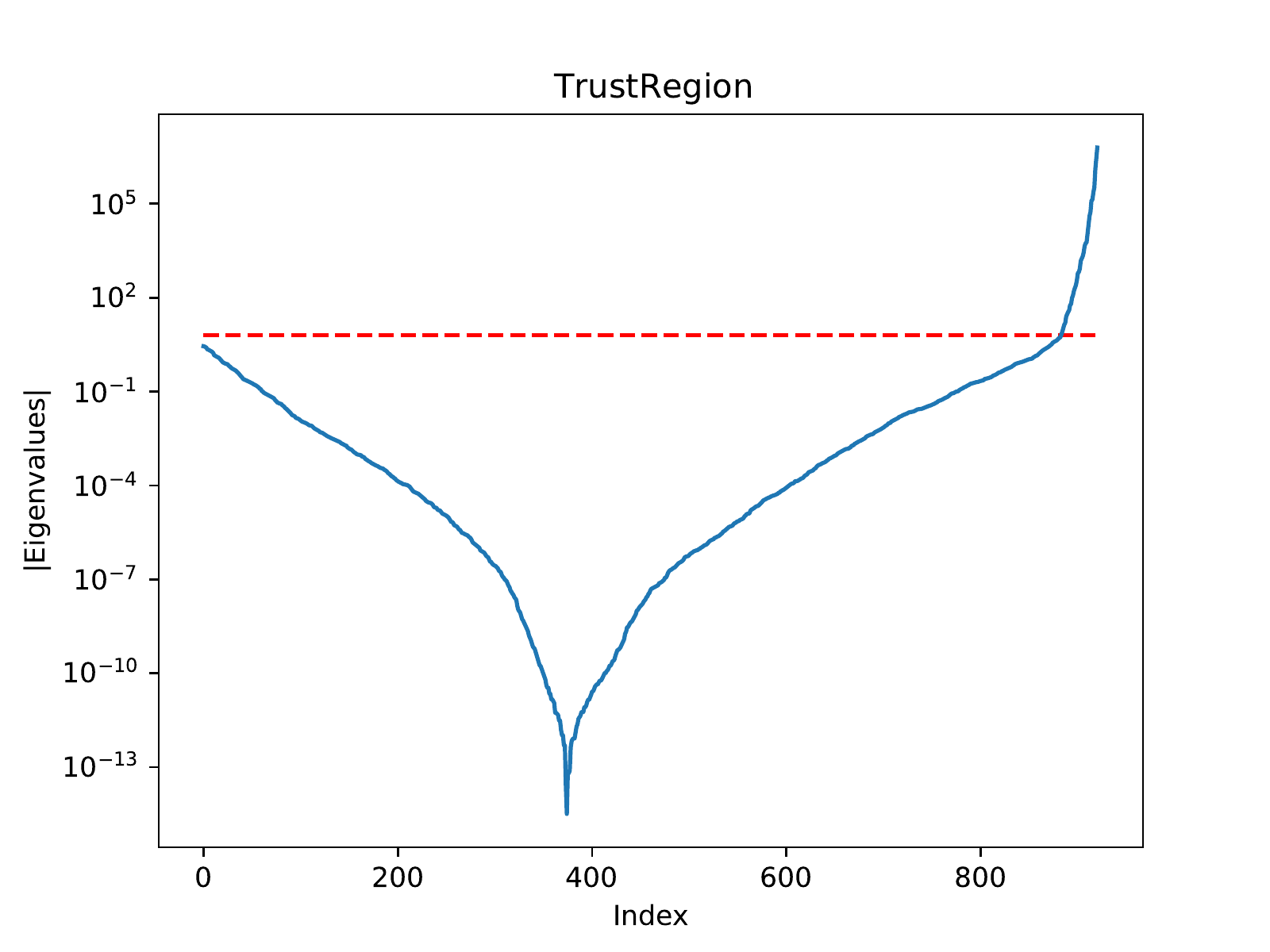}
  \caption{Magnitudes of eigenvalues at the last iteration for \Cref{equ:second-order-ex2-optimization}. The first row corresponds to the loss function in \Cref{equ:second-order-ex2-optimization} (coupled DNN and PDE solver). The second row corresponds to the loss function in \Cref{equ:second-order-ex2-dnn-loss}, where the PDE solver is removed from the loss function.}
  \label{fig:second-order-ex2-eigenvalues}
\end{figure}

\paragraph{Effect of PDEs} In \Cref{fig:second-order-ex2-eigenvalues}, we show the magnitudes of eigenvalues for the Hessian matrix of the loss function in \Cref{equ:second-order-ex2-optimization} and the DNN loss function
\begin{equation}\label{equ:second-order-ex2-dnn-loss}
    l(\theta) = \sum_{i,j} (\kappa_\theta(x_{i,j}, y_{i,j}) - \kappa(x_{i,j}, y_{i,j}))^2
\end{equation}
where $(x_{i,j}, y_{i,j})$ is the grid points. \Cref{fig:second-order-ex2-eigenvalues} exhibits a similar pattern as \Cref{fig:second-order-ex1-eigenvalues} and \Cref{fig:second-order-ex1-eigenvalue-dnn}: a few positive eigenvalues accompanied by zero eigenvalues for both Hessian matrices. The difference between Hessian matrices for DNN-PDE and DNN is that the number of positive eigenvalues are slightly larger than the loss function that couples DNNs and PDEs (see \Cref{tab:second-order-ex1-eigenvalues}). This implies that PDEs restricts effective DOFs. We attribute the diminished effective DOFs to the physical constraints imposed by PDEs.

\begin{table}[htpb]
\caption{Number of positive eigenvalues (effective DOFs) for different optimizers. ``DNN-PDE'' corresponds to the loss function in \Cref{equ:second-order-ex2-optimization} and ``DNN Only'' corresponds to the loss function \Cref{equ:second-order-ex2-optimization}.}
\label{tab:second-order-ex1-eigenvalues}
\centering
\begin{tabular}{@{}lllll@{}}
\toprule
         & ADAM & BFGS & LBFGS & Trust Region \\ \midrule
DNN-PDE  & 50   & 31   & 22    & 35           \\
DNN Only & 132  & 34   & 41    & 38           \\ \bottomrule
\end{tabular}
\end{table}

\paragraph{Different Optimizers} Now we consider the difference between different optimizers. The numerical evidence in this and last sections show that ADAM, BFGS, LBFGS, and trust region methods all converge to a local minimum. But why do different local minima have different effective DOFs? What's their implication on the goodness of local minimum?

\begin{figure}[htbp]
  \centering
  \includegraphics[width=0.6\textwidth]{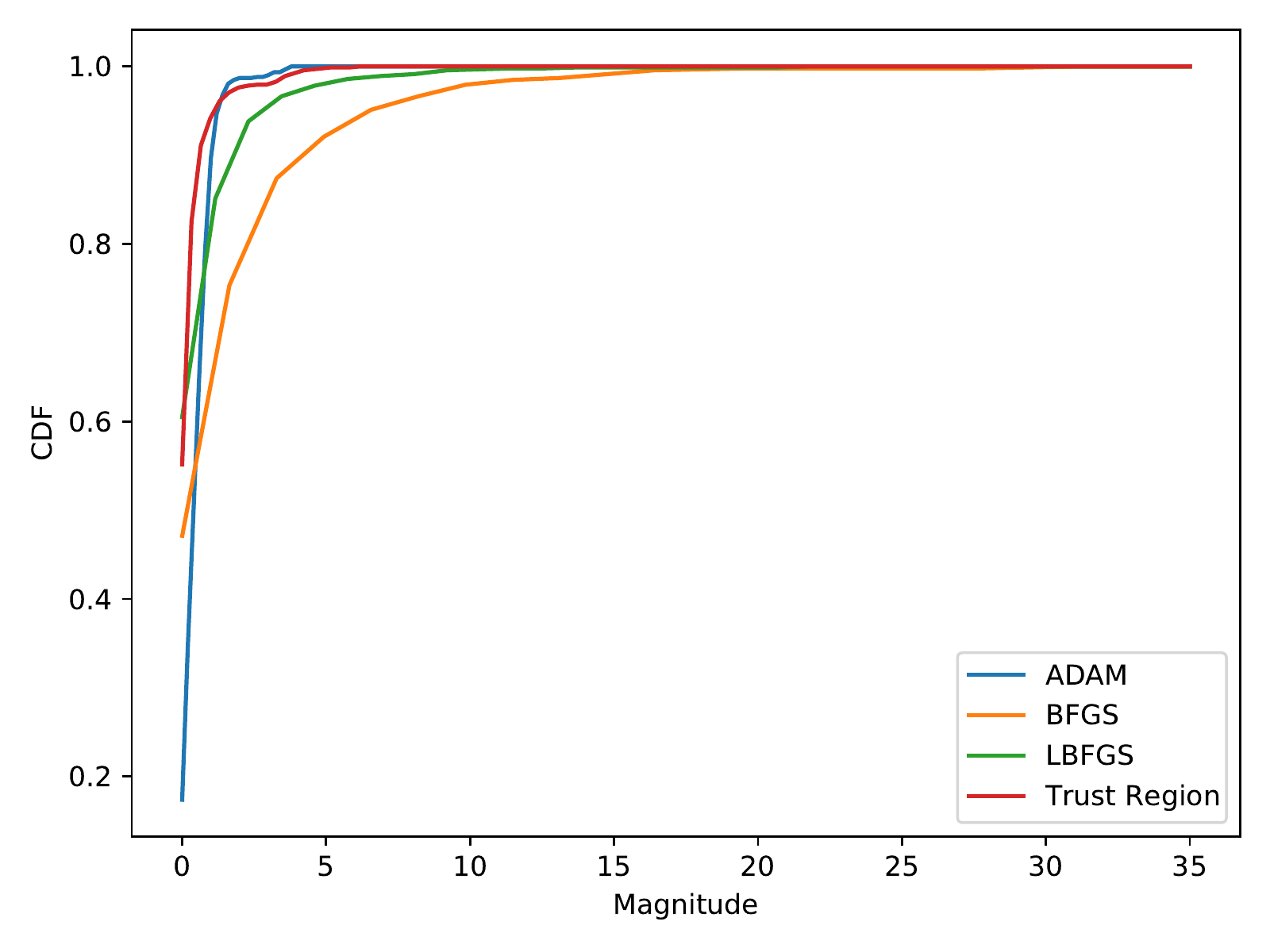}
  \caption{The cumulative distribution of the maginitude for different optimizers. BFGS and LBFGS are much more concentrated  than ADAM and trust region methods in terms of maginitudes.}
  \label{fig:second-order-weight-dist}
\end{figure}

To answer this question, we show the cumulative distribution of the magnitude of weights and biases in \Cref{fig:second-order-weight-dist}.

The plot shows that parameters in BFGS and LBFGS are much more concentrated than ADAM and trust region methods in terms of magnitudes. Because we use $\tanh$ as activation values, for fixed intermediate activation values, large weights and biases are more likely to cause saturation of activation values, i.e., the inputs to $\tanh$ are large or small and thus the outputs are close to 1. To illustrate the idea, consider a simple function
$$y = w_1 \tanh(w_2 x + b_2) + b_1$$
Given a reasonable $x$ (e.g., $x\approx 0.1$), if $|w_2|$ or $|b_2|$ is large, $y \approx b_1 \pm w_1$, and thus the effective DOF is 2; if $w_2$ and $b_2$ is close to 0, $y\approx w_1 w_2 x + w_1 b_2 + b_1$, perturbation of all four parameters $w_1$, $w_2$, $b_1$, $b_2$ may contribute to the change of $y$, and thus the effective DOF is 4. In sum, trust region methods yield weights and biases with smaller magnitudes compared to BFGS/LBFGS in general, and thus achieve more effective DOFs.

\begin{figure}[htbp]
  \centering
  \includegraphics[width=0.6\textwidth]{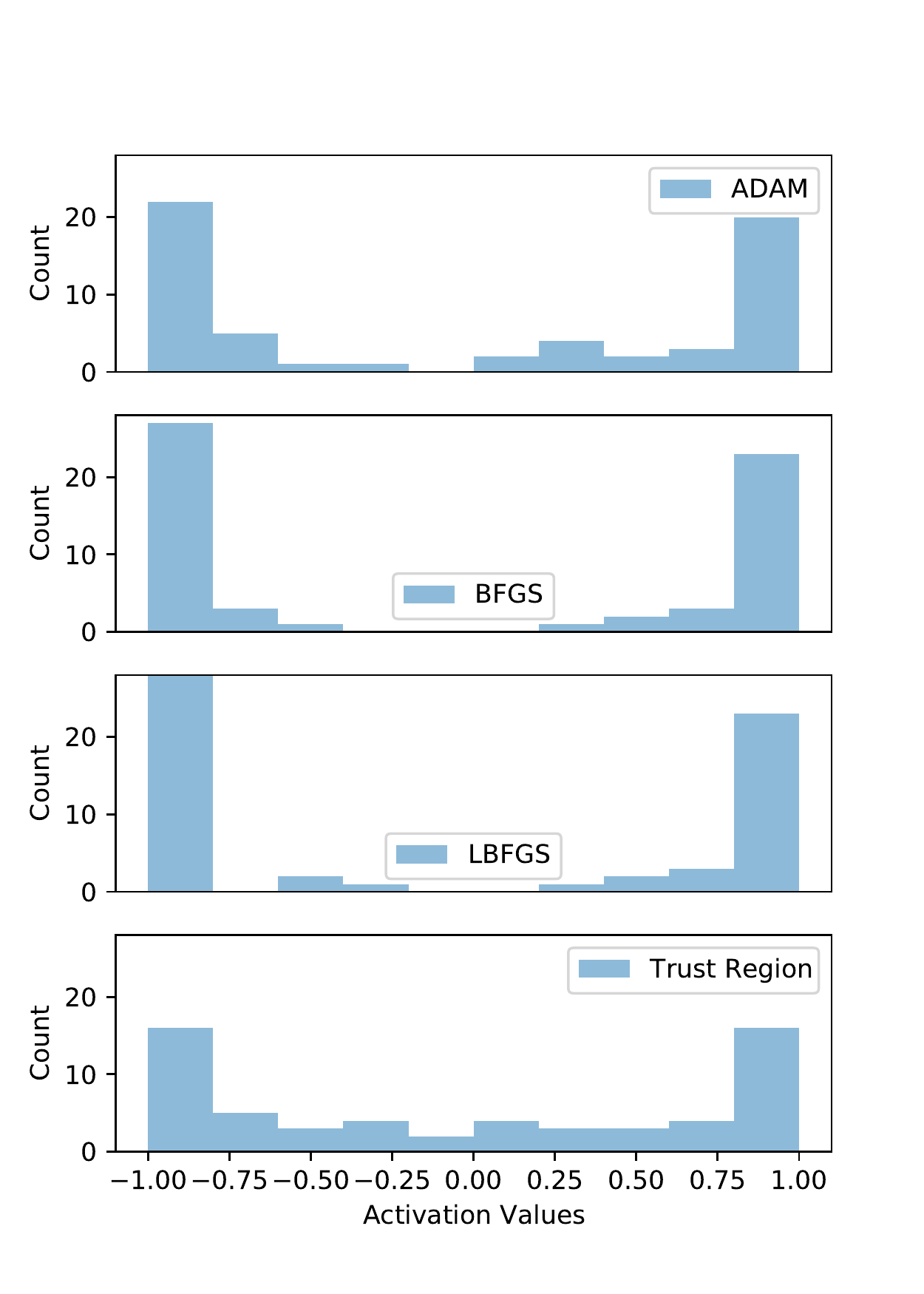}
  \caption{The histogram of intermediate activation values. Note the range of $\tanh$ is between $-1$ and $1$. $\tanh$ is considered saturated if the output is close to $-1$ or $1$.}
  \label{fig:second-order-activation-dist}
\end{figure}

This conjecture is confirmed by \Cref{fig:second-order-activation-dist}, which shows the histogram of intermediate activation values. We fixed the input $x = (0.5,0.5)$ (the midpoint of the computational domain), and collected all the outputs of the $\tanh$ function within the DNN. The figure shows that compared to the trust region method, the activation values of ADAM, BFGS and LBFGS are more concentrated near the extreme values $-1$ and $1$.

How can trust region methods manage the magnitudes of the weights and biases? The benefit is intrinsic to how the trust region method works: it only searches for ``optimal solution'' with a small neighborhood of the current state. However, BFGS and LBFGS search for ``optimal solution'' along a direction aggressively. Given so many local minima, it is very likely that BFGS and LBFGS get trapped in a local minimum with smaller effective DOFs. In this perspective, trust region methods are useful methods for avoiding (instead of ``escaping'') bad local minima.

The take-away message here is that for optimization involving neural networks, aggressive line searching methods is undesirable. Good local minimum can be found with a one-step-at-a-time approach. For this reason, we consider trust region methods to be very competitive and appropriate for solving inverse problems involving DNNs. 

\paragraph{Widths and Depths of DNNs} We now consider the effects of width and depths of DNNs on the eigenvalue distributions of Hessians. In this experiment, we fixed the initial guess of weights and biases and solved the optimization problem \Cref{equ:second-order-ex2-optimization} using different optimizers. 

\begin{figure}[htbp]
  \centering
  \includegraphics[width=0.32\textwidth]{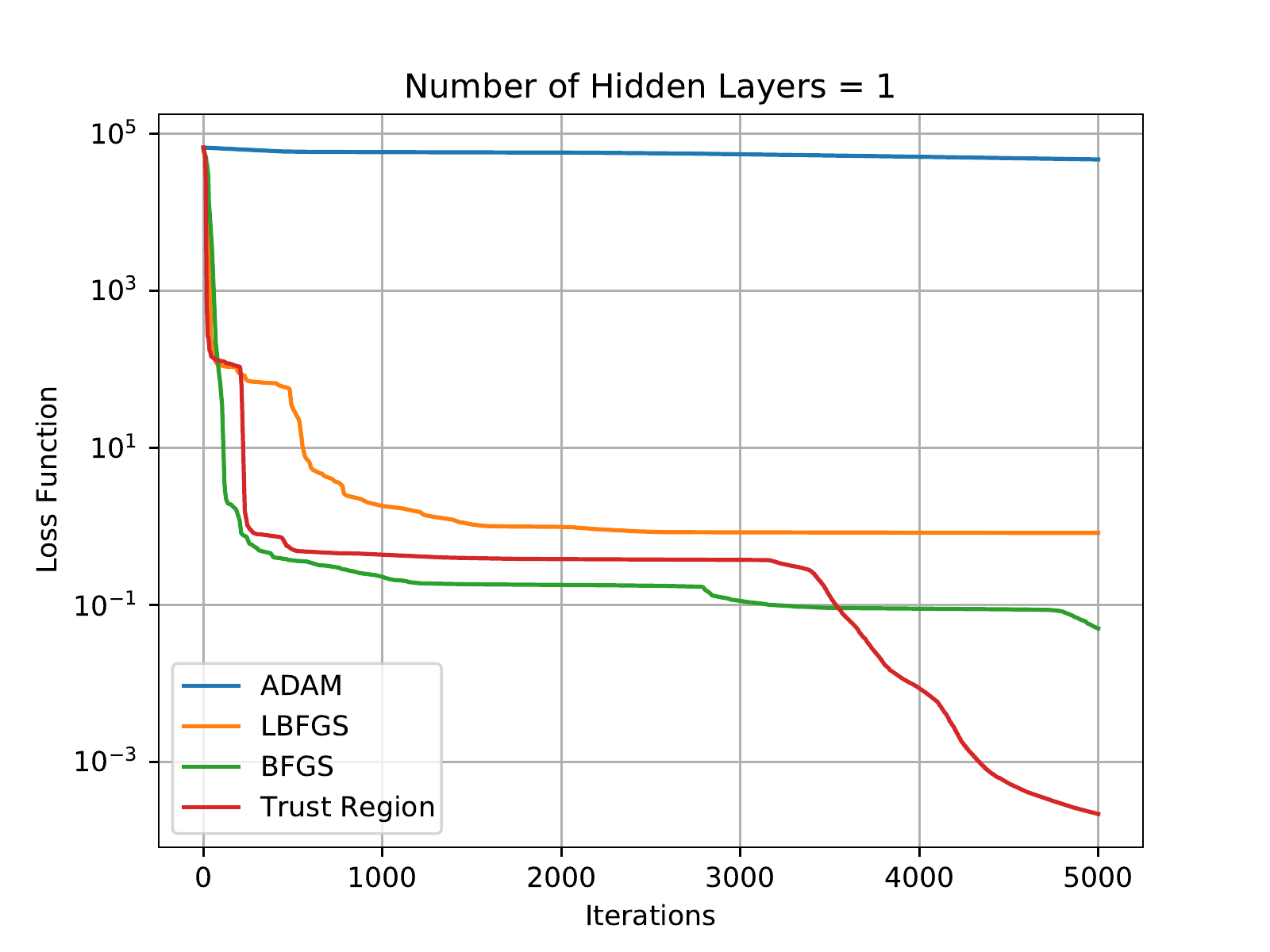}
  \includegraphics[width=0.32\textwidth]{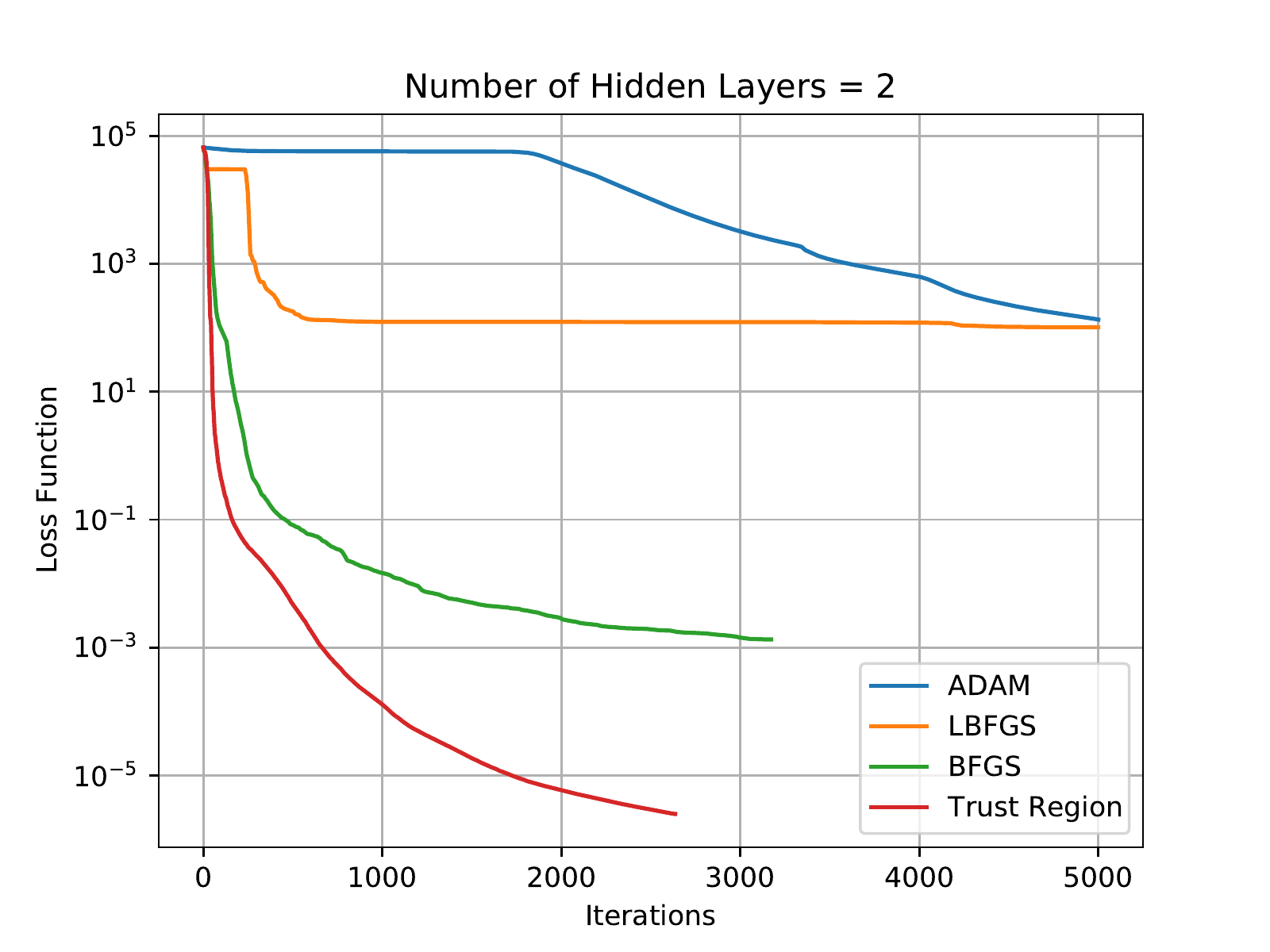}
  \includegraphics[width=0.32\textwidth]{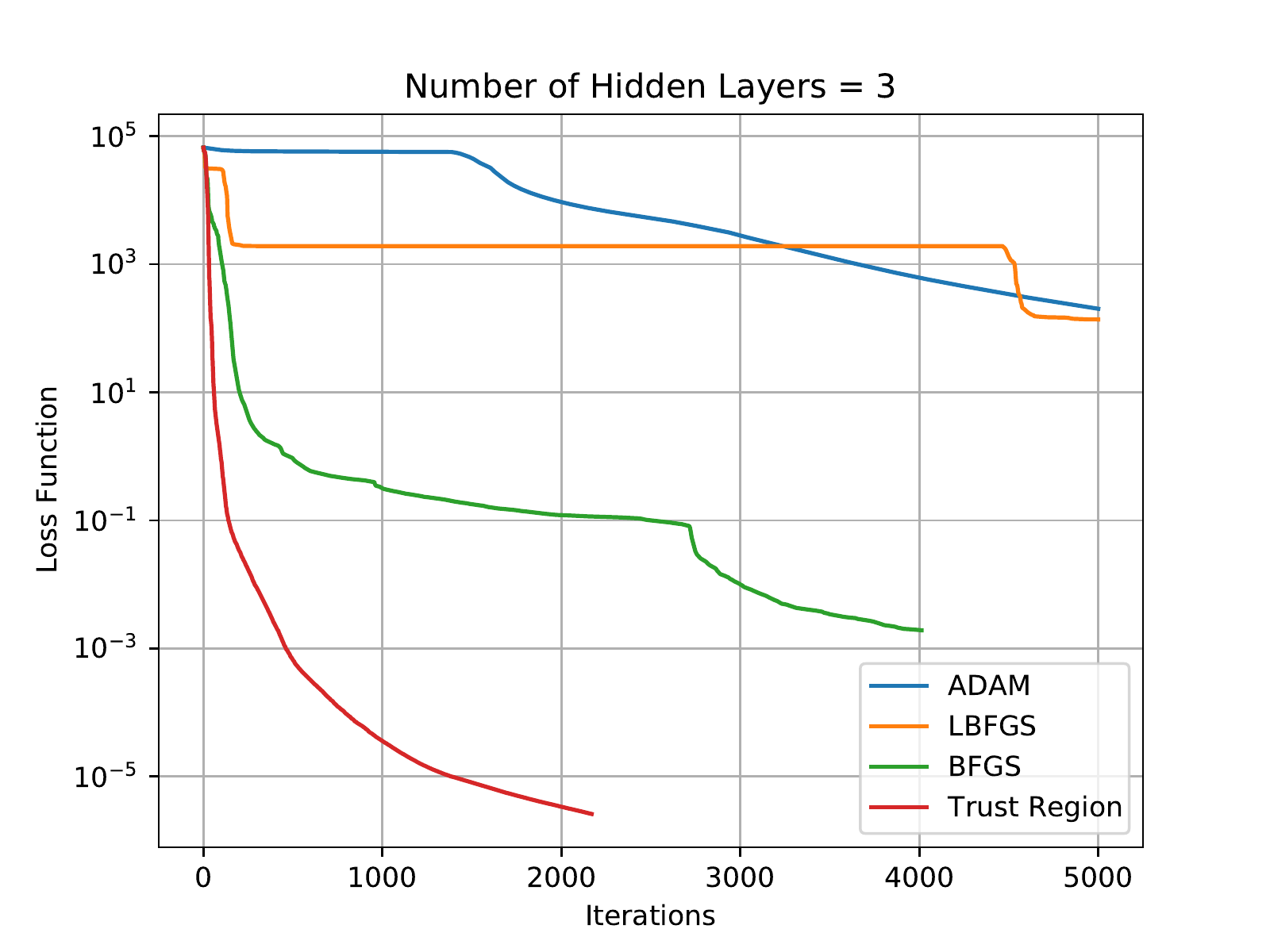}
  \caption{Loss functions for different number of layers. Trust region methods have the best performance, while ADAM barely converged.}
  \label{fig:second-order-depth}
\end{figure}

\Cref{fig:second-order-depth} shows the results for a different number of hidden layers (1, 2, and 3). Each layer has widths 20 and the activation functions are $\tanh$. We can see that in all cases, trust region methods have the best performance, followed by BFGS, LBFGS, and then ADAM. This order is consistent with previous findings. ADAM optimizer barely converged in all cases. Therefore, despite that we report the statistics of ADAM optimizers in the following for completeness, the results should be carefully interpreted. We report the ratios of zero eigenvalues using the following threshold: $|\lambda| < 10^{-6}\lambda_{\max}$ is treated as non-positive eigenvalues, where $\lambda_{\max}$ is the maximum eigenvalue. \Cref{tab:second-order-depth} shows the results. As we can see, when the number of hidden layers increases, the portion of zero eigenvalues increases. This indicates that with overparametrization through increasing depths of DNNs, the minimizer lies on a much higher dimensional manifold of the parameter space. This generally makes the optimizer easier\footnote{For example, in 3D, finding a minimizer on a hyper-plane (number of zero eigenvalues in Hessian is 2) is easier than finding a minimizer on a line (number of zero eigenvalues in Hessian is 2).}

\begin{table}[htpb]
\centering
\caption{Ratios (\%) of zero eigenvalues. The number of neurons per layer is 20 and the activation function is tanh.}
\label{tab:second-order-depth}
\begin{tabular}{@{}lllll@{}}
\toprule
Number of Hidden   Layers & ADAM  & LBFGS & BFGS  & Trust Region \\ \midrule
1                         & 0     & 76.54 & 72.84 & 77.78        \\
2                         & 69.46 & 98.2  & 94.41 & 93.21        \\
3                         & 85.99 & 98.7  & 98.15 & 96.09        \\ \bottomrule
\end{tabular}
\end{table}

\begin{figure}[htbp]
  \centering
  \includegraphics[width=0.32\textwidth]{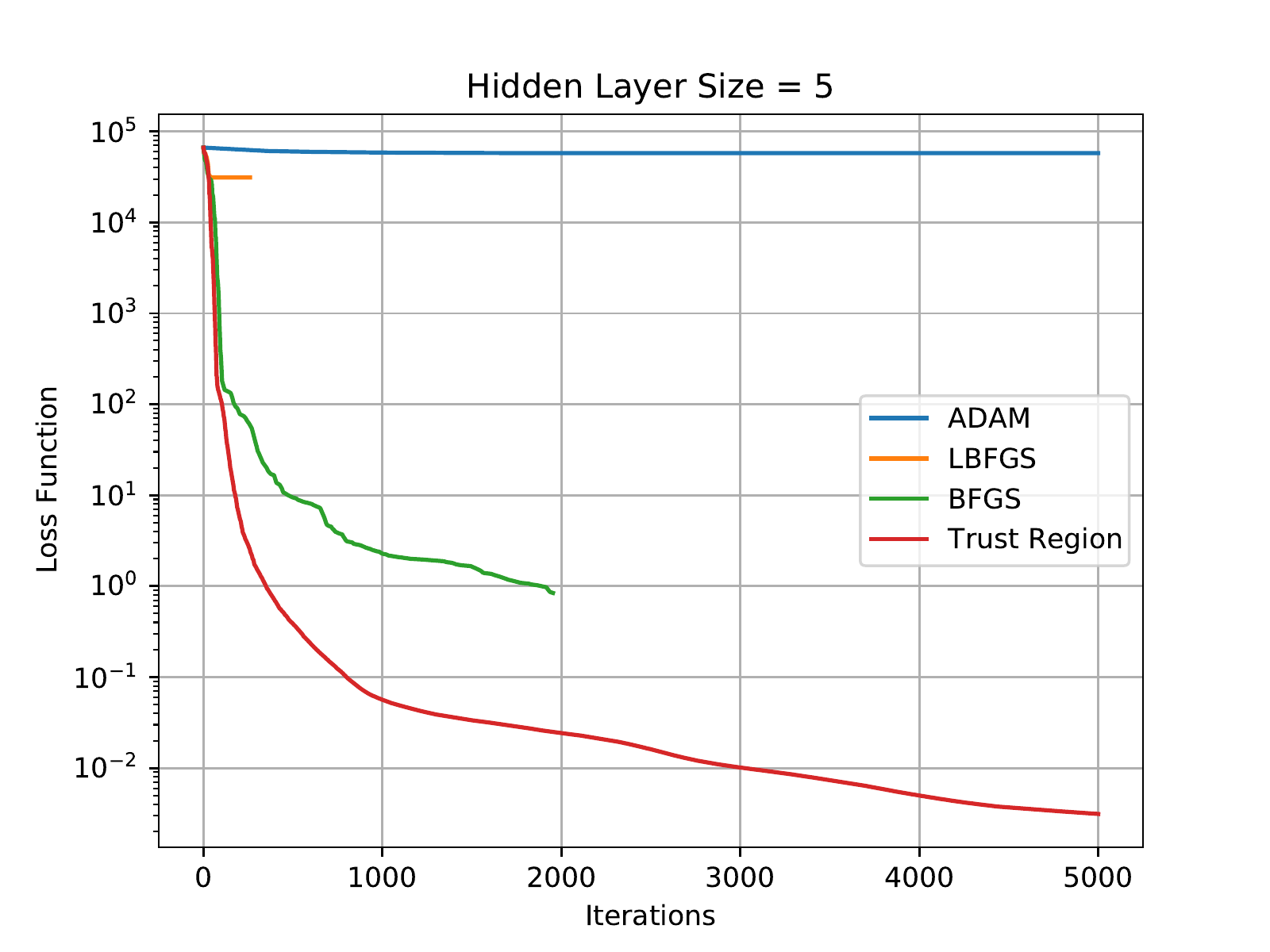}
  \includegraphics[width=0.32\textwidth]{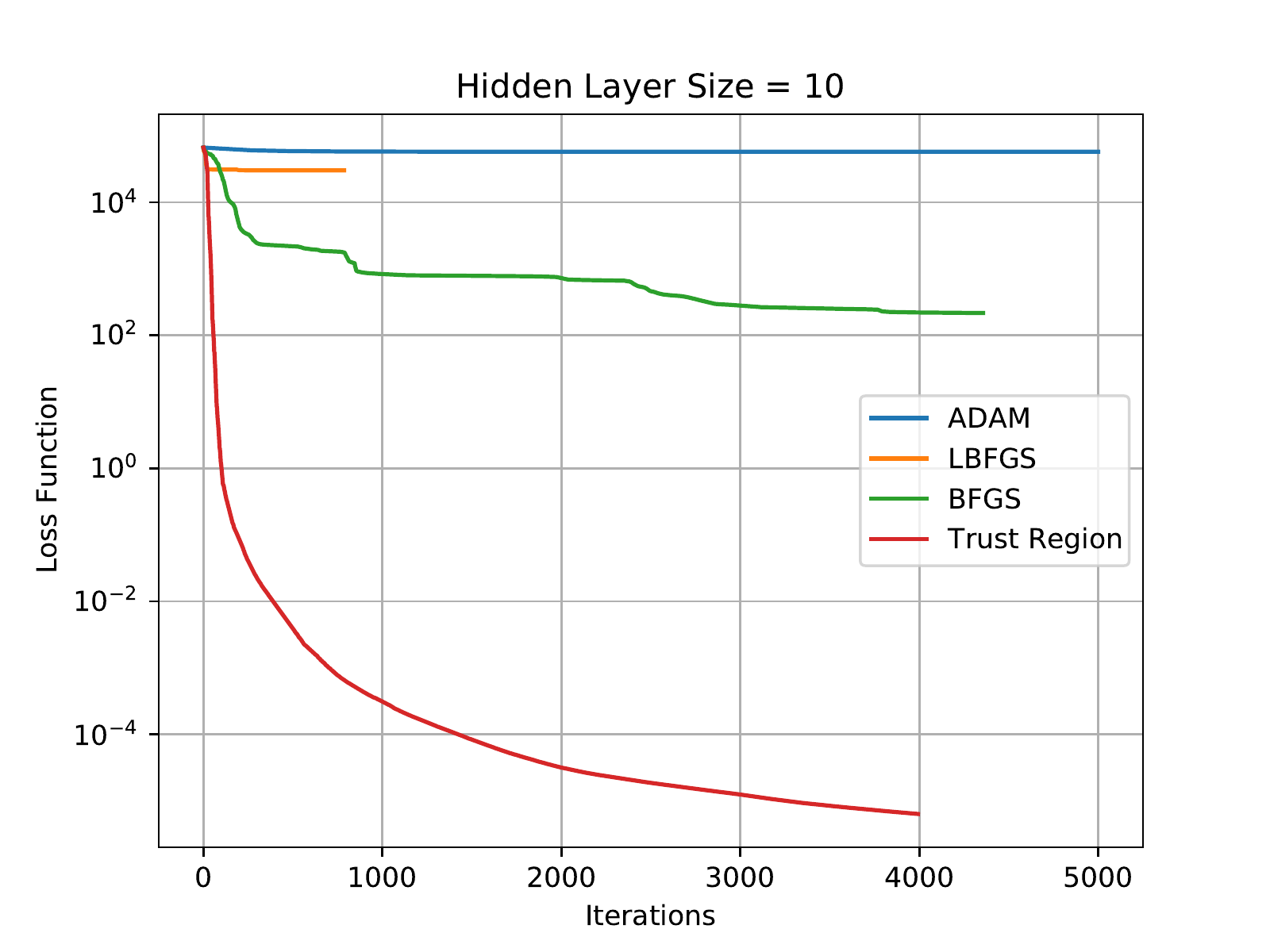}
  \includegraphics[width=0.32\textwidth]{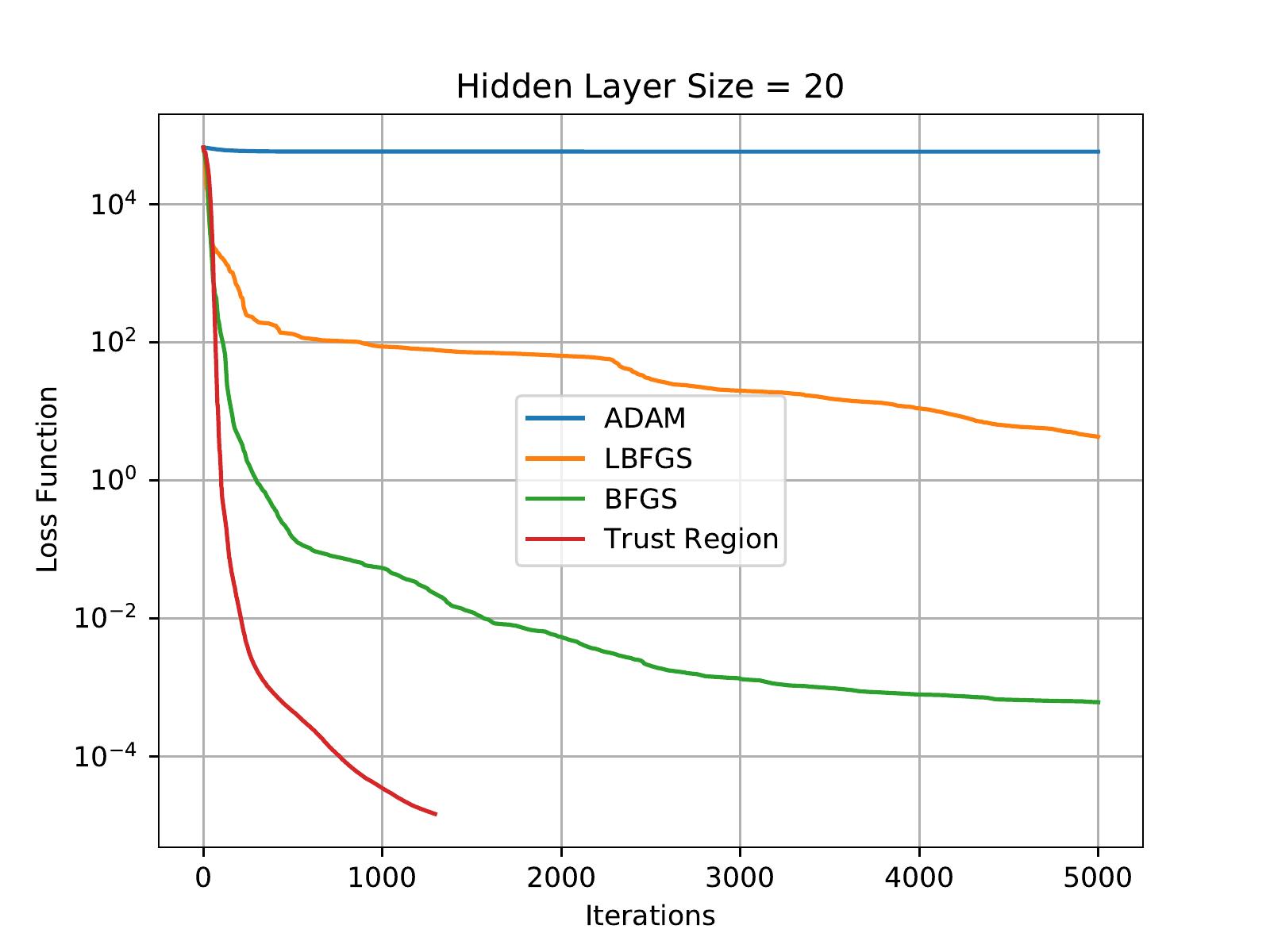}
  \caption{Loss functions for different number of layers. Trust region methods have the best performance, while ADAM barely converged.}
  \label{fig:second-order-width}
\end{figure}

\Cref{fig:second-order-width} shows the results for different number of neurons per layer (5, 10, and 20). The DNN has 3 layers and the activation functions are $\tanh$. \Cref{tab:second-order-width} shows the portion of zero eigenvalues. We see that by increasing widths, we can also increase the portion of zero eigenvalues consistently for trust region methods. This is not the case for LBFGS, where width 10 has the highest portion of zero eigenvalues. This indicates we can also gain the benefit of overparametrization via expanding widths.

\begin{table}[htpb]
\centering
\caption{Ratios (\%) of zero eigenvalues. The number of hidden layers is 3 and the activation function is tanh.}
\label{tab:second-order-width}
\begin{tabular}{@{}lllll@{}}
\toprule
Hidden Layer Size & ADAM  & LBFGS & BFGS  & Trust Region \\ \midrule
5                         & 24.69 & 93.83 & 85.19 & 69.14        \\
10                         & 50.19 & 97.7  & 83.52 & 89.66        \\
20                        & 76.98 & 96.2  & 97.39 & 96.42        \\ \bottomrule
\end{tabular}
\end{table}

These results reveal why overparametrization via increasing numbers of widths and depths makes optimization easier: \textbf{the minimizer lies on a relatively higher-dimensional manifold of the parameter space}. 

\subsection{Static Poisson's Equation: Finite Element Solver}

\begin{figure}[htbp]
  \centering
  \includegraphics[width=0.6\textwidth]{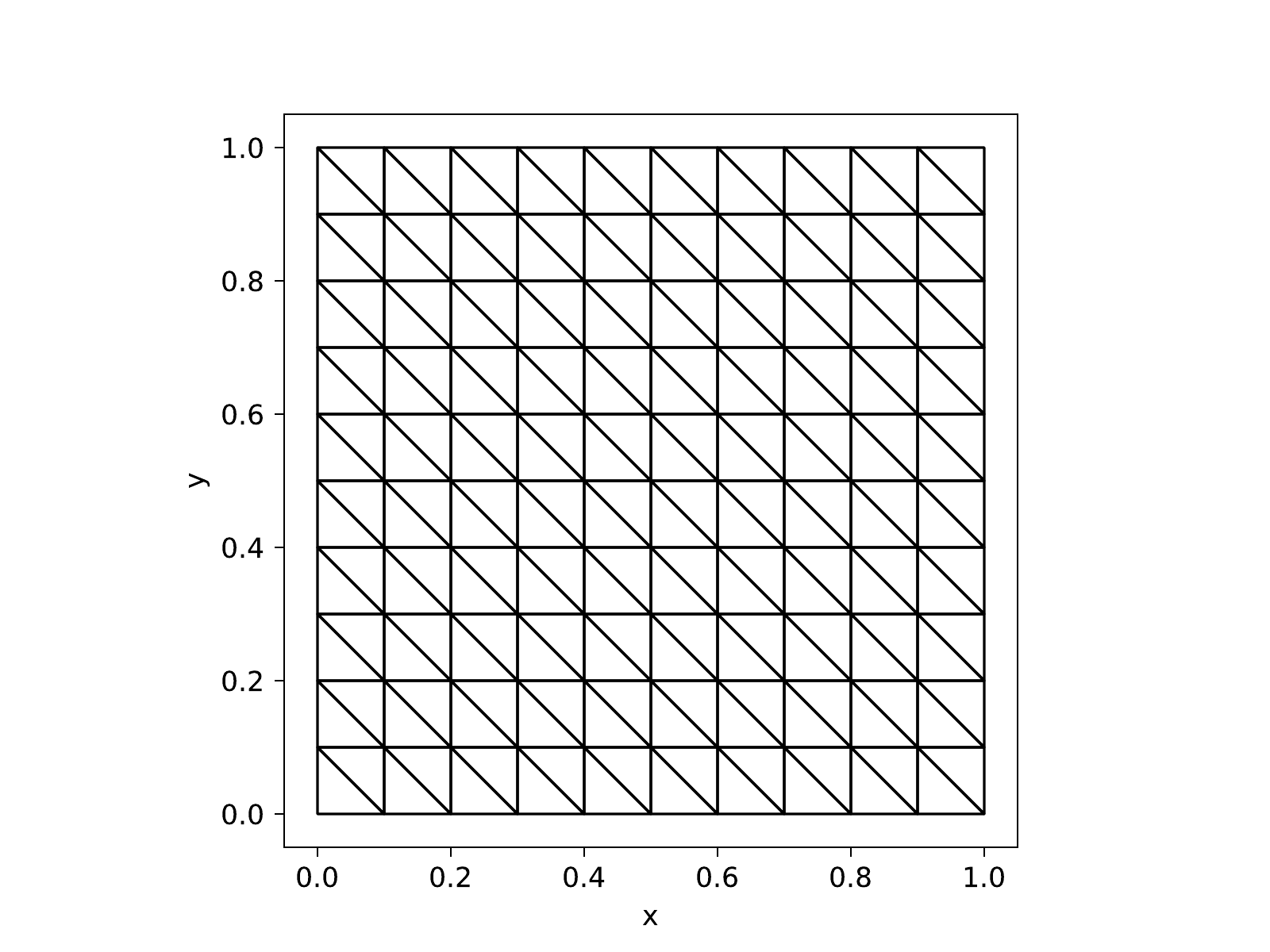}
  \caption{Finite element mesh used in \Cref{equ:second-order-optimization3}.}
  \label{fig:second-order-finite-element-mesh}
\end{figure}

In the last example, we consider the Poisson's equation with a spatially-varying diffusivity coefficient $\kappa(x,y)$
\begin{equation}\label{equ:second-order-ex3}
\begin{aligned}
     \nabla \cdot (\kappa(x,y) \nabla u)) &= f(x) & (x,y)\in \Omega\\ 
     u &= 0 & (x,y)\in  \partial\Omega
\end{aligned}
\end{equation}
Here 
$$\kappa(x, y) = \frac{1}{1+x^2+y^2}+1$$
We use the linear finite element method to solve \Cref{equ:second-order-ex3}, and the finite element mesh is shown in \Cref{fig:second-order-finite-element-mesh}. 

\begin{figure}[htbp]
  \centering
  \includegraphics[width=0.2\textwidth]{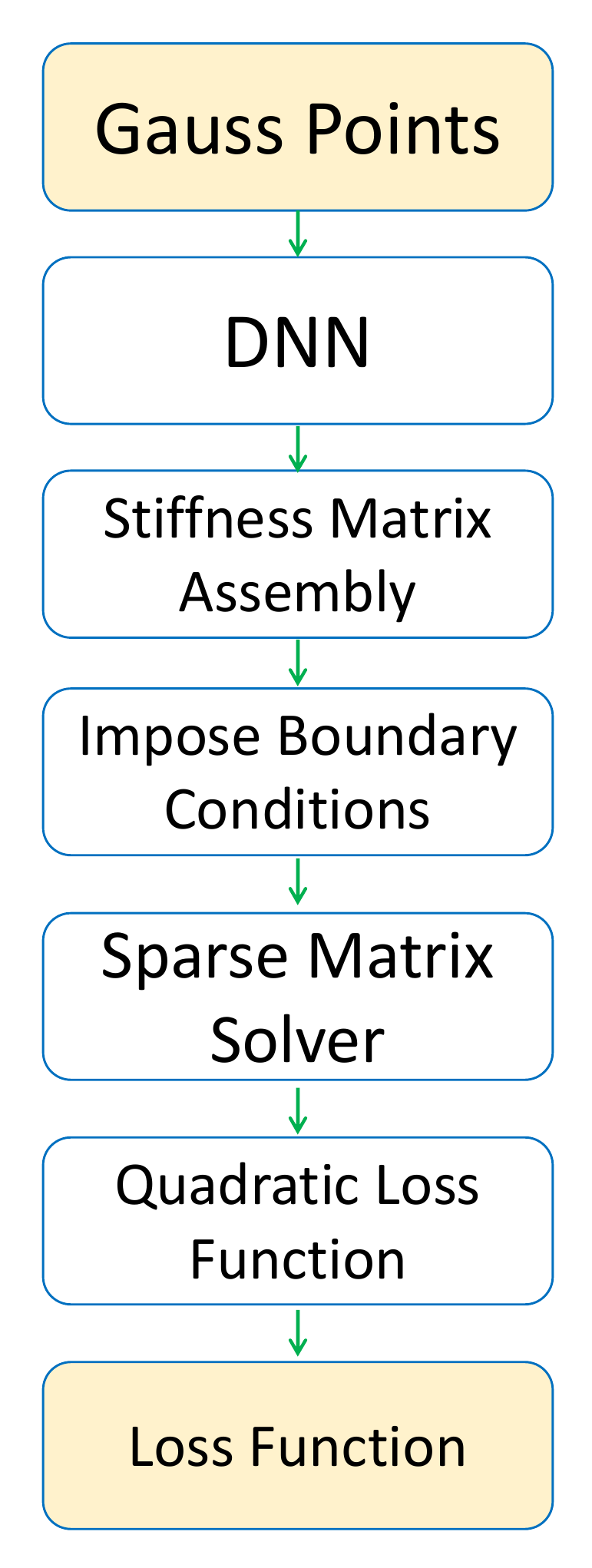}
  \caption{Workflow for implementing \Cref{equ:second-order-optimization3}. This workflow involves a matrix solver, for which we need to use techniques in \Cref{sect:sparse_linear_solver} to derive Hessian update rules.}
  \label{fig:second-order-impl}
\end{figure}

We assume that $u$ can be observed at all finite element nodes $\{(x_i,y_i)\}_{i=1}^n$, and the observation is denoted as $u_i$. We want to use this information to estimate $\kappa(x,y)$. We approximate $\kappa(x,y)$ using a deep neural network $\kappa_\theta(x,y)$ ($\theta$ is the weights and biases) and discretize \Cref{equ:second-order-ex3} after replacing $\kappa$ with $\kappa_\theta$. We formulate the inverse problem as an optimization problem
\begin{equation}\label{equ:second-order-optimization3}
    \begin{aligned}
    \min_{\theta} \; \sum_{i=1}^n (\mathbf{u}_i - u_i)^2\\ 
    \text{s.t.} \; A(\theta)\mathbf{u} = \mathbf{f} 
    \end{aligned}
\end{equation}
Here $A(\theta)\mathbf{u} = \mathbf{f}$ is the finite element discretization for \Cref{equ:second-order-ex3} and $\mathbf{u}$ is the discretized solution vector. $A(\theta)$ is the stiffness matrix, which is sparse and whose entries depend on $\theta$. Different from the residual minimization approach in the last two sections, \Cref{equ:second-order-optimization3} minimizes the discrepancy between the state variable and the observation directly. The advantage is that \Cref{equ:second-order-optimization3} can easily deal with sparse observations (i.e., only a subset of $u_i$ can be observed). However, the implementation of second order PCL for \Cref{equ:second-order-optimization3} is much more challenging as it involves a matrix solver, where we need to build a Hessian update rule introduced in \Cref{sect:sparse_linear_solver} (see \Cref{fig:second-order-impl}).

\begin{figure}[htbp]
  \centering
  \includegraphics[width=0.33\textwidth]{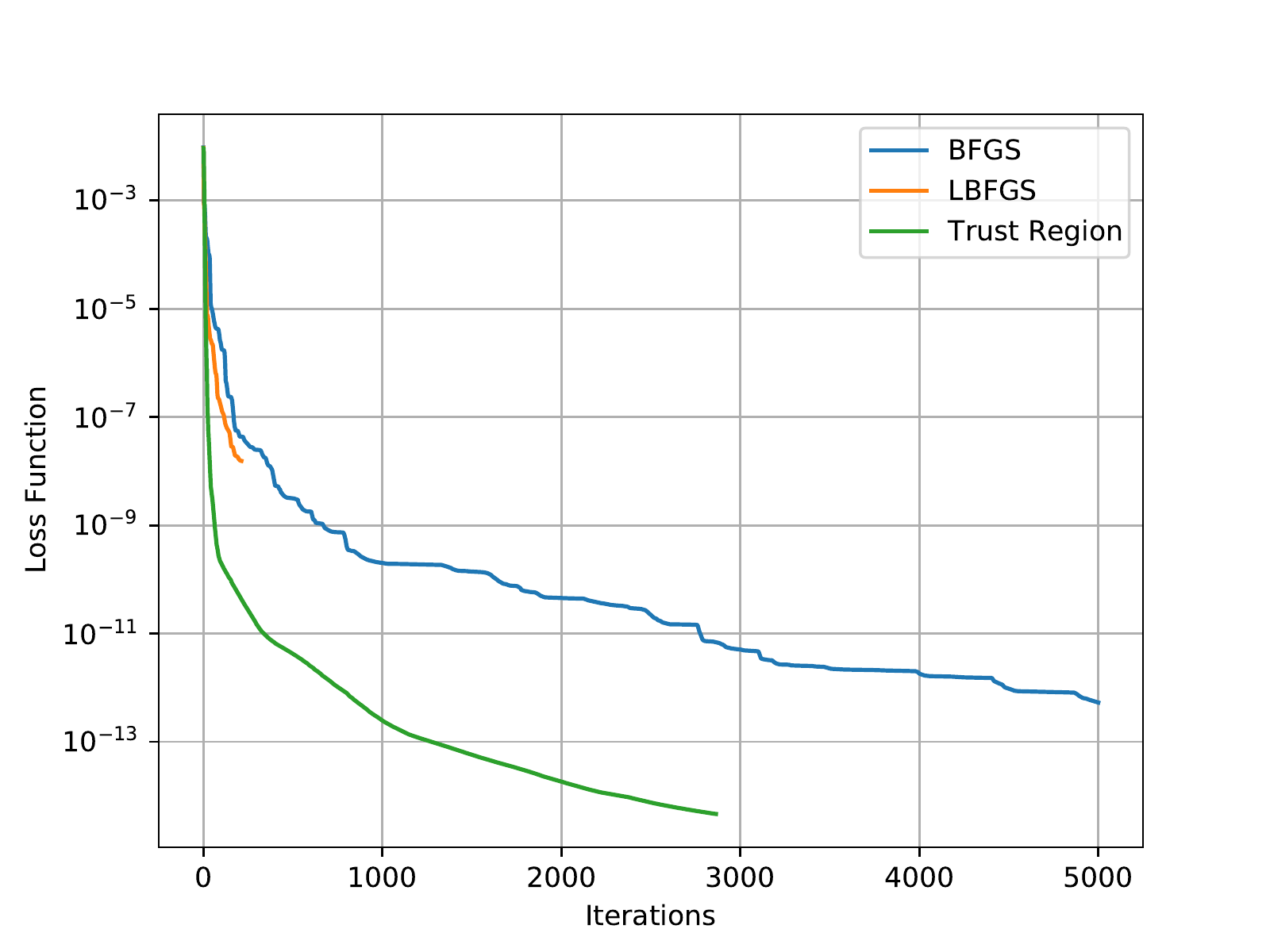}~
  \includegraphics[width=0.33\textwidth]{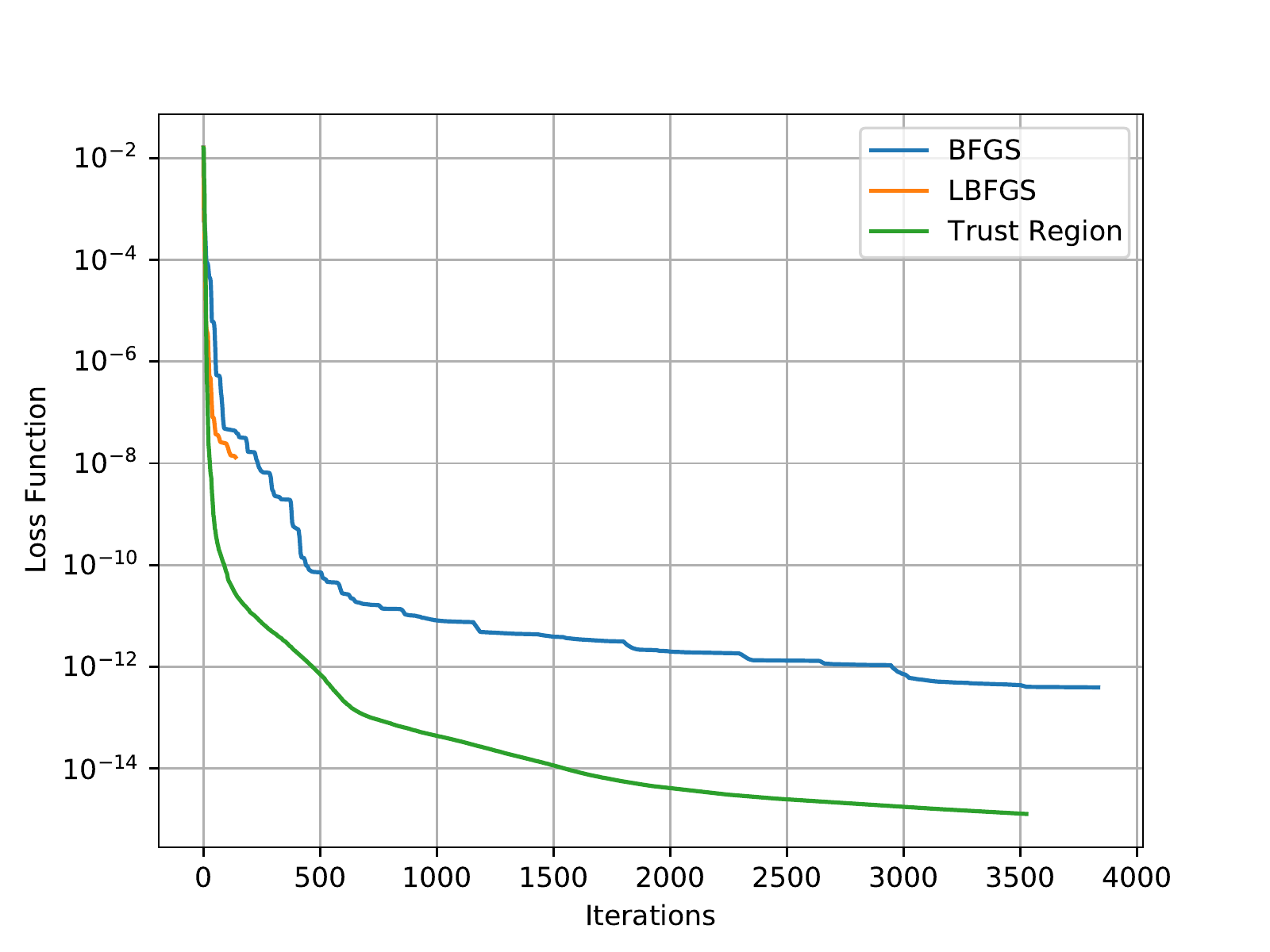}~
  \includegraphics[width=0.33\textwidth]{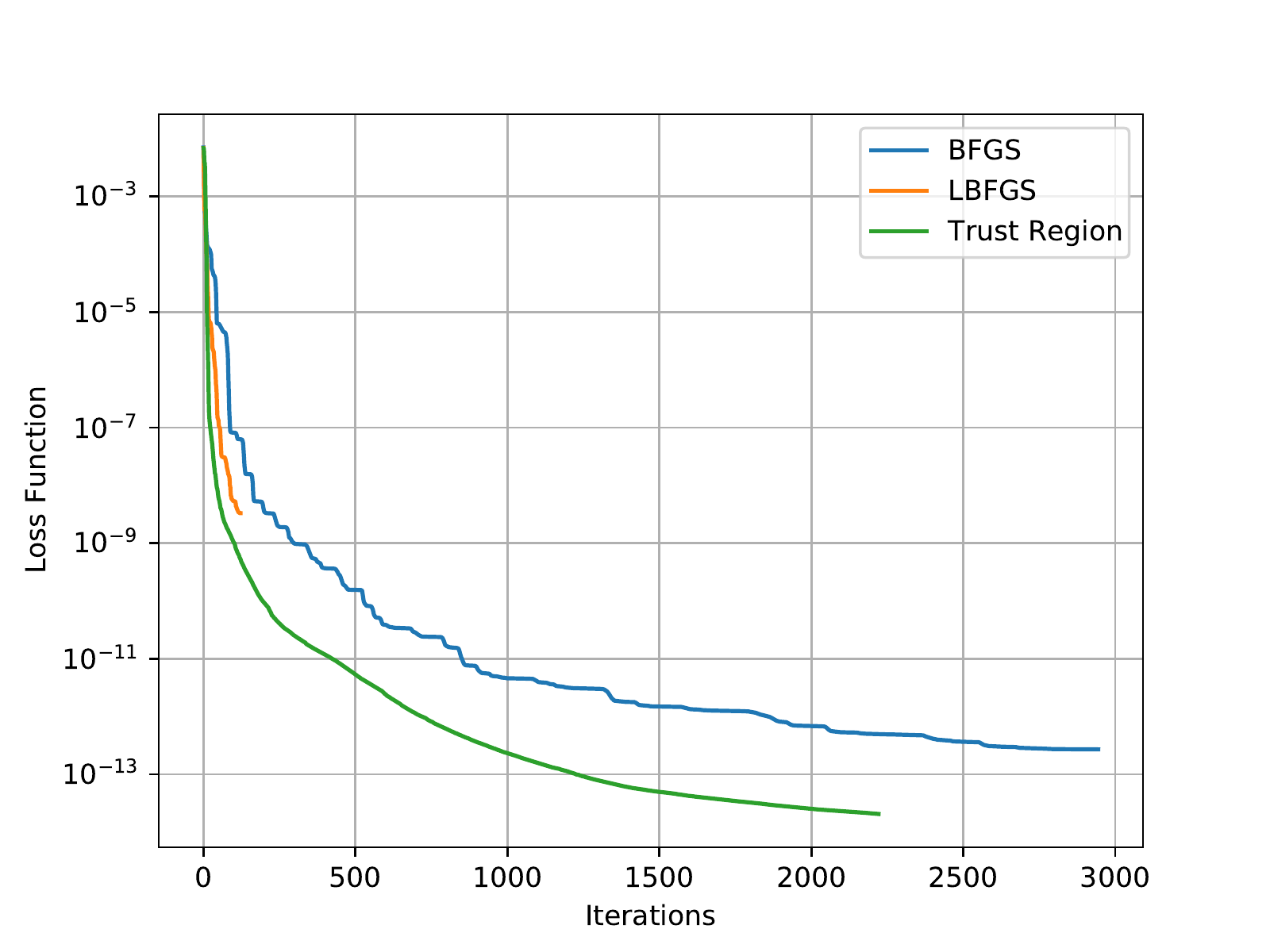}
  \caption{Loss functions for \Cref{equ:second-order-optimization3}. Trust region methods outperform all other methods.}
  \label{fig:second-order-ex3-loss}
\end{figure}

\begin{figure}[htbp]
  \centering
  \includegraphics[width=0.8\textwidth]{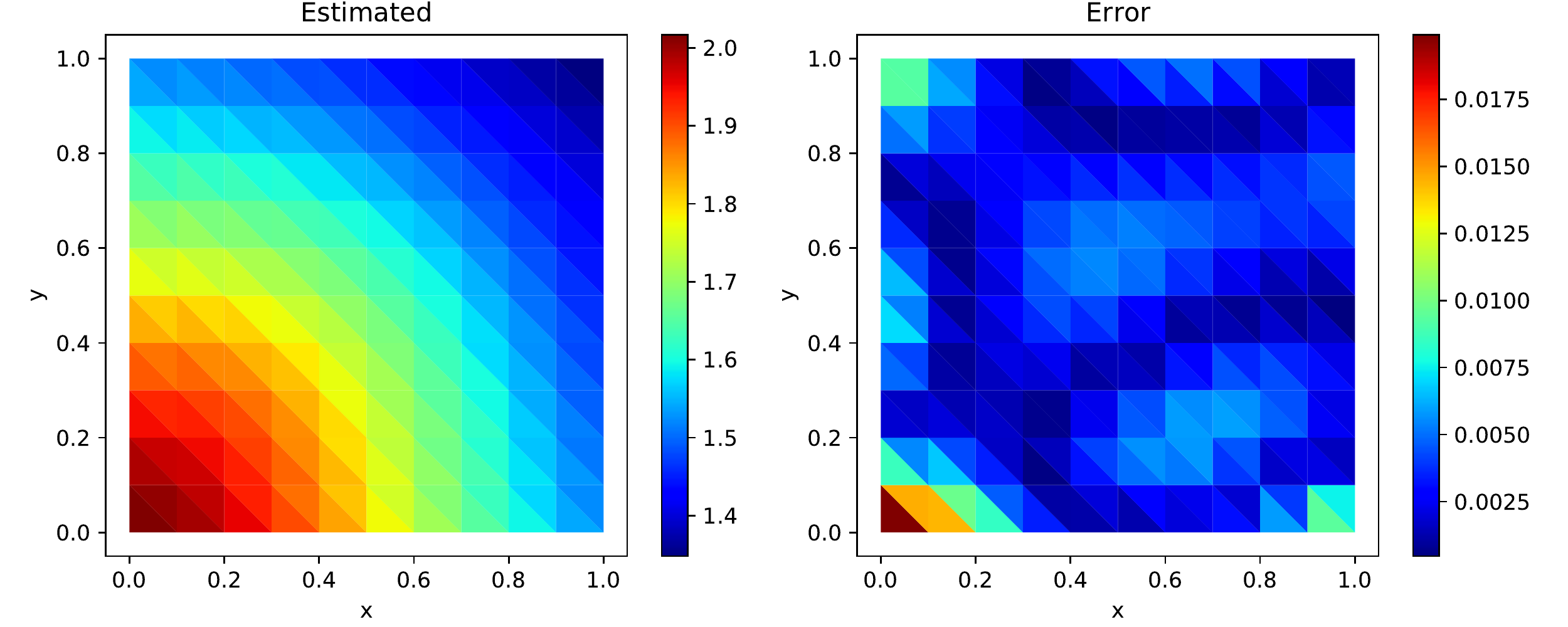}
  \includegraphics[width=0.8\textwidth]{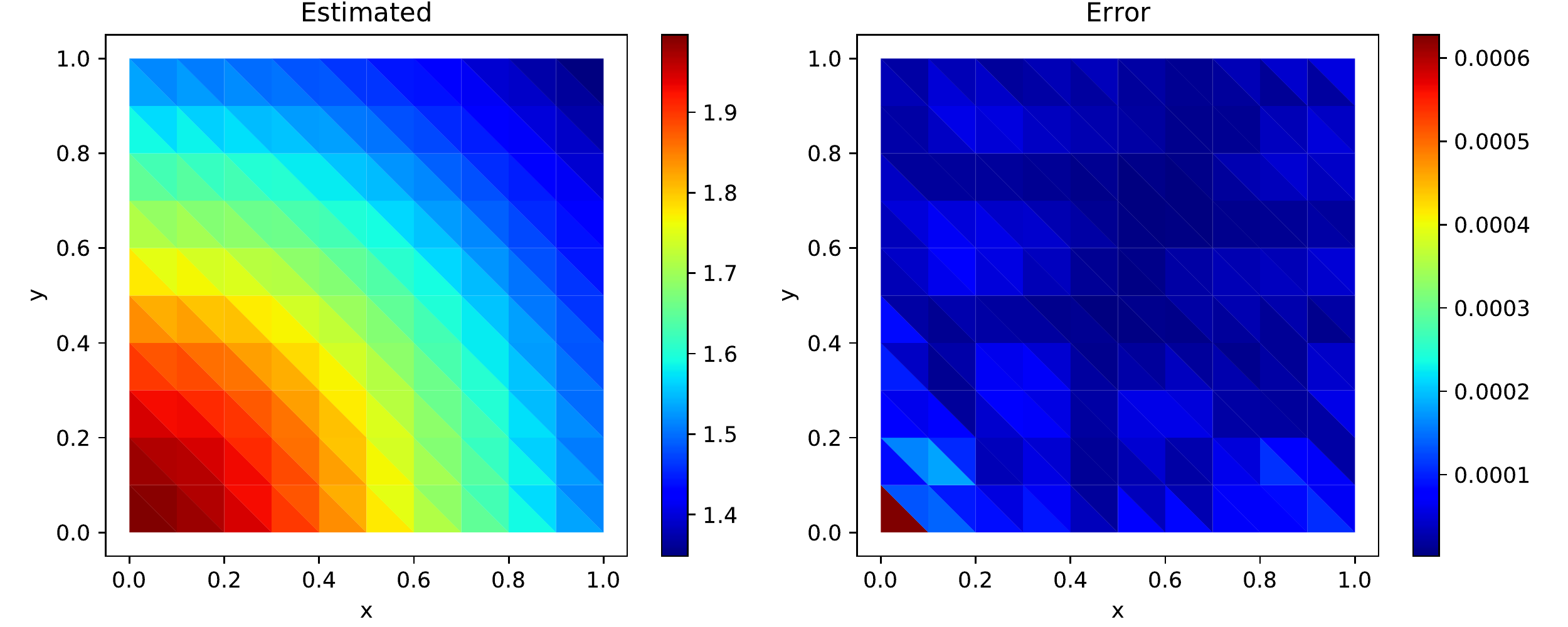}
  \includegraphics[width=0.8\textwidth]{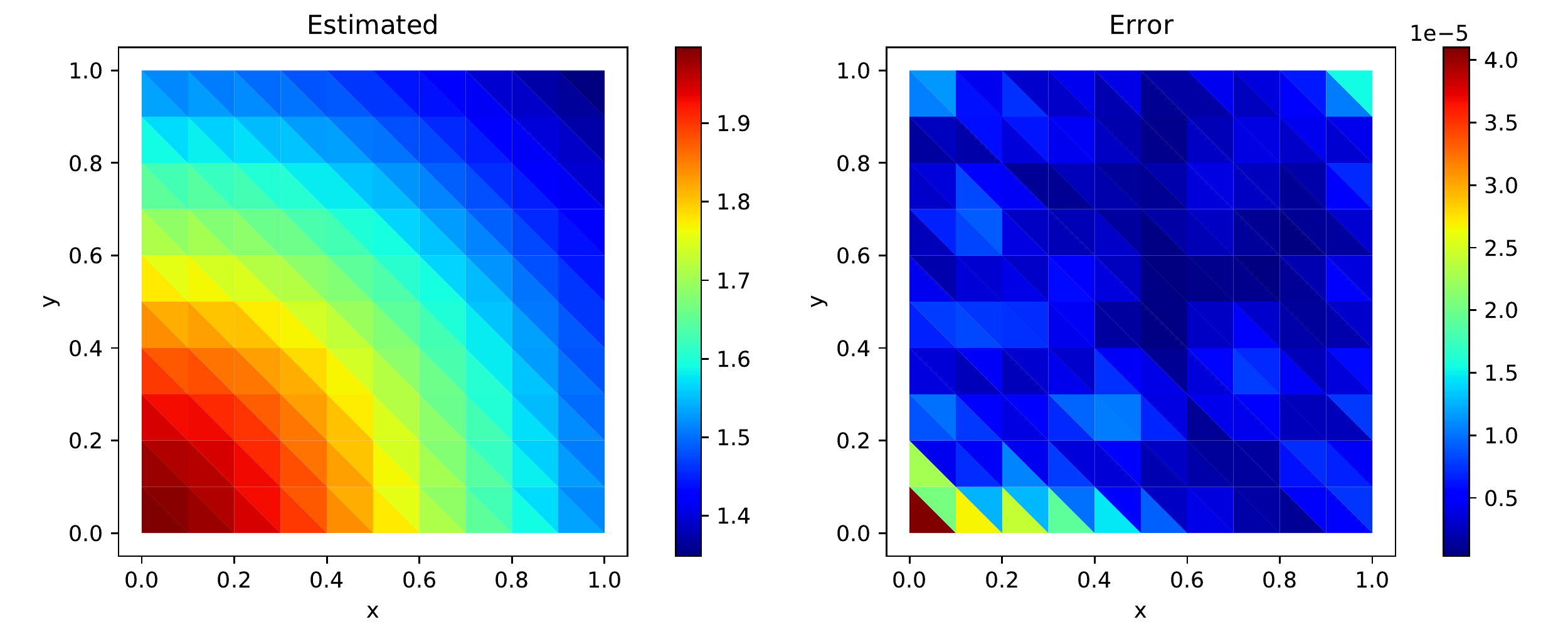}
  \caption{Estimated $\kappa_\theta$ and pointwise errors $|\kappa-\kappa_\theta|$. From top to bottom: LBFGS, BFGS, and trust region.}
  \label{fig:second-order-compare}
\end{figure}

\Cref{fig:second-order-ex3-loss} shows the loss function profiles for three different initialized guesses. We see that trust region methods outperform all other methods. Specifically, \Cref{fig:second-order-compare} shows the estimated $\kappa_\theta$ and the pointwise error $|\kappa-\kappa_\theta|$ for the first case in \Cref{fig:second-order-ex3-loss}, which shows that our estimation is quite accurate. 

\begin{figure}[htbp]
  \centering
  \includegraphics[width=0.33\textwidth]{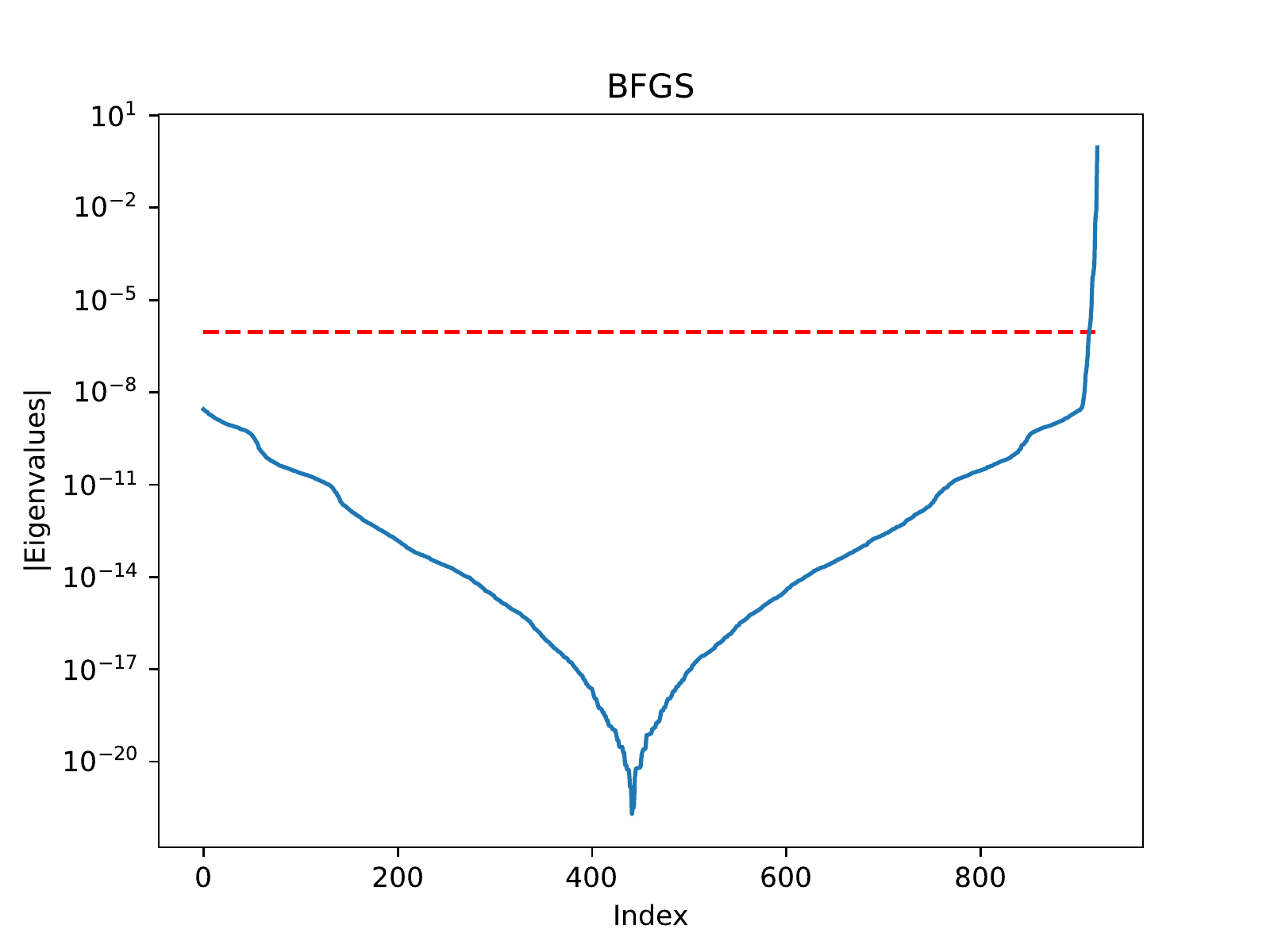}~
  \includegraphics[width=0.33\textwidth]{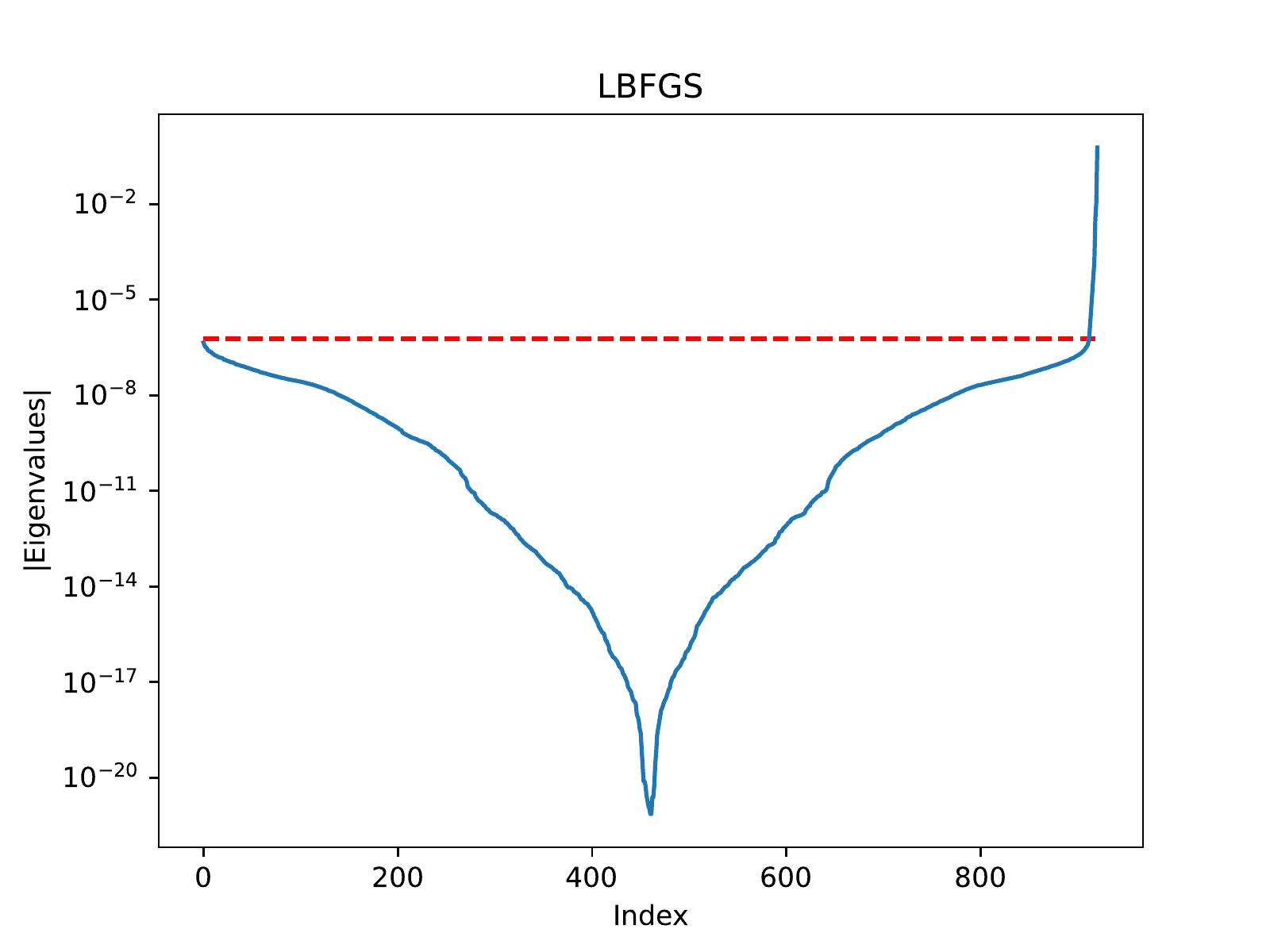}~
  \includegraphics[width=0.33\textwidth]{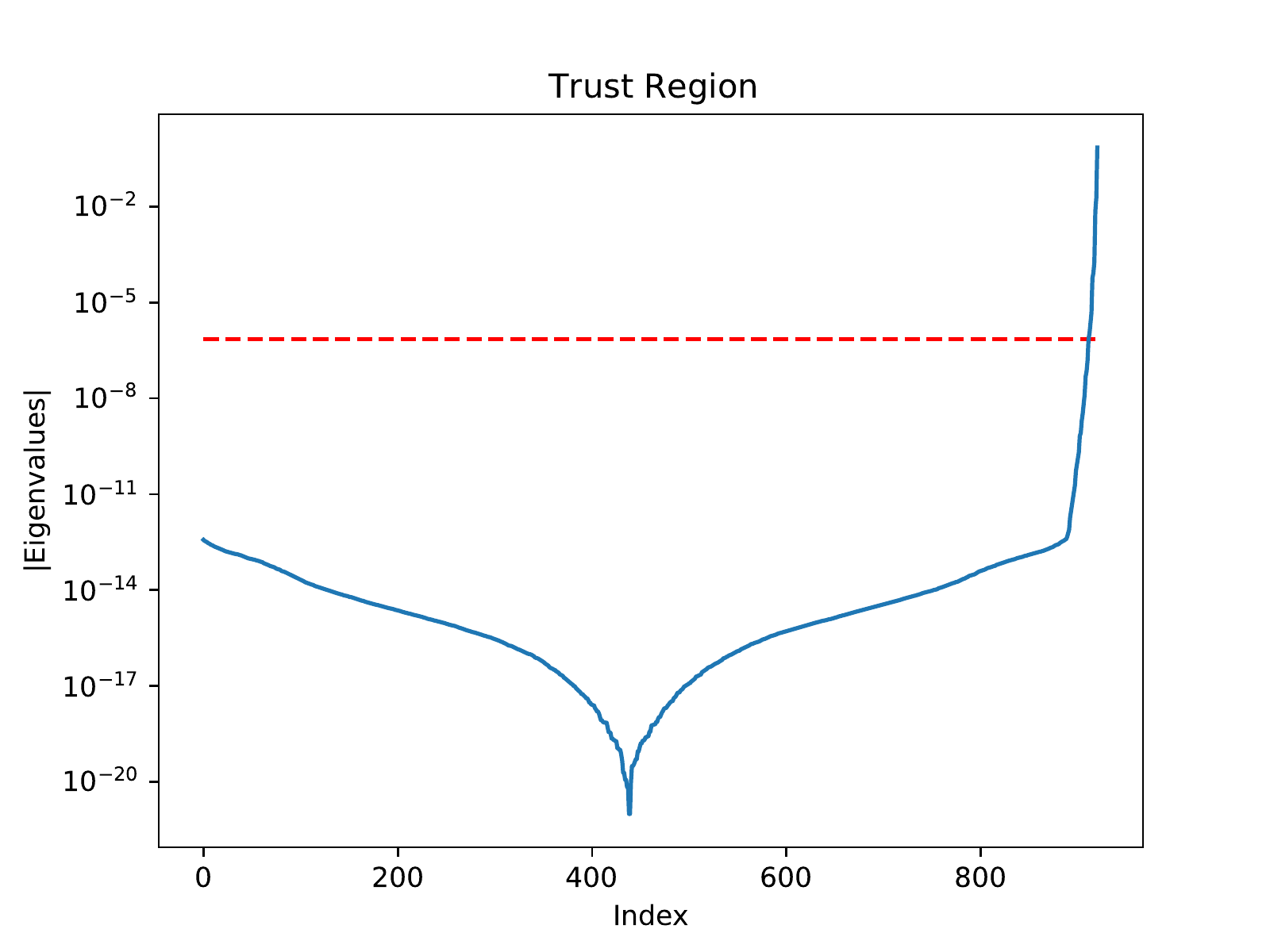}
  \caption{The magnitudes of eigenvalues for different optimizers. For all cases, only a few eigenvalues are nonzero.}
  \label{fig:second-order-ex3-pcl-eigs}
\end{figure}

\begin{figure}[htbp]
  \centering
  \includegraphics[width=0.45\textwidth]{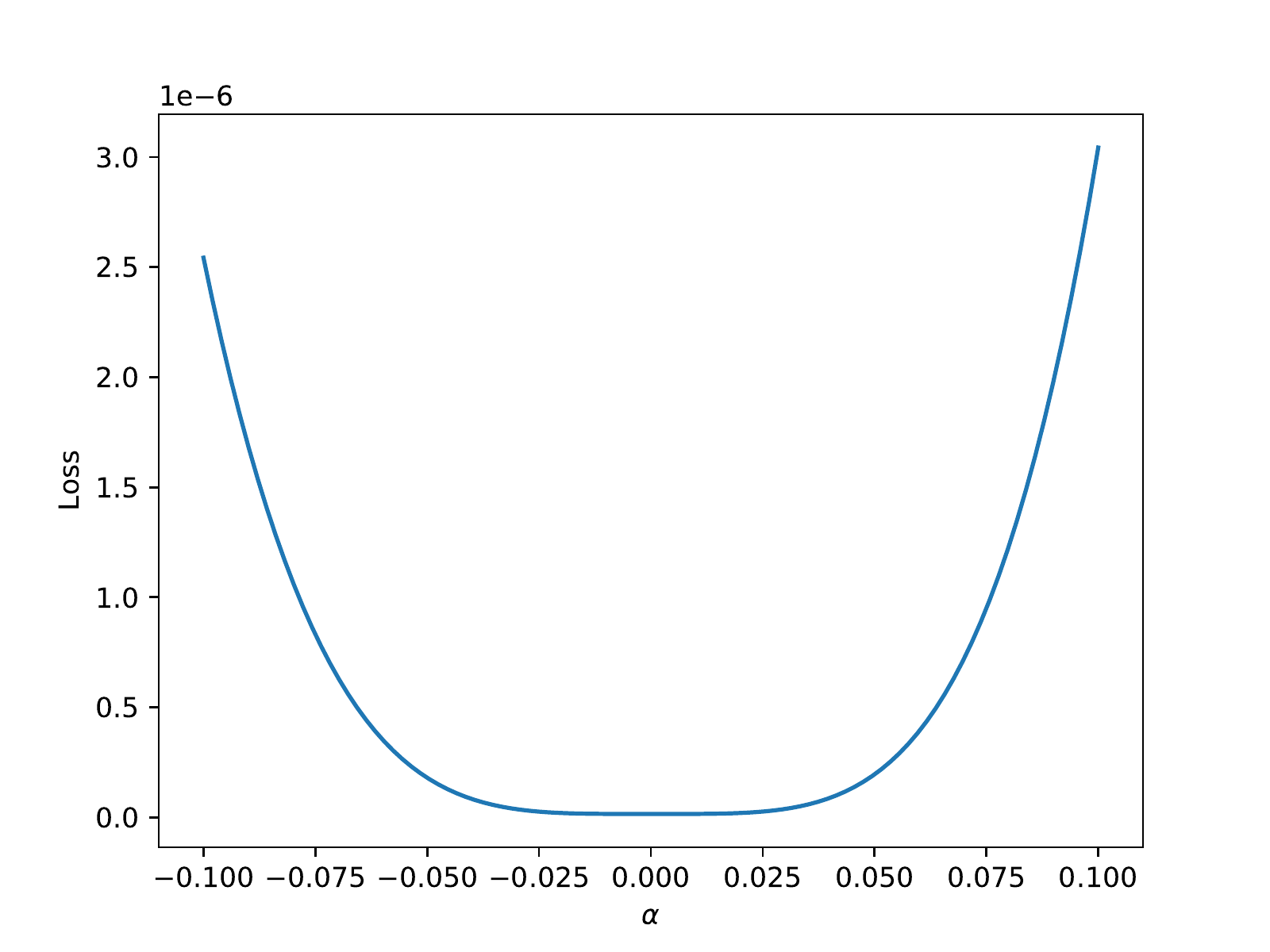}~
  \includegraphics[width=0.45\textwidth]{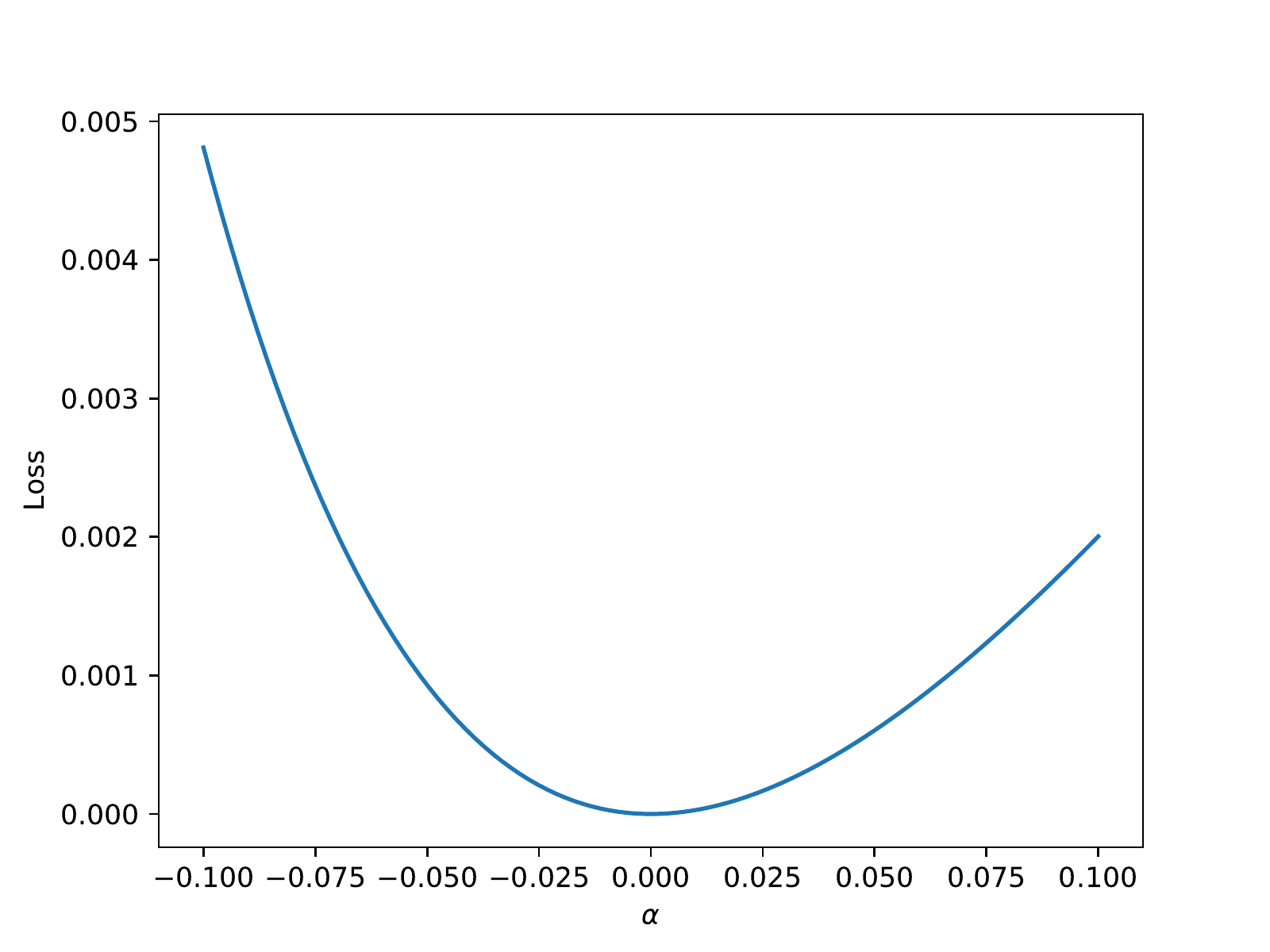}
  \caption{Loss function profiles for $F(\alpha) = L(x^* + \alpha v)$. Left: $v$ is the eigenvector associated with $\lambda_{\min}$; we can see that the local landscape at $\alpha = 0$ is almost flat. Right: $v$ is the eigenvector associated with $\lambda_{\max}$.}
  \label{fig:second-order-ex3-pcl-evs}
\end{figure}

\Cref{fig:second-order-ex3-pcl-eigs} shows the magnitudes of eigenvalues for different optimizers. These plots confirm our findings in the last two sections again: at convergence, only a small subset of eigenvalues are positive, and all other eigenvalues are nearly zero. Eigenvalues that lie below the red dashed line can be treated as zero. This means that for BFGS and the trust region method, the optimizers find local minima. We show
$$F(\alpha) = L(x^* + \alpha v)$$
in \Cref{fig:second-order-ex3-pcl-evs}, where $x^*$ is the converged point for LBFGS, $v$ is the eigenvector corresponding to either the minimum or maximum eigenvalues of the Hessian. The profile for the former case is quite flat, indicating that small perturbation along the eigenvector direction makes little change to the loss function. 

\section{Limitations}
Despite many promising features of the trust region method, it is not without limitations, which we want to discuss here. The current trust region method requires calculating the Hessian matrix. Firstly, computing the Hessian matrix can be technically difficult, especially when DNNs are coupled with a sophisticated numerical PDE solver. There are many existing techniques for computing the Hessian. The TensorFlow backend supports Hessian computation concurrently, but it requires users to implement rules for calculating ``gradients of gradients''. Additionally, TensorFlow uses reverse-mode automatic differentiation to evaluate the Hessian. This means that TensorFlow loops over each gradient component and calculating a row of Hessian at a time. This does not leverage the symmetry of Hessians and can be quite inefficient if the number of unknowns is large. Another approach, the edge pushing algorithm, which we use in this work, uses one backward pass to evaluate the Hessian. This approach takes advantage of the symmetry of Hessians. However, a general purpose implementation can be quite convoluted and computations can be expensive in some scenarios. Instead, we develop specialized algorithms for coupled systems of DNNs and PDEs, which enables us to leverage problem structures. An interesting direction in the future is to reduce the computational cost but still captures all the benefits brought about by trust region methods; for example, subsampling from the Hessian matrices is a promising direction.

\section{Conclusion}

Trust region methods are a class of global optimization techniques. They are less popular in the deep learning approach because the DNNs tend to be huge and the computation of Hessians is expensive. However, they are very suitable for many computational engineering problems, where DNNs are typically small, and convergence, as well as accuracy, is a critical concern. The problems themselves are nonconvex and have many local minima ---different from the common belief that in deep learning, stationary points are usually saddle points if they are not a global minimum. Trust region methods do not guarantee that we can find a global minimum, or even a ``good'' local minimum. However, because trust region methods show faster convergence and superior accuracy in many cases, it never harms to add trust region methods into the optimization toolbox. Additionally, the Hessian calculated using the second order PCL is a powerful weapon for diagnosing the convergence and provides curvature information for more sophisticated optimizers.
Our point of view is that although the Hessian computations are expensive, they are quite rewarding. Future researches will focus on efficient computation and automation of Hessian computations.

\bibliographystyle{unsrt}
\bibliography{ref.bib}

\end{document}